\pdfoutput=1
\RequirePackage{ifpdf}
\ifpdf 
\documentclass[pdftex]{sigma}
\else
\documentclass{sigma}
\fi

\newtheorem{thm}{Theorem}[section]
\newtheorem{prop}[thm]{Proposition}

\newtheorem{lem}[thm]{Lemma}
\theoremstyle{definition}
\newtheorem{defn}[thm]{Definition}
\newtheorem{rem}[thm]{Remark}
\newtheorem{rems}[thm]{Remarks}
\newtheorem{aside}[thm]{Aside}

\numberwithin{equation}{section}

\DeclareMathOperator\iden{id}
\DeclareMathOperator\vol{vol}
\DeclareMathOperator\kernel{ker}

\DeclareMathOperator\trace{tr}

\DeclareMathOperator\divg{div}

\DeclareMathOperator\rank{rank}

\DeclareMathOperator\End{End}
\DeclareMathOperator\Hom{Hom}
\DeclareMathOperator\Sym{Sym}

\DeclareMathOperator\scal{scal}
\DeclareMathOperator\GL{GL}
\DeclareMathOperator\SL{SL}

\DeclareMathOperator\SO{SO}

\DeclareMathOperator\Un{U}
\DeclareMathOperator\SU{SU}
\DeclareMathOperator\Diff{Dif\/f}
\DeclareMathOperator\SDiff{SDif\/f}

\DeclareMathOperator\Vect{Vect}

\newcommand{\abrack}[1]{[\mkern-3mu[#1]\mkern-3mu]}
\newcommand{\Abrack}[1]{\bigl[\mkern-5mu\bigl[#1\bigr]\mkern-5mu\bigr]}
\newcommand{\ip}[1]{\langle#1\rangle}
\newcommand{\Ip}[1]{\bigl\langle#1\bigr\rangle}
\newcommand{\lie}[1]{\mathfrak{#1}}

\newcommand{\CP}[1]{{\mathbb{C}}{\rm P}^{#1}}
\newcommand{\sd}{{\raise1pt\hbox{$\scriptscriptstyle +$}}}
\newcommand{\asd}{{\raise1pt\hbox{$\scriptscriptstyle -$}}}
\newcommand{\sdasd}{{\raise1pt\hbox{$\scriptscriptstyle\pm$}}}
\newcommand{\asdsd}{{\raise1pt\hbox{$\scriptscriptstyle\mp$}}}

\begin{document}

\allowdisplaybreaks

\renewcommand{\thefootnote}{$\star$}

\renewcommand{\PaperNumber}{034}

\FirstPageHeading

\ShortArticleName{Integrable Background Geometries}

\ArticleName{Integrable Background Geometries\footnote{This paper is a~contribution to the Special Issue on Progress in
Twistor Theory. The full collection is available at
\href{http://www.emis.de/journals/SIGMA/twistors.html}{http://www.emis.de/journals/SIGMA/twistors.html}}}

\Author{David M.J.~CALDERBANK}

\AuthorNameForHeading{D.M.J.~Calderbank}

\Address{Department of Mathematical Sciences, University of Bath, Bath BA2 7AY, UK}

\Email{\href{mailto:D.M.J.Calderbank@bath.ac.uk}{D.M.J.Calderbank@bath.ac.uk}}
\URLaddress{\url{http://people.bath.ac.uk/dmjc20/}}

\ArticleDates{Received January 21, 2014, in f\/inal form March 18, 2014; Published online March 28, 2014}

\Abstract{This work has its origins in an attempt to describe systematically the integrable geometries and gauge
theories in dimensions one to four related to twistor theory.
In each such dimension, there is a~nondegenerate integrable geometric structure, governed by a~nonlinear integrable
dif\/ferential equation, and each solution of this equation determines a~background geometry on which, for any Lie group
$G$, an integrable gauge theory is def\/ined.
In four dimensions, the geometry is selfdual conformal geometry and the gauge theory is selfdual Yang--Mills theory,
while the lower-dimensional structures are nondegenerate (i.e., non-null) reductions of this.
Any solution of the gauge theory on a~$k$-dimensional geometry, such that the gauge group $H$ acts transitively on an
$\ell$-manifold, determines a~$(k+\ell)$-dimensional geometry ($k+\ell\leqslant4$) f\/ibering over the $k$-dimensional
geometry with $H$ as a~structure group.
In the case of an $\ell$-dimensional group $H$ acting on itself by the regular representation, all
$(k+\ell)$-dimensional geometries with symmetry group $H$ are locally obtained in this way.
This framework unif\/ies and extends known results about dimensional reductions of selfdual conformal geometry and the
selfdual Yang--Mills equation, and provides a~rich supply of constructive methods.
In one dimension, generalized Nahm equations provide a~uniform description of four pole isomonodromic deformation
problems, and may be related to the $\SU(\infty)$ Toda and dKP equations via a~hodograph transformation.
In two dimensions, the $\Diff(S^1)$ Hitchin equation is shown to be equivalent to the hyperCR Einstein--Weyl equation,
while the $\SDiff(\Sigma^2)$ Hitchin equation leads to a~Euclidean analogue of Plebanski's heavenly equations.
In three and four dimensions, the constructions of this paper help to organize the huge range of examples of
Einstein--Weyl and selfdual spaces in the literature, as well as providing some new ones.
The nondegenerate reductions have a~long ancestry.
More recently, degenerate or null reductions have attracted increased interest.
Two of these reductions and their gauge theories (arguably, the two most signif\/icant) are also described.}


\Classification{53A30; 32L25; 37K25; 37K65; 53B35; 53C25; 58J70; 70S15; 83C20; 83C80}

\renewcommand{\thefootnote}{\arabic{footnote}} \setcounter{footnote}{0}

\section{Introduction}

\subsection*{Background}

This paper concerns dif\/ferential geometry in dimensions one to four, and (primarily) four kinds of geometric structure,
one in each dimension, governed by four nonlinear integrable dif\/ferential equations.
Associated to each manifold carrying one of these geometric structures, and to each Lie group, is an integrable gauge
theory, generalizing a~well-known gauge theory on f\/lat space: the selfdual Yang--Mills equation on ${\mathbb{R}}^4$ or
${\mathbb{R}}^{2,2}$ and the gauge f\/ield equations arising as reductions by a~(nondegenerate) group of translations to
lower-dimensional f\/lat spaces.
Specif\/ically, reduction by a~single non-null translation gives the Bogomolny equation for monopoles on ${\mathbb{R}}^3$,
or its Lorentzian analogue on ${\mathbb{R}}^{2,1}$~\cite{AtHi:gmm,Ward:iss}.
Reduction by two such translations yields the Hitchin equation for Higgs pairs on a~Riemann surface, harmonic maps from
a~Riemann surface to a~Lie group, or the principal chiral model on a~two-dimensional space-time~\cite{Hit:sde,Ward:iss}.
Reduction by three non-null translations leads to the Nahm equation~\cite{Nahm:eqns}.

In addition to the physical motivation, the selfdual Yang--Mills equation has attracted interest because of its good
integrability properties, which are inherited by the Bogomolny, Hitchin and Nahm equations, and their analogues.
Further integrable systems may be obtained by reducing the selfdual Yang--Mills equation by other groups of conformal
transformations of ${\mathbb{R}}^4$ or ${\mathbb{R}}^{2,2}$, and many such reductions have been
investigated~\cite{MaWo:ist}.
For example, hyperbolic monopoles arise from the reduction to three dimensions by a~rotation, the Ernst equation is
a~reduction to two dimensions by a~translation and rotation, while reductions to one dimension may be interpreted as
isomonodromic deformation problems with four poles, governed (in the generic case) by the Schlesinger
equation~\cite{Hit:tem,MaWo:ist,MMW:sbm}.

Twistor theory gives one explanation for this integrability: there is a~Ward correspondence (see~\cite{AtWa:iag})
between solutions of the selfdual Yang--Mills equation on f\/lat space and holomorphic vector bundles on (suitable open
subsets of) $\CP3$.
This suggests that the selfdual Yang--Mills equation will continue to be integrable on other spaces $M$ so long as there
is still a~Ward correspondence between solutions and holomorphic vector bundles on some complex $3$-mani\-fold~$Z$, the
\emph{twistor space} of~$M$.
Such curved twistor spaces~$Z$ were introduced by Penrose~\cite{Pen:nlg} to study selfdual vacuum metrics.
(Note that in Euclidean signature, selfdual Ricci-f\/lat metrics are locally hyperk\"ahler.)

Deep relationships between gauge f\/ield equations and selfdual vacuum metrics have been observed in a~number of places.
In~\cite{Ward:suc}, Ward considered gauge theories on ${\mathbb{R}}^{4-\ell}$ with the gauge group being a~transitive
group of dif\/feomorphisms of an $\ell$-manifold $\Sigma^\ell$.
He focused on the group of dif\/feomorphisms preserving a~f\/ixed volume form on $\Sigma^\ell$ and observed that gauge
f\/ields then give rise to selfdual vacuum metrics on ${\mathbb{R}}^{4-\ell}\mathbin{{\times}}\Sigma^\ell$.
When $\ell=1$, the group of volume preserving dif\/feomorphisms of a~circle or a~line is $\Un(1)$ or ${\mathbb{R}}$ acting
by translation, so the monopoles are Abelian, and Ward's construction reduces to the Gibbons--Hawking
Ansatz~\cite{GiHa:gmi}.
For $\ell=2$, f\/inite-dimensional subgroups of $\SDiff(\Sigma^2)$ yield other interesting constructions of selfdual
vacuum metrics~\cite{DMW:2dg}, while $\ell=3$ gives the Ashtekar--Jacobson--Smolin description of selfdual vacuum
metrics in terms of the Nahm equation~\cite{AJS:hfe}.
One can also view the (closely related) Mason--Newman formulation~\cite{MaNe:eym} as the case $\ell=4$.

In Ward's construction, gauge f\/ields on lower-dimensional f\/lat spaces give rise to curved $4$-manifolds.
On the other hand, the example of hyperbolic monopoles shows that the lower-dimensional spaces can also be curved.

One of the goals of this paper is to place these miscellaneous results and observations in a~geo\-metric framework which
simultaneously unif\/ies and generalizes them.
Selfdual vacuum metrics do not provide the right setting for this.
For example, the natural generalization of the Gibbons--Hawking Ansatz to Abelian monopoles on hyperbolic space is
LeBrun's hyperbolic Ansatz for scalar-f\/lat K\"ahler  metrics~\cite{LeBr:cp2}, and even on ${\mathbb{R}}^{4-\ell}$, one
f\/inds that without the volume-preserving condition, Ward's construction leads to hypercomplex, rather than hyperk\"ahler,
structures~-- see Hitchin~\cite{Hit:hcm} and Joyce~\cite{Joy:esd} (or Dunajski~\cite{Dun:tph}) for the analogue of the
Ashtekar--Jacobson--Smolin and Mason--Newman description respectively.

In order to incorporate these constructions, it is essential to take into account a~second basic feature of the selfdual
Yang--Mills equation, in addition to integrability: \emph{conformal invariance}.
The selfdual Yang--Mills equation makes sense on any oriented conformal $4$-manifold $M$, and the purely conformal part
of Penrose's nonlinear graviton shows that there is a~curved twistor space $Z$ as long as the conformal structure on $M$
is selfdual~\cite{AHS:sd4, Pen:nlg}.

\looseness=-1
The thesis of this paper is that reductions of selfdual conformal geometry to lower dimensions provide ``integrable
background geometries'' on which curved versions of the Bogomolny, Hitchin and Nahm equations (with good integrability
properties) can be def\/ined.
This thesis can be illustrated by the three-dimensional case: the Jones--Tod correspondence~\cite{JoTo:mew} shows that
the reduction of a~selfdual conformal structure to three dimensions is an Einstein--Weyl structure, and selfdual
Yang--Mills f\/ields reduce to generalized monopoles on such Einstein--Weyl spaces, pla\-cing Euclidean and hyperbolic
monopoles in a~common framework.
Selfdual spaces with symmetry over a~given Einstein--Weyl space are built out of Abelian monopoles, and this provides
many constructions (generalizing the Gibbons--Hawking and hyperbolic Ans{\"a}tze) of hyperk\"ahler, scalar-f\/lat
K\"ahler  and selfdual Einstein metrics, and also of hypercomplex structures~\cite{CaPe:sdc,GaTo:hms,Hit:cme,LeBr:cp2}.

The role of gauge f\/ields and Abelian monopoles in Ward's construction and the Jones--Tod construction respectively
suggests an underlying principle relating gauge f\/ield equations and constructions of selfdual spaces.
This is amplif\/ied by the following two generalizations of the Jones--Tod construction.

First, in~\cite{Cal:sde}, the Jones--Tod construction and Ward's construction on ${\mathbb{R}}^3$ were simultaneously
generalized by considering the Einstein--Weyl Bogomolny equation with the gauge group acting transitively by
dif\/feomorphisms on a~circle or a~line.
The volume preserving case gives the usual Jones--Tod correspondence between selfdual spaces with symmetry and Abelian
monopoles on Einstein--Weyl spaces, but more general gauge groups lead to new constructions of selfdual conformal
structures and metrics, including, in special cases, hyperk\"ahler  and selfdual Einstein metrics.

Second, note that the Jones--Tod construction gives rise to a~procedure for constructing a~new selfdual space from an
invariant selfdual Maxwell f\/ield on a~given selfdual space with one-dimensional symmetry group: the selfdual Maxwell
f\/ield descends to an Abelian monopole on the quotient Einstein--Weyl space, out of which a~new selfdual space may be
built.
This was generalized, by Maszczyk, Mason and Woodhouse~\cite{MMW:sbm}, to any freely acting symmetry group, using
a~construction they call the ``switch map''~\cite{MaWo:ist}: given a~selfdual space with freely acting group of
conformal transformations $H$, and an invariant selfdual Yang--Mills connection on a~bundle $P$ with a~gauge group $G$
of the same dimension as $H$, the quotient of $P$ by $H$ is another selfdual space, with symmetry group $G$.
Hence, for example, $T^3$-invariant $\SU(2)$ Yang--Mills f\/ields on ${\mathbb{R}}^4$ give rise to selfdual conformal
structures with $\SU(2)$ symmetry, such as the scalar-f\/lat K\"ahler, hypercomplex and selfdual Einstein metrics
of~\cite{Dan:sfk,Hit:tem,Hit:hcm,PePo:kzs,Tod:p3}.

\subsection*{Overview} Ward's construction, the (generalized) Jones--Tod correspondence, and the switch map all point to
the following framework for dimensional reduction of selfdual conformal geometry and the selfdual Yang--Mills equation.
\begin{enumerate}\itemsep=0pt
\item[(i)]
There are background geometries in each dimension less than four obtained by dimensional reduction of the selfduality
condition for conformal structures.
\item[(ii)]
The nonlinear dif\/ferential equations def\/ining these background geometries have the surprising feature that they do not
depend on the symmetry group that one reduces by.
Therefore, including selfdual conformal geometry, there are only four kinds of (non\-de\-ge\-ne\-rate) background geometry.
\item[(iii)]
Instead, the symmetry group enters as a~gauge group for a~gauge theory def\/ined on the background geometry, and the gauge
f\/ield equation is the dimensional reduction of the selfdual Yang--Mills equation.
Hence the gauge f\/ield equations play a~remarkable dual role: solutions give rise both to selfdual conformal
$4$-manifolds and also to selfdual Yang--Mills f\/ields on such manifolds.
\end{enumerate}
In this paper, the above framework is established and studied in full generality.
Furthermore, the dif\/ferent geometries are related not just by symmetry reduction, but by a~more general form of
dimensional reduction, of which Ward's construction and the generalized Jones--Tod construction are examples
(cf.~also~\cite{GKPS:frsd}).
Such constructions also relate the lower-dimensional geometries to each other.
More precisely:
\begin{itemize}\itemsep=0pt
\item
For $1\leqslant k<k+\ell\leqslant4$, $(k+\ell)$-dimensional geometries are obtained from a~$k$-dimensional geometry by
solving the gauge f\/ield equation on that background where the gauge group acts transitively by dif\/feomorphisms on an
$\ell$-manifold.
\end{itemize}

In four dimensions, the gauge f\/ields and background geometries are, of course, selfdual Yang--Mills f\/ields on selfdual
spaces, while in three dimensions, one obtains monopoles on Einstein--Weyl spaces.
The one and two-dimensional stories are new, although these structures have been implicitly studied in many places, at
least in special cases~\cite{HYMO:aag,Taf:2dr} and the twistor theory of the two-dimensional geometry has been developed
independently by Donaldson and Fine~\cite{DoFi:tasd,Fin:tasde}.
I~f\/irst describe brief\/ly the two-dimensional geometries, in the form obtained by reduction from Euclidean signature.
The two reductions from Kleinian signature $(2,2)$ are similar.

Given a~complex line bundle ${\mathcal{W}}$ over a~Riemann surface $N$ with a~Hermitian metric on
${\mathcal{W}}^{*}\otimes TN$, the geometric structure is a~triple consisting of a~$\Un(1)$-connection $a$ on
${\mathcal{W}}^{*}\otimes TN$, a~section $\psi$ of ${\mathcal{W}}^{*}$ and a~section ${\mathcal{C}}$ of
${\mathcal{W}}^{*}\otimes{\mathcal{W}}^{*}\otimes TN$ satisfying the following equations:
\begin{gather*}
\overline\partial{}^a {\mathcal{C}}=0,
\qquad
\overline\partial{}^a \psi=-3{\mathcal{C}}\overline\psi,
\qquad
{*F^a}=|\psi|^2-2|{\mathcal{C}}|^2.
\end{gather*}
The f\/lat geometry is obtained by setting ${\mathcal{W}}=TN$, with the trivial connection on ${\mathcal{W}}^{*}\otimes
TN$ and ${\mathcal{C}}=\psi=0$.
This is simply a~Riemann surface with no additional structure.
The general two-dimensional background geometry was found in joint work with Lionel Mason~\cite{CaMa:svm}.

The gauge f\/ields on this background, with gauge group $G$, are pairs $(A,\Phi)$ consisting of a~$G$-connection $A$ and
a~section $\Phi$ of ${\mathcal{W}}^{*}\otimes\lie{g}_N$, where $\lie{g}_N$ is the associated Lie algebra bundle.
These pairs satisfy the following equations:
\begin{gather*}
F^A-[\Phi,\overline\Phi]=\psi\wedge\overline\Phi+\overline\psi\wedge\Phi,
\qquad
\overline\partial{}^{a,A}\Phi={\mathcal{C}}\overline\Phi.
\end{gather*}
On the f\/lat geometry, where ${\mathcal{C}}=0=\psi$, these are the Hitchin equations for stable pairs, but in the general
case, the equations are coupled to the background geometry.

In one dimension, the background geometry is governed by a~symmetric traceless $(3\times 3)$-matrix ${\mathcal{B}}$
satisfying the Riccati equation ${\mathcal{B}}_r=2({\mathcal{B}}^2)_0$, where $r$ is an af\/f\/ine coordinate and the
subscript zero denotes the traceless part.
More invariantly, ${\mathcal{B}}$ is a~section of $T^{*} C\otimes\Sym_0{\mathcal{E}}$ where ${\mathcal{E}}$ is a~rank
$3$ conformal vector bundle over a~curve $C$, with $\wedge^3{\mathcal{E}}=(TC)^3$, and one f\/ixes compatible
connections on ${\mathcal{E}}$ and $TC$ to def\/ine the $r$-derivative ${\mathcal{B}}_r$.

The gauge f\/ields on this background are sections $\Phi$ of ${\mathcal{E}}^{*}\otimes\lie{g}_C$ (where $\lie{g}_C$ is
a~Lie algebra bundle over $C$) satisfying the equation $\Phi_r-{*[\Phi,\Phi]}={\mathcal{B}}\mathinner{\cdot}\Phi$, where
$*$ denotes the star operator on ${\mathcal{E}}$ and one f\/ixes a~$G$-connection to def\/ine $\Phi_r$.
The f\/lat geometry is the trivial solution ${\mathcal{B}}=0$, in which case this gauge f\/ield equation reduces to the Nahm
equation.

There is one further important property of the Jones--Tod correspondence that continues to hold in the general
framework: it is constructive, i.e., a~selfdual conformal structure is \emph{explicitly} determined by an Einstein--Weyl
structure together with a~monopole, and conversely.
The correspondence between additional monopoles and selfdual Yang--Mills f\/ields is equally explicit.
The same remarks hold for the switch map and Ward's construction.
I therefore present explicit formulae for the constructions of the paper.

For convenience of exposition, I will concentrate on the reductions from Euclidean signature and the notation will be
adapted to this case.
However, the \emph{non-null} reductions from Kleinian signature $(2,2)$ are completely analogous, as are such reductions
of complex geometries: I indicate throughout the nondegeneracy assumptions that need to be made, and any changes in
notation that are needed.
On the other hand, it remains an interesting open project to study the integrable background geometries arising from
\emph{null} reductions: only two such reductions are considered here.

\subsection*{User guide}

This is a~long paper to read from start to f\/inish, so I give a~detailed guide to the sections, both to draw attention to
the highlights and to enable the reader to dip into the paper more easily.

The general theory of the paper is developed in Sections~\ref{s:sdbg}, \ref{s:gfe}, \ref{s:sgf} and~\ref{s:bgf}.
Section~\ref{s:sdbg} concerns selfdual spaces with a~freely acting symmetry group and presents the background geometry
equations in dimensions $1$--$3$.
In each case the main result, Theorems~\ref{th:ric},~\ref{th:sv} or~\ref{th:JnT}, identif\/ies the selfduality condition on
a~conformal structure with the background geometry equations: in particular Theorem~\ref{th:JnT} is the Jones--Tod
correspondence~\cite{JoTo:mew}.
Section~\ref{s:gfe} deals with invariant selfdual Yang--Mills f\/ields on selfdual spaces with a~freely acting symmetry group.
Here the gauge f\/ield equations on the background geometries are computed.

The main general theorems are in Sections~\ref{s:sgf} and~\ref{s:bgf}, where the inverse constructions to
Section~\ref{s:sdbg} are established and generalized.
Theorem~\ref{th:sdgf} proves that selfdual spaces may be constructed from solutions to gauge f\/ield equations on
$k$-dimensional background geometries where the gauge group acts transitively on a~$(4-k)$-manifold: in the case of
a~group acting on itself by the regular representation, this theorem reconstructs the selfdual spaces with freely acting
symmetry group of Section~\ref{s:sdbg}~-- although, as mentioned already, in the general case, the group need not even be
f\/inite-dimensional.

Theorem~\ref{th:bgf} generalizes all this to gauge groups acting on $\ell$-manifolds, giving explicit constructions of
$(k+\ell)$-dimensional background geometries from gauge f\/ields on $k$-dimensional geometries for $k+\ell\leqslant4$.
Unfortunately, the calculations here are too complicated to present in full.
However, I do provide explicit formulae: in examples arising in practice, it is usually not too hard to verify that the
$(k+\ell)$-dimensional geometry satisf\/ies the background equations, once one has a~formula.

The relation between selfdual spaces, and the three-dimensional background geometries, Einstein--Weyl spaces, is
a~longstanding one, and as a~consequence the latter have been extensively studied~\cite{Cal:sde,CaPe:sdc,CaTo:emh,
Gau:swe,GaTo:hms, Hit:cme,JoTo:mew,LeBr:cp2,Tod:p3,Tod:sew,Ward:sut}.
By contrast, the background geometries in one and two dimensions, Riccati spaces, and spinor-vortex spaces, although
implicitly underlying previous work, have only been introduced and investigated relatively recently (see,
e.g.,~\cite{CaMa:svm,DoFi:tasd}).
Assuming that one is not interested in zero-dimensional dif\/ferential geometry (see Remark~\ref{r:zdg}), then these
geometries form the foundation for the more well-known higher-dimensional structures (using Theorems~\ref{th:sdgf}
and~\ref{th:bgf}).

In Section~\ref{s:rs}, the geometry of Riccati spaces is described.
Although there are only six Riccati spaces up to local isomorphism, they have a~rather rich structure, which is most
easily revealed in a~complexif\/ied setting, since only three of the Riccati spaces arise as reductions from Euclidean
signature.
The six solutions correspond to the six types of quartic polynomial, i.e., the f\/ive types of conf\/iguration of four
points on $\CP1$ together with the zero polynomial.
These types are denoted (I,~II, III, D, N, 0), following the well-known application of this classif\/ication to Weyl tensors.
It is easy to see that the gauge f\/ield equation on the trivial (0) Riccati space is the Nahm equation.
One can also observe that the gauge f\/ield equations on the nontrivial Riccati spaces are equivalent to (strong)
isomonodromy equations for a~connection with four poles on~$\CP1$, the conf\/iguration of the poles corresponding to the
type of Riccati space.
This leads to a~uniform Lax pair for these problems, which reduces to the usual Lax pair for the Nahm equation in the
trivial case.

The theory of this paper would be very dry without examples and applications, so I intersperse the main
development with \emph{interludes}, which motivate or illustrate the theory, yet are, to varying degrees,
self-contained.
The f\/irst interlude, Section~\ref{s:bm}, relates the approach of Section~\ref{s:sdbg} to other studies of selfdual
Bianchi metrics~\cite{Dan:sfk,DaSt:cok,PePo:kzs,Tod:p6,Tod:com,Tod:p3}.

The second interlude, Section~\ref{s:svhe}, provides a~simple \emph{a priori} explanation for the conformal invariance
of Hitchin's selfduality equation on a~Riemann surface~\cite{Hit:sde} (and of course, reductions from Kleinian signature
also give conformally invariant equations: harmonic maps into a~Lie group, and the principal chiral model).
This demystif\/ies this conformal invariance: it is, after all, a~consequence of conformal invariance in four dimensions.

The third interlude, Section~\ref{s:hchk}, is a~unif\/ied treatment of various constructions (or at least interpretations)
of hypercomplex and hyperk\"ahler  structures.
My aim here is three-fold:
\begin{enumerate}\itemsep=0pt
\item[(i)]
to provide geometrical descriptions of the well-known Mason--Newman~\cite{MaNe:eym},
Ashtekar--Jacobson--Smolin~\cite{AJS:hfe}, and Park--Ward~\cite{Par:sdg,Ward:suc} constructions of hyperk\"ahler
metrics from lower-dimensional gauge f\/ields;
\item[(ii)]
to give, at the same time, hypercomplex generalizations, following~\cite{Dun:tph,GrSt:his,Hit:hcm,Joy:esd};
\item[(iii)]
to prove that all hypercomplex and hyperk\"ahler  structures are locally obtained from any of these constructions, in
a~way that is manifestly compatible with any reality conditions, and that clarif\/ies the extra choice that needs to be
made to reduce a~hypercomplex or hyperk\"ahler  structure to a~solution of the relevant gauge f\/ield equation.
\end{enumerate}
I hope that the overview provided here is useful, at least to the reader who is not familiar with the treatments in the
physics literature.
In addition this work answers~-- and extends to the hypercomplex case~-- a question of Ward~\cite{Ward:suc}, who
conjectured that any hyperk\"ahler  metric could be obtained from Hitchin's selfduality equation, with gauge group
$\SDiff(\Sigma^2)$.
I also show how this description gives a~Euclidean analogue of Plebanski's heavenly equations, which is well-adapted to
the study of hyperk\"ahler  metrics on elliptic f\/ibrations.

The fourth interlude, Section~\ref{s:s1Hit}, is a~two-part analysis of hyperCR Einstein--Weyl spaces, the three-dimensional analogue of hypercomplex structures.
In the f\/irst part, following an approach of Tod~\cite{Tod:sew}, the hyperCR Einstein--Weyl equation is shown to be
equivalent to the $\Diff(S^1)$ Hitchin equation (revealing a~hidden $\SO(3)$ symmetry in the latter).
This can be viewed as a~three-dimensional version of the constructions of Section~\ref{s:hchk} although it is remarkable
that the f\/ibres are only one-dimensional, since such constructions are not suf\/f\/iciently general in four dimensions.
In the second part, Einstein--Weyl spaces admitting a~dimensional reduction with geodesic one-dimensional f\/ibres are
studied.
They are shown to be hyperCR, and the quotient spaces are trivial or spherical spinor-vortex geometries.
Conversely, any shear-free geodesic congruence on a~hyperCR Einstein--Weyl space def\/ines a~generalized dimensional
reduction.

The f\/inal interlude again consists of two-parts.
The f\/irst shows that hyperCR Einstein--Weyl spaces may also be constructed from the $\Diff(\Sigma^2)$ Nahm equations.
In the second part, I present a~proof that the well-known $\SU(\infty)$ Toda f\/ield equation $u_{xx}+u_{yy}+(e^u)_{zz}=0$
and the dKP equation $u_{yy}=(u_t-uu_x)_x$ are both equivalent to generalized Nahm equations, via a~hodograph
transformation (cf.~\cite{DuTo:p12} in the dKP case).
Even these results are to a~large extent self-contained, although the the work of Section~\ref{s:bgf} shows that the
construction is a~special case of the general theory, while Section~\ref{s:rs} shows that the backgrounds for these
generalized Nahm equations are the type (D) and (N) Riccati spaces respectively.
To the best of my knowledge, all examples of hodograph solutions to the $\SU(\infty)$ Toda f\/ield equation or dKP
equation arise in this way from explicit solutions of a~generalized Nahm equation.
I end the interlude by discussing these examples.

\subsection*{Addenda}

The majority of this paper was written in the period 1999--2001, and the present content is not
substantially dif\/ferent from a~January 2002 version which has been posted on my academic home page since that time.
The intervening 12 years have seen many advances in the f\/ield (for instance, by Dunajski and his
collaborators~\cite{Dun:hfcq,Dun:sitt,DGS:misd,DuKr:ewhv,DuSp:dis1}), and I have collected a~few of the most closely
related papers as ``Additional references'' at the end of the bibliography.
Some of these works develop ideas from the 2002 version of this paper, or discover related ideas independently.

There has also been much work on null reductions and geometries, which were only touched upon in the original version of
this paper.
A major driving force has been the introduction of a~global twistor theory of holomorphic discs by LeBrun and
Mason~\cite{LeMa:zmcs,LeMa:nghd,LeMa:zmbchd,Nak:ssdz}, both for Kleinian signature selfdual conformal structures and
$2$-dimensional projective structures.
Dunajski and West~\cite{DuWe:acs,DuWe:acns} established a~relationship between these structures by considering selfdual
conformal structures with a~null conformal Killing vector f\/ield.
As explained in~\cite{Cal:sdp}, the natural context for their construction is the null reduction of selfdual conformal
structures along a~$\beta$-surface foliation.
The quotient is a~surface with a~natural projective structure structure, and the gauge f\/ield equations are projectively
invariant on this background geometry.
In~\cite{Nak:sdza} the global Mason--LeBrun theories are related by this construction.
More recently, similar methods have been applied in Einstein--Weyl geometry~\cite{LeMa:ewhd,Nak:ewlm}.

One of the most intriguing recent developments has been the introduction, by Ferapontov and his coworkers, of the method
of hydrodynamic reductions to analyse integrability.
This was applied in~\cite{FHZ:cqim} to study the central quadric Ansatz and its relation to the Painlev\'e equations.
In doing so, they independently rediscovered equations equivalent to the $\SDiff(\Sigma^2)$ generalized Nahm equations
on Riccati spaces.
More recently, in~\cite{FeKr:disew}, Ferapontov and Kruglikov construct a~Weyl structure from the formal linearization
of a~second-order PDE in three dimensions, and show that it is Einstein--Weyl for all solutions if and only if the
system is integrable.
This and related conjectures for second-order PDEs in four dimensions suggest deep connections with the integrable
background geometry concept.

For these reasons, the main (odd-numbered) sections have been updated with addenda which place the above works in the
framework of this paper.
Section~\ref{section12} and Subsection~\ref{s:sdsn}
have also been updated to ref\/lect the exciting new directions that are currently unfolding.

\subsection*{Notation}

In order to manipulate, in a~tensorial way, the objects and structures entering into the equations and constructions of
this paper, it will be convenient to employ the formalism of densities.
If $V$ is a~real $n$-dimensional vector space and $w$ any real number, then the oriented one-dimensional linear space
$L^w=L^w_V$ carrying the representation $A\mapsto|\det A|^{w/n}$ of $\GL(V)$ is called the space of \emph{densities of
weight $w$} or \emph{$w$-densities}.
It can be constructed canonically as the space of maps $\rho\colon(\wedge^nV)\setminus0\to{\mathbb{R}}$ such that
$\rho(\lambda\omega)=|\lambda|^{-w/n}\rho(\omega)$ for all $\lambda\in{\mathbb{R}}^{\times}$ and
$\omega\in(\wedge^nV)\setminus0$.

For a~vector bundle ${\mathcal{V}}\to M$ this construction yields an oriented real line bundle $L^w_{\mathcal{V}}$,
a~\emph{density line bundle}.
If ${\mathcal{V}}$ is oriented and of rank $n$, then $L^{-n}_{\mathcal{V}}$ is canonically isomorphic to
$\wedge^n{\mathcal{V}}$; indeed an orientation may be def\/ined as an \emph{orientation form}
$*1\in\mathrm{C}^\infty(M,L^n_{\mathcal{V}} \wedge^n{\mathcal{V}}^{*})$.
(Here and elsewhere, when tensoring with a~density line bundle, I shall often omit the tensor product sign.) More
generally, the \emph{Hodge star operator} is the isomorphism
\begin{gather*}
{*}\colon \ L^{w-k}_{\mathcal{V}} \wedge^k{\mathcal{V}}\to L^{w+n-k}_{\mathcal{V}} \wedge^{n-k}{\mathcal{V}}^{*}
\end{gather*}
determined by the nondegenerate pairing $\wedge^k{\mathcal{V}}^{*}\otimes\wedge^{n-k}{\mathcal{V}}^{*}\to
\wedge^n{\mathcal{V}}^{*}\cong L^{-n}_{\mathcal{V}}$.

A \emph{conformal structure} on ${\mathcal{V}}$ is a~nondegenerate symmetric bilinear form on ${\mathcal{V}}$ with
va\-lues in~$L^2_{\mathcal{V}}$, or equivalently a~metric on $L^{-1}_{\mathcal{V}} {\mathcal{V}}$.
The conformal inner product of sections~$X$,~$Y$ is $\ip{X,Y}\in\mathrm{C}^\infty(M,L^2_{\mathcal{V}})$ and the conformal
structure itself may be viewed as a~section $\mathsf{c}\in\mathrm{C}^\infty(M,L^2_{\mathcal{V}}
S^2{\mathcal{V}}^{*})$.
I shall make free use of the isomorphism between $L^{w-k}_{\mathcal{V}}\wedge^k{\mathcal{V}}$ and
$L^{w+k}_{\mathcal{V}} \wedge^k{\mathcal{V}}^{*}$ given by a~conformal structure.

When ${\mathcal{V}}$ is the tangent bundle of $M$, $L^w_M=L^w_{TM}$ is called the bundle of $w$-densities of $M$,
denoted $L^w$ when $M$ is understood.
The line bundles $L^w$ are trivializable and a~nonvanishing (usually positive) section $\mu$ of $L=L^1$ will be called
a~\emph{length scale} or \emph{gauge}.
I shall also say that tensors in $L^w\otimes(TM)^j\otimes(T^{*} M)^k$ have \emph{weight} $w+j-k$.

A \emph{Weyl derivative} is a~covariant derivative $D$ on $L$.
It induces covariant derivatives on~$L^w$ for all $w$.
If $M$ is conformal, i.e., there is a~conformal structure $\mathsf{c}$ on~$TM$, then any Weyl derivative induces (via
the Koszul formula) a~\emph{Weyl connection}: a~torsion-free connection~$D$ on~$TM$ with $D\mathsf{c}=0$.
Compatible Riemannian metrics~$g$ correspond to length scales $\mu$, and the Weyl connection induced by the Weyl
derivative preserving $\mu$ is, of course, the Levi-Civita connection of $g=\mu^{-2}\mathsf{c}$.

The curvature $R^D$ of a~Weyl connection $D$, as a~$\mathop{\mathfrak{co}}(TM)$-valued $2$-form, decomposes as
\begin{gather*}
R^D_{X,Y}=W_{X,Y}-\abrack{r^D(X),Y}+ \abrack{r^D(Y),X}.
\end{gather*}
Here $W$ is the \emph{Weyl curvature} of the conformal structure, an $\mathop{{\rm so}}(TM)$-valued $2$-form
independent of the choice of $D$, and $r^D$ is a~covector valued $1$-form, the \emph{$($normalized$)$ Ricci
curvature} of~$D$.
For a~$1$-form $\gamma$ and tangent vector~$X$, $\abrack{\gamma,X}=\gamma(X)\iden+\gamma\mathinner{\vartriangle}X$, where
$(\gamma\mathinner{\vartriangle}X)(Y) =\gamma(Y)X-\ip{X,Y}\gamma$.
This bracket is part of a~Lie algebra structure on $TM\oplus\mathop{\mathfrak{co}}(TM)\oplus T^{*} M$ and the same
notation will be used for the commutator bracket in $\mathop{\mathfrak{co}}(TM)$.

The normalized Ricci curvature decomposes into a~symmetric traceless part $r^D_0$, a~scalar part $\scal^D$ (the scalar
curvature) and a~skew part, which is just a~multiple of the curvature of~$D$ on~$L$ called the \emph{Faraday curvature}~$F^D$.
In practice a~Weyl derivative is described by its connection $1$-form $\omega$ relative to a~length scale: $\omega$ is
called the \emph{Weyl $1$-form}, and $F^D={\rm d}\omega$.
If $F^D=0$ then $D$ is said to be \emph{closed}.
It follows that there are local length scales $\mu$ with $D\mu=0$.
If such a~length scale exists globally then $D$ is said to be \emph{exact}.

The above constructions can also be carried out locally on complex manifolds, except that~$L^1$ is now a~choice of local
$n$th root of $\wedge^nTM$.

Although confusion with the index $i$ is unlikely, as a~courtesy to the reader, I denote the (chosen) square root of~$-1$ by ${\boldsymbol i}$.

\subsection*{Twistors and Lax pairs} Several of the results in this paper were motivated by twistor or integrable
systems methods: in particular the idea of generalized dimensional reduction arises naturally when one considers
holomorphic foliations of twistor spaces.
However, I have deliberately suppressed discussion of twistor spaces and Lax pairs, for at least two reasons: f\/irst, to
make the paper accessible to the reader not familiar with these ideas; second, because I believe it is a~useful to
present all calculations and formulae in purely dif\/ferential geometric terms~-- it is often impossible to carry out
twistor constructions in practice.
An unfortunate consequence is that some of the results and formulae appear miraculous: the twistor point of view
provides a~quick way to see why such results are true, while the Lax pair formalism provides one way to carry out more
detailed calculations.

\section{Selfdual spaces and the background geometries}
\label{s:sdbg}

On a~conformal $4$-manifold, the Hodge star operator is an involution on $2$-forms, so there is a~decomposition
$\wedge^2 T^{*} M= \wedge^2_\sd T^{*} M\oplus \wedge^2_\asd T^{*} M$, and the eigenspaces $\wedge^2_\sdasd T^{*}
M$ are called the selfdual and antiselfdual $2$-forms.
This induces a~similar decomposition
$\mathop{{\rm so}}(TM)=\mathop{{\rm so}}_\sd(TM)\oplus\mathop{{\rm so}}_\asd(TM)$ of the skew
endomorphisms of $TM$.
The Weyl curvature $W$ splits as a~sum of selfdual and antiselfdual $2$-forms $W^\sdasd$ with values in
$\mathop{{\rm so}}(TM)$: in fact $W^\sd$ is $\mathop{{\rm so}}_\sd(TM)$-valued and $W^\asd$ is
$\mathop{{\rm so}}_\asd(TM)$-valued (essentially because $W$ is traceless).

A \emph{selfdual space} is a~conformal $4$-manifold with selfdual Weyl curvature, i.e., $W^\asd=0$.

In~\cite{Joy:esd}, Joyce studied selfdual spaces with a~surface-orthogonal action of the torus $T^2$ by conformal
transformations, and found new explicit selfdual conformal metrics on connected sums of complex projective planes.
The key idea in his approach is the use of conformal connections with torsion, and the following observation.

\begin{lem}[see Joyce~\cite{Joy:esd}]\label{lem:J}
Let $(M,\mathsf{c})$ be an oriented conformal $4$-manifold and ${\mathcal{D}}$ a~conformal connection $($i.e.,
${\mathcal{D}}\mathsf{c}=0$, but ${\mathcal{D}}$ may have torsion$)$.
Suppose that the antiselfdual part of the torsion of ${\mathcal{D}}$ is tracelike.
Then $(M,\mathsf{c})$ is selfdual if and only if the Weyl part of the curvature of~${\mathcal{D}}$ is selfdual.
\end{lem}

When the torsion is selfdual and traceless, this lemma follows easily, since ${\mathcal{D}}$ then dif\/fers from a~Weyl
connection by an $\mathop{{\rm so}}_\sd(TM)$-valued $1$-form, so that the $\mathop{{\rm so}}_\asd(TM)$-valued
part of the curvature of ${\mathcal{D}}$ agrees with that of the Weyl connection.
The general case is a~consequence of the fact that the trace parts of the torsion (which are $1$-forms) cannot
contribute to the Weyl part of the curvature.

I shall refer to this result as Joyce's lemma: although simple, and perhaps previously known, its application by
Joyce~\cite{Joy:esd} was one of the main motivations for the present work.
Indeed, for a~conformal manifold with a~surface-orthogonal $T^2$-action, Joyce constructed, on the open set where the
torus acts freely, a~conformal connection with torsion, and hence separated the selfduality equation for the conformal
structure into a~nonlinear equation for a~quotient geometry and a~linear equation def\/ined on this background.
He then showed that the quotient geometry in this case is the hyperbolic plane, and superposed known solutions of the
linear equation to f\/ind new explicit metrics.

In this section, Joyce's lemma will be applied to selfdual spaces with \emph{any} freely acting group of conformal
transformations, and a~large class of dimensional reductions of the selfduality equation will be obtained.
The restriction to freely acting groups of symmetries will also be relaxed later, leading to a~generalized version of
dimensional reduction.

Let $M$ be a~conformal manifold with a~free proper action of a~group $H$, so that $M$ is a~principal $H$-bundle over the
orbit space $Q=M/H$.
The generators of the action form a~Lie algebra $\lie{h}$ of vector f\/ields on $M$ and pointwise evaluation def\/ines an
isomorphism $M\mathbin{{\times}}\lie{h}\to VM$, where $VM$ is the vertical bundle of $M\to Q$.
Now suppose that $H$ acts conformally with nondegenerate orbits.
Then $Q$ is a~conformal manifold and the horizontal distribution $VM^\perp\leqslant TM$ def\/ines a~principal
$H$-connection $\alpha$ on $\pi\colon M\to Q$; furthermore, $VM$ is isomorphic to the pullback of a~conformal vector
bundle ${\mathcal{V}}\to Q$ with $L^1_{{\mathcal{V}}}=L^1_Q$.
Although $M\mathbin{{\times}}\lie{h}$ and $\pi^*{\mathcal{V}}$ are both isomorphic to~$VM$, it will be crucial in the
following to distinguish between them, since the trivialization of~$VM$ given by the $H$-action is not, in general,
compatible with the conformal structure.
However, the isomorphism $\pi^*{\mathcal{V}}\to M\mathbin{{\times}}\lie{h}$ is $H$-equivariant, so it may be viewed as
a~bundle isomorphism $\varphi\colon{\mathcal{V}}\to\lie{h}_Q$ over $Q$, where $\lie{h}_Q=M\mathbin{{\times}}_H\lie{h}$.
To summarize, the conformal geometry of $M$ is encoded by:
\begin{itemize}\itemsep=0pt
\item
a conformal structure on $Q$;
\item
a conformal vector bundle ${\mathcal{V}}\to Q$ with $L^1_{{\mathcal{V}}}=L^1_Q$;
\item
a principal $H$-connection $\alpha$ on $M\to Q$;
\item
a bundle isomorphism $\varphi\colon{\mathcal{V}}\to\lie{h}_Q$ over $Q$.
\end{itemize}
The data $(\alpha,\varphi)$ identify the tangent bundle $TM$ with the pullback of ${\mathcal{V}}\oplus TQ$: sections of~${\mathcal{V}}$ or~$TQ$ will be denoted $U$, $V$, $W$ or $X$, $Y$, $Z$ respectively, and identif\/ied with invariant vector f\/ields on~$M$.

The curvature of the principal connection $\alpha$ on $\pi\colon M\to Q$ is an $\lie{h}_Q$-valued $2$-form $F^\alpha$,
given by minus the Frobenius curvature of the horizontal distribution: $F^\alpha(X,Y)=-\varphi([X,Y])$, where~$\varphi$
extended by zero from ${\mathcal{V}}$ to ${\mathcal{V}}\oplus TQ$.

Choose a~conformal connection $D$ on ${\mathcal{V}}$ over $Q$.
This induces a~Weyl derivative on $L^1_{\mathcal{V}}=L^1_Q$, hence a~torsion-free conformal connection on $TQ$ and
a~direct sum connection on ${\mathcal{V}}\oplus TQ$.
These conformal connections will be denoted by $D$, as will the pullback connection on $TM=\pi^*({\mathcal{V}}\oplus
TQ)$, which is conformal, but not torsion-free in general: the f\/ibres of $\pi\colon M\to Q$ need not be umbilic unless
$Q$ is three-dimensional, the nonlinear connection on $M\to Q$ need not be f\/lat unless $Q$ is one-dimensional, and the
sections of ${\mathcal{V}}$ parallel along the f\/ibres of $\pi$ will have nontrivial Lie brackets unless $\lie{h}$ is
Abelian.

In order to apply Joyce's lemma, the torsion needs to be reduced.
To do this, introduce $\psi\colon\wedge^2TQ\to{\mathcal{V}}$ and ${\mathcal{C}}\colon S^2_0{\mathcal{V}}\to TQ$, and
def\/ine ${\mathcal{D}}=D+\widehat{\mathcal{C}}+\widehat\psi$, where $\widehat{\mathcal{C}}$ and $\widehat\psi$ are the
sections of $T^{*} M\otimes\mathop{{\rm so}}(TM)$ given by
\begin{gather*}
\widehat {\mathcal{C}}_{U+X} (V+Y)={\mathcal{C}}(U,V)-\ip{{\mathcal{C}}(U,\cdot),Y},
\\
2\widehat\psi_{U+X} (V+Y)= \ip{\psi(X,\cdot),V}+\ip{\psi(Y,\cdot),U}-\psi(X,Y).
\end{gather*}
The idea is that $\psi$ will compensate for the curvature of the horizontal distribution, while ${\mathcal{C}}$ will
of\/fset the traceless second fundamental form of the f\/ibres.
It will then be possible to make the torsion selfdual by the choice of $D$.

The torsion $T^{\mathcal{D}}$ of this modif\/ied conformal connection may be computed by applying it to invariant vector
f\/ields, as long as one is careful that the Lie bracket on invariant vertical vector f\/ields is minus the Lie bracket on
the Lie algebra $\lie{h}$ of generators of the action, hence minus the bracket on the associated Lie algebra bundle
$\lie{h}_Q$ over $Q$
\begin{gather}
\label{eq:t1}
T^{\mathcal{D}}(U,V)=\varphi^{-1}[\varphi(U),\varphi(V)]_{\lie{h}},
\\
T^{\mathcal{D}}(U,X)=\varphi^{-1}D^\alpha_X\varphi(U)-\ip{{\mathcal{C}}(U),X},
\\
T^{\mathcal{D}}(X,Y)=\varphi^{-1}F^\alpha(X,Y)-\psi(X,Y).
\label{eq:t3}
\end{gather}
Here $D^\alpha$ denotes the connection on $\pi^*({\mathcal{V}}\otimes\lie{h}_Q)$ induced by~$D$ and $\alpha$.
Note that the torsion is vertical-valued.

The computation of the curvature is a~little more complicated.
First of all the curvatures of~${\mathcal{D}}$ and~$D$ (on~$TM$) are related by
\begin{gather*}
\nonumber R^{\mathcal{D}}_{U+X,V+Y}=R^D_{U+X,V+Y}+{\rm d}^D(\widehat {\mathcal{C}}+\widehat\psi)_{U+X,V+Y}
+\Abrack{(\widehat{\mathcal{C}}+\widehat\psi)_{U+X}, (\widehat{\mathcal{C}}+\widehat\psi)_{V+Y}}.
\end{gather*}
The second term is the twisted exterior derivative of $\mathop{{\rm so}}(TM)$-valued $1$-forms and must be handled
carefully, since $D$ has torsion.
Relating this to the torsion of ${\mathcal{D}}$ gives
\begin{gather*}
{\rm d}^D(\widehat{\mathcal{C}}+\widehat\psi)_{U+X,V+Y}= D_{U+X}(\widehat{\mathcal{C}}+\widehat\psi)_{V+Y}
-D_{V+Y}(\widehat{\mathcal{C}}+\widehat\psi)_{U+X}
\\
\hphantom{{\rm d}^D(\widehat{\mathcal{C}}+\widehat\psi)_{U+X,V+Y}=}{}
+(\widehat{\mathcal{C}}+\widehat\psi)_{T^{\mathcal{D}}(U+X,V+Y)} +(\widehat{\mathcal{C}}+\widehat\psi)_{\psi(X,Y)}+
(\widehat{\mathcal{C}}+\widehat\psi)_{\ip{{\mathcal{C}}(U),Y}-\ip{{\mathcal{C}}(V),X}}.
\end{gather*}
Notice that the derivatives in vertical directions are zero, since the connection is a~pullback connection and
$\widehat{\mathcal{C}},\widehat\psi$ are invariant.
Also, $R^D$ is horizontal.
Hence putting everything together leads to the following formula for the ``torsion-adjusted'' curvature of
${\mathcal{D}}$
\begin{gather*}
R^{{\mathcal{D}},ta}_{U+X,V+Y}:=
R^{\mathcal{D}}_{U+X,V+Y}-(\widehat{\mathcal{C}}+\widehat\psi)_{T^{\mathcal{D}}(U+X,V+Y)}
\\
\phantom{R^{{\mathcal{D}},ta}_{U+X,V+Y}}~
=R^D_{X,Y}+D_X\widehat{\mathcal{C}}{}_V-D_Y\widehat{\mathcal{C}}{}_U +D_X\widehat\psi{}_V-D_Y\widehat\psi{}_U
\\
\phantom{R^{{\mathcal{D}},ta}_{U+X,V+Y}=}~
{}+\widehat{\mathcal{C}}_{\psi(X,Y)} +\widehat\psi_{\ip{{\mathcal{C}}(U),Y}}-\widehat\psi_{\ip{{\mathcal{C}}(V),X}}
+\abrack{\widehat\psi_{U+X},\widehat{\mathcal{C}}_V} +\abrack{\widehat{\mathcal{C}}_U,\widehat\psi_{V+Y}}
\\
\phantom{R^{{\mathcal{D}},ta}_{U+X,V+Y}=}~
{}+\widehat{\mathcal{C}}_{\ip{{\mathcal{C}}(U),Y}}
-\widehat{\mathcal{C}}_{\ip{{\mathcal{C}}(V),X}}+\abrack{\widehat{\mathcal{C}}_U,\widehat{\mathcal{C}}_V}
+\widehat\psi_{\psi(X,Y)} +\abrack{\widehat\psi_{U+X},\widehat\psi_{V+Y}}.
\end{gather*}
Now suppose that $T^{\mathcal{D}}$ is selfdual.
Then, by Joyce's lemma, the Weyl curvature $W$ of $M$ is selfdual if and only if the Weyl part of this formula is
selfdual.
This condition is a~nonlinear dif\/ferential equation on $Q$ def\/ining a~reduced background geometry.
A key feature of the formula is that the right hand side is manifestly independent of $(\alpha,\varphi)$ and so the
group structure of $H$ decouples from the reduced background geometry.

The details depend on the dimension; in particular, ensuring that the torsion is selfdual constrains $D$.
In the following subsections, I will explain these constraints and obtain the background geometries explicitly.
For later use I will also introduce individual notations for the geometries in each dimension.
As far as possible, the notation will be chosen to be consistent with existing usage.
The case of a~four-dimensional group acting on itself is left as an exercise: a~f\/ield
${\mathcal{Y}}\in\wedge^2_\asd{\mathcal{V}}^{*}\otimes{\mathcal{V}}$ is needed here to make the torsion selfdual; see
Remark~\ref{r:zdg}.

\subsection{Reduction to one dimension}

In this case $Q$ is a~oriented curve, with weightless unit tangent $\xi$, and $\psi=0$.
Write ${\mathcal{C}}(U,V)=\ip{{\mathcal{B}}(U),V}\xi$, so that ${\mathcal{B}}$ is a~symmetric traceless endomorphism of
weight $-1$.
I will denote the curve $Q$ by $C$ and the bundle ${\mathcal{V}}$ by ${\mathcal{E}}$.
Thus ${\mathcal{E}}$ is a~rank $3$ conformal vector bundle over $C$ with $L^1_{\mathcal{E}}=L^1_C=TC$, the last
identif\/ication being given by $\xi$.
The Hodge star operator is an isomorphism between $\wedge^2{\mathcal{E}}$ and $TC\otimes{\mathcal{E}}=L^1{\mathcal{E}}$.

The two components of the torsion of ${\mathcal{D}}$ are given by
\begin{gather*}
T^{\mathcal{D}}(U,V)=\varphi^{-1}[\varphi(U),\varphi(V)]_{\lie{h}},
\qquad
T^{\mathcal{D}}(U,\xi)=\varphi^{-1}D^\alpha_\xi\varphi(U)-{\mathcal{B}}(U).
\end{gather*}
This will be selfdual if and only if
\begin{gather}
\label{eq:riccdef}
\varphi^{-1}\bigl(D^\alpha_\xi\varphi-{*[\varphi,\varphi]_{\lie{h}}}\bigr) ={\mathcal{B}}.
\end{gather}
Conformal connections $D$ on ${\mathcal{E}}$ form an af\/f\/ine space modelled on $T^{*} C
\otimes\mathop{\mathfrak{co}}({\mathcal{E}})$.
Hence $D$ can be chosen uniquely so that the left hand side is symmetric and traceless and this in turn def\/ines
${\mathcal{B}}$.

The two components of the torsion-adjusted curvature of ${\mathcal{D}}$ are
\begin{gather*}
R^{{\mathcal{D}},ta}_{U,V} =\abrack{\widehat{\mathcal{B}}_U,\widehat{\mathcal{B}}_V},
\qquad
R^{{\mathcal{D}},ta}_{U,\xi} =-D_\xi\widehat{\mathcal{B}}_U+\widehat{\mathcal{B}}_{{\mathcal{B}}(U)}.
\end{gather*}
A straightforward calculation shows that the Weyl part is selfdual if and only if
$D_\xi{\mathcal{B}}=2({\mathcal{B}}^2)_0$, where the subscript zero denotes the traceless part.
Since $D$ is f\/lat, this is really just a~Riccati equation for a~$3\times 3$ symmetric traceless matrix.
On the other hand, the orientation of $C$ is not needed if ${\mathcal{B}}$ is viewed as an endomorphism-valued $1$-form
on $C$.

\begin{defn}
Suppose that ${\mathcal{E}}$ is a~rank $3$ conformal vector bundle over a~curve $C$ with $L^1_{\mathcal{E}}=L^1_C$.
Equip ${\mathcal{E}}$ with a~conformal connection $D$ and a~section ${\mathcal{B}}$ of $T^{*}
C\otimes\Sym_0{\mathcal{E}}$.
Then $(C,{\mathcal{E}})$ is said to be a~\emph{Riccati space} if $(D,{\mathcal{B}})$ satisfy the equation
\begin{gather*}
D{\mathcal{B}}=2\big(\mathcal{B}^2\big)_0.
\end{gather*}
Equivalently, with respect to a~conformal trivialization of ${\mathcal{E}}$, the connection $D$ is given by a~section
$(a,\Theta)$ of $T^{*} C\otimes\mathop{\mathfrak{co}}({\mathcal{E}})=T^{*} C\oplus(T^{*}C\otimes\mathop{{\rm so}}({\mathcal{E}}))$ and
\begin{gather*}
\dot{\mathcal{B}}-a{\mathcal{B}}+\abrack{\Theta,{\mathcal{B}}}=2\big(\mathcal{B}^2\big)_0,
\end{gather*}
where the dot denotes dif\/ferentiation with respect to a~compatible coordinate.
In particular, using a~$D$-parallel trivialization of ${\mathcal{E}}$ and an af\/f\/ine coordinate $r$ (i.e.,
$D{\rm d} r=0$), ${\mathcal{B}}_r=2({\mathcal{B}}^2)_0$.
\end{defn}

Joyce's lemma now gives the following result.
\begin{thm}\label{th:ric}
Let $M$ be an oriented conformal manifold with a~$3$-dimensional Lie algebra $\lie{h}$ of linearly independent conformal
vector fields such that the projection $\pi$ onto the space of orbits is a~submersion over a~curve $C$.
Let ${\mathcal{E}}$ be a~rank $3$ conformal vector bundle on $C$ such that $\pi^*{\mathcal{E}}$ is the vertical bundle
of $M$ $($trivialized along the fibres by invariant vector fields$)$ and define $(D,{\mathcal{B}})$ by
equation~\eqref{eq:riccdef} as explained above.

Then $M$ is selfdual if and only if $(C,{\mathcal{E}})$ is a~Riccati space.
\end{thm}

\subsection{Reduction to two dimensions}

The two-dimensional geometry is perhaps the richest.
In the Euclidean case, $Q$ is a~Riemann surface, which will be denoted $N$, and the conformal vector bundle
${\mathcal{V}}$ has rank two, and so will be viewed as a~complex line bundle ${\mathcal{W}}\to N$.
It will be convenient to view $\varphi\colon{\mathcal{W}}\to\lie{h}_N$ as a~complex linear map
$\frac12(\varphi-{\boldsymbol i}\varphi\circ J)\colon{\mathcal{W}}\to\lie{h}_N\otimes{\mathbb{C}}$.
The constraint $L^1_{\mathcal{W}}=L^1_N$ means that there is a~Hermitian metric on
${\mathcal{W}}^{*}\otimes_{\mathbb{C}} TN$.
On f\/ixing orientations, this may also be interpreted as an identif\/ication
${\mathcal{W}}\otimes_{\mathbb{C}}\overline{\mathcal{W}}= TN\otimes_{\mathbb{C}}\overline{TN}$, or
$\wedge^{1,1}{\mathcal{W}}=\wedge^{1,1}TN$.
Here the orientations will be chosen so that the induced almost complex structure $J$ on $TM=\pi^*({\mathcal{W}}\oplus
TN)$ is selfdual, i.e., the corresponding weightless $2$-form $\Omega\in\mathrm{C}^\infty(M,L^2 \wedge^2T^{*} M)$ is selfdual.

The Hodge star operator on $\wedge^2TM$ interchanges the vertical and horizontal components using the identif\/ication
above, and acts on mixed bivectors by $*(U\wedge X)=-JU\wedge JX$.
In particular the torsion of the modif\/ied connection ${\mathcal{D}}$ is selfdual if and only if
\begin{gather}
\varphi^{-1}\bigl(F^\alpha-[\varphi,\overline\varphi]_{\lie{h}}\bigr) =\psi,
\nonumber
\\
\varphi^{-1}\bigl(D^\alpha_X\varphi(U)+D^\alpha_{JX}\varphi(JU)\bigr)
=\ip{{\mathcal{C}}(U),X}+\ip{{\mathcal{C}}(JU),JX}
\label{eq:svdef}
\end{gather}
for any vector f\/ield $X$ and section $U$ of ${\mathcal{W}}$.
The f\/irst equation def\/ines $\psi$ uniquely.
For the second equation, note that conformal connections $D$ on ${\mathcal{W}}$ form an af\/f\/ine space modelled on $T^{*}
N\otimes_{\mathbb{R}}\mathop{\mathfrak{co}}({\mathcal{W}})$, and therefore $D$ can be chosen uniquely such that
$\varphi^{-1}D^\alpha_X\varphi$ is symmetric and traceless for all $X$.
The remaining ambiguity in ${\mathcal{C}}$ is f\/ixed by supposing it is complex linear, i.e., a~section of
$S^2_0{\mathcal{W}}^{*}\otimes_{\mathbb{C}} TN$, so that the second equation determines it uniquely.
Note though, that the second equation does not use all of $D$: only the induced holomorphic structure of ${\mathcal{W}}$
is needed.

For the selfduality of the Weyl curvature, it now suf\/f\/ices to compute $R^{{\mathcal{D}},ta}_{U,Y}$ and one additional
component.
To see this, note f\/irst that $\smash{\widehat\psi}$ and $\smash{\widehat{\mathcal{C}}}$ are
$\mathop{{\rm so}}(TM)$-valued, and therefore no information is lost by considering $R^{\mathcal{D}}$ and $R^D$ to
be the curvatures on $L^{-1}  TM$ rather than $TM$, so that the curvature equation is an identity
between $\mathop{{\rm so}}(TM)$-valued $2$-forms.
The selfduality condition is obtained by requiring that the Weyl part of $R^{{\mathcal{D}},ta}$ is selfdual.
This amounts to considering the traceless part of the induced bundle map from $\wedge^2_\asd T M$ to
$\mathop{{\rm so}}_\asd(TM)$.
Let $\overline\Omega$ and $\overline J$ be the sections of these bundles obtained by reversing the orientation of
$\Omega$ and $J$ on $TN$.
There are therefore three equations to f\/ind:
\begin{enumerate}\itemsep=0pt
\item[(i)]
The traceless part of the $\overline\Omega^\perp\otimes\overline J^\perp$ component of $R^{{\mathcal{D}},ta}$ should be
zero.
\item[(ii)]
The part of $R^{{\mathcal{D}},ta}$ in $\overline\Omega^\perp\otimes\langle\overline J\rangle$ should be zero.
\item[(iii)]
The multiple of $\overline\Omega\otimes\overline J$ should equal half the trace in
$\overline\Omega^\perp\otimes\overline J^\perp$.
\end{enumerate}
All of these except (iii) involve considering only the part of $R^{{\mathcal{D}},ta}$ in
$\overline\Omega^\perp\otimes\mathop{{\rm so}}(TM)$, which involves evaluating $R^{{\mathcal{D}},ta}$ on bivector
f\/ields of the form $U\wedge Y-{*(U\wedge Y)}=U\wedge Y+JU\wedge JY$.
An easy calculation gives
\begin{gather*}
R^{{\mathcal{D}},ta}_{U,Y}
=-D_Y{\mathcal{C}}(U,\cdot)+D_Y{\mathcal{C}}(U,\cdot)^{\scriptscriptstyle\mathrm T}-\tfrac14D_Y\psi(JU)(J-\overline J)
\\
\phantom{R^{{\mathcal{D}},ta}_{U,Y}=}
{}+{\mathcal{C}}\bigl(\ip{{\mathcal{C}}(U),Y},\cdot\bigr)
-{\mathcal{C}}\bigl(\ip{{\mathcal{C}}(U),Y},\cdot\bigr)^{\scriptscriptstyle\mathrm T}
-\tfrac14\ip{{\mathcal{C}}(U,J\psi),Y}(3J-\overline J)-\tfrac14J\psi(U)J\psi\mathinner{\vartriangle}Y,
\end{gather*}
where ${}^{\scriptscriptstyle\mathrm T}$ denotes the transpose, and I have contracted horizontal and vertical skew
endomorphisms using $X\mathinner{\vartriangle}Y=\frac12\Omega(X,Y)(J-\overline J)$ and $U\mathinner{\vartriangle}V
=\frac12\Omega(U,V)(J+\overline J)$.
Also, using $\Omega$, $\psi$ may be viewed as a~section of ${\mathcal{W}}^{*}$, so that $\psi(X,Y)=2\Omega(X,Y)J\psi$.

One readily obtains conditions (i) and (ii):{\samepage
\begin{gather*}
D_X{\mathcal{C}}(U,\cdot)+D_{JX}{\mathcal{C}}(JU,\cdot)=0,
\\
\tfrac12\bigl( D_X\psi(U)+D_{JX}\psi(JU)\bigr)=-3\ip{{\mathcal{C}}(\psi,U),X},
\end{gather*}
for any vector f\/ield $X$ and section $U$ of ${\mathcal{W}}$.}

It remains to compute condition (iii).
For this note that $R^D$ is given by $\frac14s_N(\Omega-\overline\Omega) \otimes(J-\overline
J)+\frac14s_{\mathcal{W}}(\Omega-\overline\Omega)\otimes(J+\overline J)$, where $s_N=\frac12\scal_N$ and
$s_{\mathcal{W}}=\frac12\scal_{\mathcal{W}}$ are the normalized scalar curvatures of $D$ on $TN$ and ${\mathcal{W}}$.
This yields, f\/inally,
\begin{gather*}
s_N-s_{\mathcal{W}}=|\psi|^2-2|{\mathcal{C}}|^2.
\end{gather*}

These equations, and the equations def\/ining ${\mathcal{C}}$ and $\psi$, depend only on $D$ through the induced
holomorphic structure $\overline\partial{}^a$ on ${\mathcal{W}}$.
The conformal connections on ${\mathcal{W}}$ and $TN$ may be viewed as Chern connections determined by the holomorphic
structures on these bundles together with the choice of a~Weyl derivative on $L^1$.
The dif\/ference $s_N-s_{\mathcal{W}}$ does not depend on the choice of Weyl derivative, since it is the normalized scalar
curvature of the Chern connection on the Hermitian line bundle ${\mathcal{W}}^{-1}  TN$.

\begin{aside}
The construction of a~``Chern--Weyl'' connection does not seem to be known and relates to a~simple coordinate-free
description (also not well-known) of the usual Chern connection, so I will sketch it here.
Let $E$ be a~complex vector bundle (over a~complex manifold $M,J$) with holomorphic structure $\overline\partial{}^E$
and compatible conformal metric $\mathsf{c}\colon S^2_{\mathbb{R}} E\to L^2_E$ with respect to which the complex
structure on $E$ is orthogonal; the latter is equivalently a~conformal Hermitian structure $\smash{\overline
E}\otimes_{\mathbb{C}} E\to L^2_E\otimes{\mathbb{C}}$.
(Note that a~complex line bundle automatically has a~conformal Hermitian structure.) Now given any covariant derivative
$D$ on $L^1_E$, let $\overline\partial{}^D_X=\frac12(D_X+{\boldsymbol i} D_{JX})$ (for all vector f\/ields $X$) be the
induced almost holomorphic structure on $L^2_E\otimes{\mathbb{C}}$.
Then there is a~unique conformal Hermitian connection $D^E$ on $E$ inducing $D$ on $L^1_E$.
It is given by the formula
\begin{gather*}
\ip{D^E_Xs_1,s_2}=\overline\partial{}^D_X\ip{s_1,s_2}-\ip{s_1,\overline\partial{}^E_X s_2}+
\ip{\overline\partial{}^E_Xs_1,s_2},
\end{gather*}
where $s_1$ and $s_2$ are sections of $E$.
The proof is immediate (the idea behind the formula is simply that $D^E=\partial^E+\overline\partial{}^E$, where
$\partial^E$ is the complex-linear part of the derivative).
In the case that the covariant derivative on $L^1_E$ is just a~trivialization, so that
$\overline\partial{}^D=\overline\partial$, this is the usual Chern connection by uniqueness.
\end{aside}

\begin{defn}
Suppose that ${\mathcal{W}}$ is a~complex line bundle over a~Riemann surface $N$ such that $L^1_{\mathcal{W}}=L^1_N$
(i.e., there is a~Hermitian metric on ${\mathcal{W}}^{-1}  TN$).
Equip ${\mathcal{W}}$ with a~holomorphic structure~$\overline\partial{}^a$, a~section ${\mathcal{C}}$ of
${\mathcal{W}}^{-2}  TN$, and a~section $\psi$ of ${\mathcal{W}}^{-1}$.
Then $(N,{\mathcal{W}})$ is said to be a~\emph{spinor-vortex space} if $(\overline\partial{}^a,{\mathcal{C}},\psi)$
satisfy the equations
\begin{gather}
\overline\partial{}^a{\mathcal{C}}=0,
\label{eq:sv1}
\\
\overline\partial{}^a \psi=-3{\mathcal{C}}\overline\psi,
\label{eq:sv2}
\\
s_{\smash{{\mathcal{W}}^{-1}  TN}} =\psi\overline\psi-2{\mathcal{C}}\overline{\mathcal{C}},
\label{eq:sv3}
\end{gather}
where $s_{\smash{{\mathcal{W}}^{-1}  TN}}$ is the normalized scalar curvature of the Chern connection
on ${\mathcal{W}}^{-1}  TN$.
\end{defn}
Joyce's lemma now gives the following result.
\begin{thm}
\label{th:sv}
Let $M$ be an oriented conformal manifold with a~$2$-dimensional Lie algebra $\lie{h}$ of linearly independent conformal
vector fields such that the projection $\pi$ onto the space of orbits is a~submersion over a~Riemann surface~$N$.
Let~${\mathcal{W}}$ be a~complex line bundle on $N$ such that~$\pi^*{\mathcal{W}}$ is the vertical bundle of~$M$
$($trivialized along the fibres by invariant vector fields$)$ and define
$(\overline\partial{}^a,{\mathcal{C}},\psi)$ by~\eqref{eq:svdef} as explained above.

Then $M$ is selfdual if and only if $(N,{\mathcal{E}})$ is a~spinor-vortex space.
\end{thm}

Of course there are only two possible $2$-dimensional Lie algebras.
In the case that $\lie{h}$ is Abelian and $\psi=0$, this result reduces to the original application of Lemma~\ref{lem:J}
by Joyce~\cite{Joy:esd}.

\begin{rem}
It is straightforward to adapt the calculations of this subsection to other signatures: in general $N$ is a~conformal
surface with two line bundles ${\mathcal{W}},\widetilde{\mathcal{W}}$ such that ${\mathcal{W}} \widetilde{\mathcal{W}}\cong T^{1,0}N  T^{0,1}N$ where $T^{1,0}N$ and $T^{0,1}N$ are the null line
subbundles of $TN\otimes{\mathbb{C}}$.
The line bundles ${\mathcal{W}}$ and $\widetilde{\mathcal{W}}$ are equipped with ``(anti)holomorphic structures''
$\tilde\partial^a,\partial^a$ and there are sections $({\mathcal{C}},\widetilde{\mathcal{C}},\psi,\widetilde\psi)$
satisfying appropriate reality conditions: in the Euclidean case $\widetilde{\mathcal{W}}=\overline{\mathcal{W}}$,
$\widetilde{\mathcal{C}}=\overline{\mathcal{C}}$ and $\widetilde\psi=\overline\psi$, but there is also a~Lorentzian
case, when the f\/ields are all real, and a~Euclidean reduction from Kleinian signature $(2,2)$ when
$\widetilde\psi=-\overline\psi$.
\end{rem}

\subsection{Reduction to three dimensions}

The reduction to three dimensions is the Jones--Tod correspondence~\cite{JoTo:mew}, which was one of the motivations for
this work.
The proof using Joyce's lemma is outlined in~\cite{Joy:esd}, so I will only recall the main ideas and def\/initions.

In this case, $Q$ is an oriented $3$-dimensional conformal manifold, which will be denoted $B$, and ${\mathcal{V}}$ is
an oriented line bundle isomorphic to $L^1=L^1_B$.
Note that ${\mathcal{C}}$ automatically vanishes, and write $\psi=*\omega$, where $\omega$ is a~$1$-form on $B$.
Then the torsion $T^{\mathcal{D}}$ will be selfdual provided that
\begin{gather}
*(D\varphi-\omega\varphi)=F^\alpha.
\end{gather}
(The connection $\alpha$ does not act on $\varphi$ since $\lie{h}$ is Abelian.)

This equation determines $(D,\omega)$ up to the gauge freedom $D\mapsto D+\gamma,\omega\mapsto\omega-\gamma$ for any
$1$-form $\gamma$ on $B$ (note that $(D+\gamma)\varphi=D\varphi-\gamma\varphi$, since $\varphi$ has weight $-1$).
In other words, it determines uniquely the Weyl derivative $D+\omega$.
This is called the \emph{Jones--Tod} Weyl structure.

The gauge freedom can be used to set $D\varphi=0$, or to set $\omega=0$.
Taking the latter point of view, the Jones--Tod Weyl structure is just $D$, determined by $* D\varphi=F^\alpha$.

The curvature of ${\mathcal{D}}$ is simply the pullback of the curvature of $D$ on $B$.
In particular the Weyl part will be selfdual if\/f it vanishes and this is readily seen to be equivalent to the vanishing
of~$r^D_0$~\cite{Joy:esd}.

\begin{defn}
Suppose that $B$ is a~conformal $3$-manifold and let $D$ be a~Weyl connection on~$B$.
Then $B$ is said to be an \emph{Einstein--Weyl space} if\/f $r^D_0=0$, i.e., the symmetric traceless Ricci tensor of~$D$
vanishes.
\end{defn}

\begin{thm}
[\cite{JoTo:mew}]
\label{th:JnT}
Let $M$ be an oriented conformal $4$-manifold and $K$ a~nonvanishing conformal vector field such that the projection
$\pi$ of $M$ onto the space of trajectories is a~submersion over a~conformal manifold $B$.
Equip $B$ with the Jones--Tod Weyl structure.

Then $M$ is selfdual if and only if $B$ is Einstein--Weyl.
\end{thm}

All three background geometries are themselves def\/ined by geometric gauge theories: conformal local trivializations of
the bundle ${\mathcal{V}}$ over $Q$ are related by gauge transformations.
In the three-dimensional case, ${\mathcal{V}}$ is simply $L^1$, so conformal trivializations are length scales, and this
is Weyl's original gauge theory~\cite{Wey:stm}.
Analogously, in the one and two-dimensional case, ${\mathcal{E}}$~and~${\mathcal{W}}$ should be regarded as part of the
geometry of the space, \emph{not} auxiliary bundles, and these background equations, although gauge-theoretic, should
not be confused with the gauge f\/ield equations on auxiliary $G$-bundles which will be studied in Section~\ref{s:gfe}.

\subsection*{Addenda: null reductions}

The constructions of this section assume nondegeneracy of the conformal structure $\mathsf{c}$ on $M$ along the orbits
of the symmetry group $H$.
In (real) Euclidean signature, this is automatic, but this is not the case when $M$ has Kleinian signature, nor when
$M$ is a~holomorphic conformal manifold.
When the conformal structure degenerates on the orbits, the reduction is said to be \emph{null}.
For the generic local considerations of this paper, I assume that the radical (or kernel) $R_V:=VM\cap VM^\perp$ of the
conformal structure along the $H$ orbits has constant rank (as before, $VM$ denotes the tangent bundle to the $H$
orbits).
There are thus three possibilities:
\begin{enumerate}\itemsep=0pt
\item[(1)] $R_V$ has rank one;
\item[(2$+$)] $R_V=U_+$ has rank two and is selfdual;
\item[(2$-$)] $R_V=U_-$ has rank
two and is antiselfdual.
\end{enumerate}
In case (1), $R_V^\perp = VM+VM^\perp$ has rank three, and $R_V^\perp/R_V$ is the sum of two null subbundles with
a~nondegenerate pairing between them.
Their inverse images $U_\pm$ in $R_V^\perp$ are selfdual and antiselfdual null 2-plane distributions.
There are several subcases to consider here, depending on the rank of $VM$.
If either $U_\pm \subseteq VM$ then $VM$ must have rank $3$, in which case $VM^\perp\subset VM$, hence $VM^\perp=R_V$
and $R_V^\perp=VM$, so $VM=U_+ + U_-$.
If $VM$ has rank two, then $TM=VM+U_++U_-$, and clearly if $\rank V_M=1$, then $R_V=VM$.

Cases (2$\pm$) are simpler, at least when $\rank VM=2$, so $VM=R_V$ is a~bundle of selfdual or antiselfdual $2$-planes.
I concentrate on these two cases here.
One justif\/ication for such a~focus is that in the other cases, one may be able to prove that the distributions $U_+$ or
$U_-$ are integrable, and hence study the reduction in terms of their leaf spaces.
For instance, Dunajski--West~\cite{DuWe:acs} establish such an integrability result when $\rank VM=1$.

Suppose then that $\pi\colon M\to Q=M/H$ is a~principal $H$-bundle over a~manifold $Q$ such that $VM$ has rank two and
is totally null.
As $VM^\perp=VM$, neither $VM$ nor $TQ$ inherit conformal structures from $M$; instead there is a~nondegenerate pairing
$VM\times\pi^* TQ\to L^2$.
This pairing in $H$-invariant, and thus identif\/ies $VM=\pi^* (T^*Q\otimes {\mathcal{V}})$, for a~line bundle
${\mathcal{V}}\to Q$ (whose pullback to $M$ will be identif\/ied with $L^2$).
As in the nondegenerate setting the isomorphism of $VM$ with $M\times\lie{h}$ descends to a~bundle isomorphism
$\varphi\colon T^*Q\otimes {\mathcal{V}} \to \lie{h}_Q:= M\times_H\lie{h}$ over $Q$.
There is no canonical splitting of the short exact sequence $0\to \pi^*T^*Q\otimes {\mathcal{V}} \to TM \to \pi^*TQ\to
0$, but such a~splitting may be chosen so that the complementary subbundle to $VM$ is $H$-invariant and null.
This yields a~principal $H$-connection $\alpha$ on $\pi\colon M\to Q$.

A torsion-free connection on $TN$ and a~connection on ${\mathcal{V}}$ together induce a~conformal connection $D$ on $TM$
with vertical-valued torsion as before.
A conformal connection ${\mathcal{D}}$ with selfdual torsion may be obtained by adding correction terms to $D$.
Compared to the nondegenerate case, the correction term ${\mathcal{C}}$ may be absorbed into the choice of $D$, and
replaced by a~vertical correction $\chi\colon \wedge^2(T^*Q\otimes{\mathcal{V}})\to T^*Q\otimes {\mathcal{V}}$.
Thus ${\mathcal{D}}=D+\widehat\chi+\widehat\psi$, where $\psi\colon \wedge^2TQ\to T^*Q\otimes {\mathcal{V}}$ as
before, and
\begin{gather*}
2\widehat\chi_{U+X} (V+Y)= \ip{\chi(U,\cdot),Y}+\ip{\chi(V,\cdot),X}-\chi(U,V),
\\
2\widehat\psi_{U+X} (V+Y)= \ip{\psi(X,\cdot),Y}+\ip{\psi(Y,\cdot),X}-\psi(X,Y).
\end{gather*}
The torsion satisf\/ies
\begin{gather}
\label{eq:null-t1}
T^{\mathcal{D}}(U,V)=\varphi^{-1}[\varphi(U),\varphi(V)]_{\lie{h}}-\chi(U,V),
\\
T^{\mathcal{D}}(U,X)=\varphi^{-1}D^\alpha_X\varphi(U),
\\
T^{\mathcal{D}}(X,Y)=\varphi^{-1}F^\alpha(X,Y)-\psi(X,Y),
\label{eq:null-t3}
\end{gather}
while the curvature computes to
\begin{gather*}
R^{{\mathcal{D}}}_{U+X,V+Y}=R^D_{U+X,V+Y}+(\widehat\chi+\widehat\psi)_{T^{\mathcal{D}}(U+X,V+Y)}
\\
\phantom{R^{{\mathcal{D}}}_{U+X,V+Y}=}
{}+D_{U+X}(\widehat\chi+\widehat\psi)_{V+Y} -D_{V+Y}(\widehat\chi+\widehat\psi)_{U+X}
\\
\phantom{R^{{\mathcal{D}}}_{U+X,V+Y}=}
{}+(\widehat\chi+\widehat\psi)_{\psi(X,Y)}+ (\widehat\chi+\widehat\psi)_{\chi(U,V)}
+\Abrack{(\widehat\chi+\widehat\psi)_{U+X}, (\widehat\chi+\widehat\psi)_{V+Y}}.
\end{gather*}
Hence the torsion-adjusted curvature is
\begin{gather}
R^{{\mathcal{D}},ta}_{U+X,V+Y}=R^D_{X,Y}+\Abrack{(\widehat\chi+\widehat\psi)_X,
(\widehat\chi+\widehat\psi)_Y}+\widehat\chi_{\psi(X,Y)}
\nonumber
\\
\phantom{R^{{\mathcal{D}},ta}_{U+X,V+Y}=}
{}+D_X(\widehat\chi+\widehat\psi)_Y -D_Y(\widehat\chi+\widehat\psi)_X +D_X\widehat\chi{}_V-D_Y\widehat\chi{}_U
\nonumber
\\
\phantom{R^{{\mathcal{D}},ta}_{U+X,V+Y}=}
{}+\Abrack{(\widehat\chi+\widehat\psi)_X,\widehat\chi{}_V} +\Abrack{\widehat\chi{}_U,(\widehat\chi+\widehat\psi)_Y}
+\abrack{\widehat\chi_U,\widehat\chi_V} +\widehat\chi_{\chi(U,V)}.
\label{eq:null-c}
\end{gather}

\subsection*{Antiselfdual ($\boldsymbol{\alpha}$-surface) reduction}

When $VM$ is antiselfdual, the f\/ibres of $\pi\colon M\to Q$ are $\alpha$-surfaces (i.e., they correspond to points in
the twistor space of $M$).
Under the isomorphism of $TM$ with $\pi^*(T^*Q\otimes{\mathcal{V}} \oplus TQ)$, the antiselfdual bivectors decompose
into three rank one subbundles: $\wedge^2(T^*Q\otimes{\mathcal{V}})$, $\wedge^2TQ$, and the tracelike part of
$TQ\otimes T^*Q\otimes{\mathcal{V}}$.
The corresponding decomposition of $\mathop{{\rm so}}_-(TM)$ has summands $\Hom_-(TQ,T^*Q\otimes{\mathcal{V}})$,
$\Hom_-(T^*Q\otimes{\mathcal{V}},TQ)$ and the span of $\iden_{T^*Q\otimes{\mathcal{V}}} - \iden_{TQ}$; here $\Hom_-$
denotes the subbundle of skew symmetric operators.

In order to interpret equations~\eqref{eq:null-t1}--\eqref{eq:null-t3} and~\eqref{eq:null-c}, it is convenient to set
${\mathcal{V}}={\mathcal{L}}^2\otimes \wedge^2TQ$ for a~line bundle ${\mathcal{L}}$, so that $\varphi$, $\psi$ and
$\chi$ may be viewed as $1$-forms on $Q$, with values in ${\mathcal{L}}^2\otimes\lie{h}_Q$, ${\mathcal{L}}^{-2}$ and
${\mathcal{L}}^2$ respectively.

Then $T^{\mathcal{D}}$ is selfdual if\/f $\chi\wedge\varphi = \frac12 [\varphi\wedge\varphi]_{\lie{h}}$,
$\psi\wedge\varphi=F^\alpha$ and ${\rm d}^{D,\alpha}\varphi=0$.
The f\/irst two equations determine $\chi$ and $\psi$ uniquely (since $\varphi\colon TN\to {\mathcal{L}}\otimes\lie{h}_Q$
is injective), while the third depends only on (and essentially determines) the connection $a$ induced by $D$ on
${\mathcal{L}}$.

In the expression~\eqref{eq:null-c} for the torsion-adjusted curvature, very few terms contribute to the antiselfdual
Weyl part.
For instance, only the trace part of f\/irst and last lines contribute, and the latter trace vanishes identically.
After some tedious computations, the background equations for the $1$-forms $\psi$ and $\chi$ and the connection $a$ on
${\mathcal{L}}$ reduce to
\begin{gather*}
{\rm d}^a\psi=0,
\qquad
{\rm d}^a\chi=0,
\qquad
F^a=\chi\wedge\psi,
\end{gather*}
which may be interpreted as the f\/latness of the connection ${\rm d}^a+\psi+\chi$ on ${\mathcal{L}}\oplus
{\mathcal{L}}^{-1}$ (where ${\rm d}^a$ is the direct sum connection, while $\psi$ and $\chi$ are viewed as
$1$-forms with values in $\Hom({\mathcal{L}}^{-1},{\mathcal{L}})$ and $\Hom({\mathcal{L}},{\mathcal{L}}^{-1})$).

\subsection*{Selfdual ($\boldsymbol{\beta}$-surface) reduction}

In the selfdual case, only the mixed part of $T^{\mathcal{D}}$ has a~antiselfdual component, so there is no loss in
setting $\psi=\chi=0$.
The selfduality of the torsion then reduces to $D^\alpha\varphi=\frac12{\rm d}^{D^\alpha}\varphi$, where
$\varphi$ is interpreted as a~$1$-form on $Q$ as before.
This equation determines $D$ up to projective transformation ($D_XY\mapsto D_XY +\gamma(X)Y+\gamma(Y)X$ for a~$1$-form
$\gamma$), and the background equations are vacuous.
Hence the background geometry is an arbitrary (torsion-free) projective surface $(Q,[D])$; this reduction was obtained
in~\cite{DuWe:acs} and~\cite{Cal:sdp}.

\section{Interlude: Bianchi metrics}
\label{s:bm}

Selfdual conformal manifolds with a~freely acting three-dimensional symmetry group have been studied in many
places~\cite{Dan:sfk,DaSt:cok,Hit:tem, Maz:csb,MMW:sbm,PePo:kzs,Tod:p6,Tod:com,Tod:p3}.
In this interlude, I will show brief\/ly how the direct approach to selfdual Bianchi metrics is related to the Riccati
space reduction of the previous section.

For simplicity, I focus on the case of diagonal Bianchi IX metrics.
Such a~metric may be written in the form
\begin{gather*}
g=w_1w_2w_3 {\rm d} t^2 +\frac{w_2w_3}{w_1}\sigma_1^{2}+\frac{w_3w_1}{w_2}\sigma_2^{2}
+\frac{w_1w_2}{w_3}\sigma_3^{2},
\end{gather*}
where $\sigma_1$, $\sigma_2$, $\sigma_3$ are the usual left-invariant $1$-forms on $\SU(2)$ and $w_1$, $w_2$, $w_3$ are functions of~$t$.
If $X_1$, $X_2$, $X_3$ are the dual vector f\/ields to $\sigma_1$, $\sigma_2$, $\sigma_3$, then the vector f\/ields
$\partial_t$, $\varphi_1=w_1 X_1$, $\varphi_2=w_2 X_2$ and $\varphi_3=w_3 X_3$ form a~conformal frame.
Notice that $\dot\varphi_1 -[\varphi_2,\varphi_3]$ is the multiple $\dot w_1-w_2w_3$ of $X_1$, where the dot denotes the
derivative with respect to $t$.
The other two components are similar.
Following~\cite{PePo:kzs,Tod:p6}, write
\begin{gather*}
\dot w_1=w_2w_3-w_1(A_2+A_3),
\\
\dot w_2=w_3w_1-w_2(A_3+A_1),
\\
\dot w_3=w_1w_2-w_3(A_1+A_2).
\end{gather*}
Comparing with equation~\eqref{eq:riccdef}, observe that the matrix ${\mathcal{B}}$ is diagonal, with eigenvalues
$-(A_2+A_3),-(A_3+A_1),-(A_1+A_2)$.
In the approach of the previous section, the conformal gauge freedom is used to set $A_1+A_2+A_3=0$, so that
${\mathcal{B}}$ is traceless with eigenvalues $A_1$, $A_2$, $A_3$.
This is not usually done in the literature, because by working with an arbitrary compatible metric, additional
(non-conformally-invariant) equations can be imposed.
In particular, vacuum metrics~-- and more generally, K\"ahler  metrics (with antiselfdual complex structure)~-- are
scalar-f\/lat.
Hence scalar-f\/latness is often used as a~gauge condition, in which case the following well-known system, originating in
work of Brioschi, Chazy, Darboux, and Halphen, is obtained:
\begin{gather*}
\dot A_1=A_2A_3-A_1(A_2+A_3),
\\
\dot A_2=A_3A_1-A_2(A_3+A_1),
\\
\dot A_3=A_1A_2-A_3(A_1+A_2).
\end{gather*}
Joyce's lemma explains the remarkable fact that this system depends only on $w_1,w_2,w_3$ through the functions
$A_1,A_2,A_3$.
The trace of this system is the scalar-f\/lat gauge condition: $\dot A_1+\dot A_2+\dot A_3=-A_2A_3-A_3A_1-A_1A_2$, while
the traceless part is the selfduality equation for the Weyl curvature~-- the latter condition is independent of the
choice of conformal gauge, and may be rewritten as a~matrix Riccati equation in the following way:
\begin{gather*}
{\mathcal{B}} =\frac13\left[
\begin{matrix}
2A_1-A_2-A_3&0&0
\\
0&2A_2-A_3-A_1&0
\\
0&0&2A_3-A_1-A_2
\end{matrix}
\right],
\\
\dot{\mathcal{B}}-a{\mathcal{B}}=2\big(\mathcal{B}^2\big)_0,
\end{gather*}
{where}
\begin{gather*}
a=-\frac23(A_1+A_2+A_3).
\end{gather*}
In other words, the matrix Riccati equation is obtained by separating the diagonal matrix ${\mathcal{B}}+a\iden$, with
eigenvalues $-(A_2+A_3),-(A_3+A_1),-(A_1+A_2)$, into its trace (which def\/ines the connection $D=\partial_t+a$ on the one-dimensional quotient geometry), and its traceless part (which is ${\mathcal{B}}$).

One advantage of this geometric interpretation is that dif\/ferent gauge conditions can be easily compared.
Evidently the condition $a=0$ f\/ixes $t$ up to af\/f\/ine transformations.
In order to interpret the scalar-f\/lat gauge condition $\dot a=\frac23(A_2A_3+A_3A_1+A_1A_2)$, it is natural and
illuminating to express the right hand side in terms of $a$ and $\trace{\mathcal{B}}^2$.
The result is
\begin{gather*}
\dot a-\tfrac12a^2=-\tfrac13\trace{\mathcal{B}}^2.
\end{gather*}
Note that $(\partial_t-\frac12 a)(\partial_t+\frac12 a)=\partial_t^2+\frac12\dot a- \frac14 a^2$, so that on sections of
$L^{1/2}$, $D^2=\partial_t^2-\frac16\trace{\mathcal{B}}^2$.
Hence the scalar-f\/lat gauge condition may be interpreted as f\/ixing $t$ to be a~projective coordinate with respect to the
projective structure $D^2+\frac16\trace{\mathcal{B}}^2$; this determines $t$ up to a~projective transformation.

\section{The gauge f\/ield equations}
\label{s:gfe}

Let $\nabla$ be a~$G$-connection, on a~vector bundle $E$ over a~conformal manifold $M$, which is invariant under an
action of a~group $H$ of conformal transformations.
Inf\/initesimally, the generators form a~Lie algebra $\lie{h}$ of conformal vector f\/ields on $M$ and there is an action of
these vector f\/ields by Lie derivative on sections of $E$ such that ${\mathcal{L}}_\xi\nabla=0$ for all $\xi\in\lie{h}$.
Each $\xi\in\lie{h}$ therefore determines a~Higgs f\/ield $S_\xi$ in the associated Lie algebra bundle
$\lie{g}_M\leqslant\End(E)$ def\/ined by $S_\xi s={\mathcal{L}}_\xi s-\nabla_\xi s$.
(Note the unusual sign convention, which is necessary for consistency later.) Since ${\mathcal{L}}_\xi\nabla=0$ and
$[{\mathcal{L}}_\xi,{\mathcal{L}}_\chi]={\mathcal{L}}_{[\xi,\chi]}$, it follows that
\begin{gather*}
{\mathcal{L}}_\xi (S_\chi) = S_{[\xi,\chi]}
\end{gather*}
and so $S\colon M\mathbin{{\times}}\lie{h}\to\lie{g}_M$ is $H$-equivariant.
Part of the curvature $F^\nabla$ of $\nabla$ is determined by the Higgs f\/ields~-- one readily computes that for any
$\xi\in\lie{h}$ and any vector f\/ield $X$,
\begin{gather*}
F^\nabla(\xi,X)=\nabla_X(S_\xi).
\end{gather*}
In particular, for two vector f\/ields $\xi$, $\chi$ in $\lie{h}$,
\begin{gather*}
F^\nabla(\xi,\chi)=[S_\xi,S_\chi]_{\lie{g}}-S_{[\xi,\chi]},
\end{gather*}
where the f\/irst bracket is the Lie bracket in $\lie{g}_M$.

Now suppose that $H$ acts freely on $M$ with nondegenerate conformal metrics on the orbits, and def\/ine
$(\alpha,\varphi)$ as in Section~\ref{s:sdbg}.
The bundle $E$ is the pullback of a~bundle, also denoted~$E$, over~$Q$.
Since $\nabla$ is $H$-invariant, it descends to a~connection~$A$ over~$Q$.
The Higgs f\/ields are also $H$-invariant; hence setting $\Phi(U)=S_{\varphi(U)}$ def\/ines a~bundle map
$\Phi\colon{\mathcal{V}}\to\lie{g}_{\smash Q}$ over $Q$.
Next introduce a~conformal connection $D$ on ${\mathcal{V}}$.
Then $\nabla_X(S_\xi)=(D^A_X\Phi)\bigl(\varphi^{-1}(\xi)\bigr)
-\Phi\bigl(\varphi^{-1}(D^\alpha_X\varphi)\varphi^{-1}(\xi)\bigr)$ and a~simple computation of the curvature of $\nabla$
yields the following equations:
\begin{gather*}
F^\nabla(U,V)=[\Phi(U),\Phi(V)]_{\lie{g}} -\Phi\bigl(\varphi^{-1}[\varphi(U),\varphi(V)]_{\lie{h}}\bigr),
\\
F^\nabla(U,X)=(D^A_X\Phi)(U) -\Phi\bigl(\varphi^{-1}D^\alpha_X\varphi(U)\bigr),
\\
F^\nabla(X,Y)=F^A(X,Y)-\Phi\bigl(\varphi^{-1}F^\alpha(X,Y)\bigr).
\end{gather*}
This formulation makes manifest the analogy between $(A,\Phi)$ and $(\alpha,\varphi)$: the former is a~$G$-connection
and $\lie{g}_Q$-valued section of ${\mathcal{V}}^{*}$, while the latter is an $H$-connection and $\lie{h}_Q$-valued
section of ${\mathcal{V}}^{*}$.

Furthermore, these formulae for $F^\nabla$ are closely analogous to the formulae~\eqref{eq:t1}--\eqref{eq:t3} for
$T^{\mathcal{D}}$ obtained in the previous section: adding torsion terms to the above equations yields
\begin{gather}
\label{eq:gfe1}
\big(F^\nabla+\Phi\circ T^{\mathcal{D}}\big)(U,V)=[\Phi(U),\Phi(V)]_{\lie{g}},
\\
\big(F^\nabla+\Phi\circ T^{\mathcal{D}}\big)(U,X)=(D^A_X\Phi)(U) -\Phi\bigl(\ip{{\mathcal{C}}(U),X}\bigr),
\\
\big(F^\nabla+\Phi\circ T^{\mathcal{D}}\big)(X,Y)=F^A(X,Y)-\Phi\bigl(\psi(X,Y)\bigr).
\label{eq:gfe3}
\end{gather}
Assuming that $T^{\mathcal{D}}$ is selfdual (which can always be arranged, using the choice of $D$, by the work of
Section~\ref{s:sdbg}), the selfduality of $F^\nabla$ is now equivalent to the selfduality of the right hand sides
of~\eqref{eq:gfe1}--\eqref{eq:gfe3}.

This reduced gauge f\/ield equation is easy to compute in each dimension (and equivalent calculations have already been
used in Section~\ref{s:sdbg}).
Notice that no assumption of selfduality on~$M$ is needed for these computations: just as the selfdual Yang--Mills
equation makes sense on any oriented conformal $4$-manifold, so also the generalized Bogomolny equation is def\/ined on
any Weyl space, and the same is true on the one and two-dimensional geometries.
However, the principal dogma underlying this work is that the natural backgrounds for the gauge f\/ield equations are
selfdual spaces, Einstein--Weyl spaces, spinor-vortex spaces and Riccati spaces.
There are two reasons for this: f\/irst, the Ward correspondence predicts that the reduced gauge f\/ield equations on these
backgrounds will be integrable; second, it will soon be apparent, if not already, that the gauge f\/ield equations and
background equations are intimately related.

\subsection{Generalized Nahm equations on Riccati spaces}

The reduction to one dimension is straightforward, since there is no curvature.
Hence the connection $A$ on $\lie{g}_C$ can be assumed trivial, and the gauge f\/ield equation for
$\Phi\in\mathrm{C}^\infty(C,{\mathcal{E}}^{*}\otimes\lie{g}_C)$ is
\begin{gather}
\label{eq:nahm}
D\Phi-{*[\Phi,\Phi]_{\lie{g}}} ={\mathcal{B}}\mathinner{\cdot}\Phi,
\end{gather}
where ${\mathcal{B}}\mathinner{\cdot}\Phi=\Phi\circ{\mathcal{B}}$ and the Lie bracket pairing
$[\Phi,\Phi]_{\lie{g}}\in\mathrm{C}^\infty(C,\wedge^2{\mathcal{E}}^{*}\otimes\lie{g}_C)$ is interpreted as a~section
of $T^{*} C\otimes{\mathcal{E}}^{*}\otimes\lie{g}_C$ using the Hodge star and conformal structure on ${\mathcal{E}}$,
together with the identif\/ications $L^1_{\mathcal{E}}=L^1_C=TC$.
When ${\mathcal{B}}=0$, this is the Nahm equation.

\subsection{Generalized Hitchin equations on spinor-vortex spaces}

For the reduction to two dimensions, it is natural, as before, to reinterpret $\Phi$ as a~complex linear map from
${\mathcal{W}}$ to $\lie{g}_N\otimes{\mathbb{C}}$.
Then the gauge f\/ield equations for $(A,\Phi)$ are
\begin{gather*}
F^A-[\Phi,\overline\Phi]_{\lie{g}} =\psi\wedge\overline\Phi+\overline\psi\wedge\Phi,
\qquad
\overline\partial{}^{a,A}\Phi={\mathcal{C}}\overline\Phi.
\end{gather*}
When ${\mathcal{C}}=\psi=0$ and ${\mathcal{W}}=TN$, with the induced holomorphic structure, these are Hitchin's
equations.

To adapt the gauge f\/ield equations to spinor-vortex spaces in general signature (when there are background f\/ields
${\mathcal{C}},\widetilde{\mathcal{C}},\psi,\widetilde\psi$), replace $\overline\Phi$ by an additional f\/ield
$\widetilde\Phi$, satisfying the analogous $\partial^{a,A}$-equation.
The Lorentzian reality condition ($\Phi,\widetilde\Phi$ real) provides a~generalized chiral model, while
a~generalization of the harmonic map equation is obtained by introducing the crucial sign change
$\widetilde\Phi=-\overline\Phi$ (recall also that $\widetilde\psi=-\overline\psi$ in this case).

\subsection{Generalized Bogomolny equations on Einstein--Weyl spaces}

The reduction to three dimensions gives the natural generalization to Weyl geometry of the Bogomolny equation for
magnetic monopoles
\begin{gather*}
{* D^A\Phi}=F^A.
\end{gather*}
The Euclidean, hyperbolic or spherical Bogomolny equation arises when the Einstein--Weyl structure is given by a~metric
of constant curvature.

\subsection*{Addenda: twisted f\/lat pencils and projective pairs}

The same methodology as in the nondegenerate cases yields gauge f\/ield equations over the null reductions, using the
equations~\eqref{eq:null-t1}--\eqref{eq:null-t3} for the torsion $T^{\mathcal{D}}$.

In the background geometry $(a,\phi,\psi)$ on $(Q,{\mathcal{L}})$ obtained by reduction along an $\alpha$-surface
foliation, the gauge f\/ields consist of a~$G$-connection $A$ and a~$1$-form $\Phi$ with values in
${\mathcal{L}}^2\otimes\lie{g}_Q$, and the gauge f\/ield equations are
\begin{gather}
\label{eq:tfc}
F^A=\psi\wedge\Phi,
\qquad
{\rm d}^{a,A}\Phi=0,
\qquad
\tfrac12[\Phi\wedge\Phi]=\chi\wedge\Phi.
\end{gather}
If $\psi=\chi=0$ (so $a$ is f\/lat, and ${\mathcal{L}}$ may be trivialized) then these are the equations for a~pencil of
f\/lat connections ${\rm d}^A + \lambda\Phi$ on a~surface (also known as the topological chiral model).
Thus solutions of~\eqref{eq:tfc} may be called \emph{twisted flat pencils}.
The special case $\psi=0$ has been studied by Tafel and W{\'o}jcik in~\cite{TaWo:nkv}.

On a~projective surface $(Q,[D])$ obtained by reduction along a~$\beta$-surface foliation, the gauge f\/ields consist of
a~$G$-connection $A$ and a~$1$-form $\Phi$ with values in ${\mathcal{O}}_Q(2)\otimes\lie{g}_Q$ (where
${\mathcal{O}}_Q(3)=\wedge^2TQ$), and the gauge f\/ield equations are
\begin{gather*}
D^A\Phi=\tfrac12 {\rm d}^{D,A}\Phi.
\end{gather*}
Solutions $(A,\Phi)$ were referred to (somewhat unimaginatively) as \emph{projective pairs} in~\cite{Cal:sdp}.

The special case of reductions of the selfdual Yang--Mills equation on ${\mathbb{R}}^{2,2}$ (or ${\mathbb{C}}^4$) by null
translations was considered already by Mason and Woodhouse~\cite{MaWo:ist}: the $\alpha$-plane reduction (yielding f\/lat
pencils) is denoted $H_{SD}$, while the $\beta$-plane reduction (yielding projective pairs) is denoted~$H_{ASD}$.

\section{Interlude: spinor-vortex spaces and Hitchin's equations}
\label{s:svhe}

In~\cite{Hit:sde}, Hitchin considered solutions of the selfdual Yang--Mills equation on ${\mathbb{R}}^4$ invariant under
two translations.
He observed that the Yang--Mills connection could be decomposed into a~connection over ${\mathbb{R}}^2$ and two Higgs
f\/ields.
He combined these Higgs f\/ields into a~complex Higgs f\/ield and then noticed that, remarkably, the reduced Yang--Mills
equation becomes conformally invariant provided this complex Higgs f\/ield is interpreted as a~$1$-form rather than
a~scalar.
This was a~surprise, because although the selfdual Yang--Mills equation is conformally invariant, the notion of
translation-invariance is not.
Furthermore, conformal invariance in two dimensions implies an inf\/inite-dimensional symmetry group.

The reduction process described here provides a~simple explanation of this phenomenon: $W^\asd=0$ is conformally
invariant, and so is the notion of torus symmetry.
Hence the equations for $(\overline\partial{}^a,\psi,{\mathcal{C}},\alpha,\varphi)$ are conformally invariant on a~f\/ixed
Riemann surface with a~complex line bundle~${\mathcal{W}}$ and a~Hermitian metric on ${\mathcal{W}}^{-1}  TN$.
In particular, it is clear that the equations~\eqref{eq:sv1}--\eqref{eq:sv3} for
$(\overline\partial{}^a,\psi,{\mathcal{C}})$ are conformally invariant.

If ${\mathcal{C}}$ is not identically zero, then (on the open set where ${\mathcal{C}}$ is nonvanishing) the freedom in
the holomorphic structure on ${\mathcal{W}}$ can be f\/ixed by declaring that ${\mathcal{C}}$ is an identif\/ication of~${\mathcal{W}}^2$ with~$TN$.
The length of ${\mathcal{C}}$ now def\/ines a~gauge, breaking conformal invariance.
More precisely, given a~Weyl derivative $D$, the Chern connection on ${\mathcal{W}}$ only agrees with the connection
induced from $TN$ if $D=D^g$, where $g$ is the metric induced by $|{\mathcal{C}}|$.
Equation~\eqref{eq:sv3} becomes $\frac12s^g=-2+4|\psi|^2_g$; in particular, if $\psi$ is zero, then $g$ has constant
negative curvature, which is the case studied by Joyce~\cite{Joy:esd}, the \emph{hyperbolic spinor-vortex geometry}.

On the other hand if ${\mathcal{C}}$ is identically zero, then $\psi$ is holomorphic, so if $\psi$ is not identically
zero, then (on the open set where $\psi$ is nonvanishing), $\psi$ trivializes ${\mathcal{W}}$.
Again this breaks conformal invariance by introducing a~natural gauge, $|\psi|$; the corresponding metric has constant
positive curvature ($s^g=4$).
This is the \emph{spherical spinor-vortex geometry}.

Finally if ${\mathcal{C}}$ and $\psi$ both vanish, then the holomorphic structure can be f\/ixed by setting
${\mathcal{W}}=TN$, so that the Chern connection on ${\mathcal{W}}^{-1}  TN$ is trivial.
This trivial spinor-vortex geometry does not break conformal invariance.

The Yang--Mills equation is reduced to two dimensions by interpreting the Higgs f\/ields as a~section $\Phi$ of
${\mathcal{W}}^{-1}\otimes_{\mathbb{R}}\End(V)$, using $\varphi$.
The resulting gauge f\/ield equations are independent of $(\alpha,\varphi)$, i.e., they are intrinsic to the spinor-vortex
space.
On a~trivial spinor-vortex space, ${\mathcal{W}}=TN$, so $\Phi$ becomes an endomorphism-valued $1$-form, and, as
remarked already, the gauge f\/ield equations for $(A,\Phi)$ are Hitchin's equations.
Thus it is the isomorphism $\varphi$, and the geometry of the trivial spinor-vortex space that are responsible for the
interpretation of the Higgs f\/ields as a~$1$-form, rather than as scalars.

\section{Selfdual spaces from gauge f\/ields}
\label{s:sgf}

In Sections~\ref{s:sdbg} and~\ref{s:gfe}, background geometries and gauge f\/ield equations were found by reducing the
selfduality equation for conformal structures and Yang--Mills f\/ields respectively.
In the approach taken, the symmetry group $H$, and the pair $(\alpha,\varphi)$ def\/ining the reduction, decoupled from
the construction, and the reduced background geometry and gauge f\/ield equations were found to be independent of these
data.

In some sense, this was a~sleight of hand, since $(\alpha,\varphi)$ were used to def\/ine the data on the quotient which
make up the background geometry (i.e., $D$, $\psi$ and ${\mathcal{C}}$).
However, the procedure may be turned around: starting with the background geometry, these def\/initions, originally
obtained by imposing selfduality on the torsion $T^{\mathcal{D}}$, may be viewed as equations for the pair
$(\alpha,\varphi)$.
Indeed, applying $\varphi$ to the formulae~\eqref{eq:t1}--\eqref{eq:t3} for $T^{\mathcal{D}}$ gives:
\begin{gather}
\label{eq:tg1}
\big(\varphi\circ T^{\mathcal{D}}\big)(U,V)=[\varphi(U),\varphi(V)]_{\lie{h}},
\\
\big(\varphi\circ T^{\mathcal{D}}\big)(U,X)=(D^\alpha_X\varphi)(U) -\varphi\bigl(\ip{{\mathcal{C}}(U),X}\bigr),
\\
\big(\varphi\circ T^{\mathcal{D}}\big)(X,Y)=F^\alpha(X,Y)-\varphi\bigl(\psi(X,Y)\bigr).
\label{eq:tg3}
\end{gather}
The right hand sides of these formulae correspond precisely to the right hand sides of~\eqref{eq:gfe1}--\eqref{eq:gfe3}.
Hence the selfduality equation for $T^{\mathcal{D}}$ coincides with the gauge f\/ield equation for $(\alpha,\varphi)$ on
the background geometry.

This leads immediately to an inverse construction of selfdual spaces with symmetry from gauge f\/ields with
$\ell$-dimensional gauge group on $k$-dimensional background geometries, where $k+\ell=4$.
However, the construction can be generalized further by noting that the selfduality of $T^{\mathcal{D}}$ is implied by
the gauge f\/ield equation as long as $\varphi\colon{\mathcal{V}}\to\lie{h}_Q$ is injective.
I will now explain this generalized construction.

Let $P\to Q$ be a~principal $H$-bundle with an $H$-connection $\alpha$ and a~Higgs f\/ield
$\varphi\colon{\mathcal{V}}\to\lie{h}_{\smash{Q}}$, where $\lie{h}_{\smash{Q}}=P\mathbin{{\times}}_H\lie{h}$.
Suppose that $H$ acts transitively on an $\ell$-manifold $\Sigma^\ell$, where $Q$ has dimension $k=4-\ell$, so that the
associated f\/ibre bundle $P\mathbin{{\times}}_H\Sigma^\ell$ is four-dimensional.

The basic example is the case that the action of $H$ is also free, in which case $\Sigma^\ell$ is a~principal
homogeneous space for $H$ and therefore there is a~commuting free transitive action of a~Lie group $\tilde H$ isomorphic
to $H$.
Choosing a~basepoint on $\Sigma^\ell$ identif\/ies $\tilde H$ and $\Sigma^\ell$ with $H$ and the two actions are the left
and right regular actions.
However, it can be useful to distinguish between the \emph{structure group} $H$, and the \emph{symmetry group} $\tilde
H$: if $P$ is a~principal $H$-bundle over $Q$, then $P\mathbin{{\times}}_H \Sigma^\ell$ is a~principal $\tilde
H$-bundle.

For general $\Sigma^\ell$, there is still an associated bundle $\pi\colon P\mathbin{{\times}}_H\Sigma^\ell\to Q$, but
this does not have any symmetries in general, since there is no longer a~commuting right action of $\tilde H$ on
$\Sigma^\ell$.
However, the f\/ibre of $\lie{h}_Q=P\mathbin{{\times}}_H\lie{h}$ is still a~Lie algebra of vertical vector f\/ields on
$P\mathbin{{\times}}_H\Sigma^\ell$, which will be called ``invariant'', but for consistency with the case
$\Sigma^\ell=H$, the Lie bracket of these vector f\/ields is opposite to the Lie bracket in $\lie{h}_Q$.
The map $\varphi\colon{\mathcal{V}}\to\lie{h}_Q$ therefore induces a~map~$\widehat\varphi$ from~$\pi^*{\mathcal{V}}$ to
the vertical bundle of $P\mathbin{{\times}}_H\Sigma^\ell$.
These bundles both have rank~$\ell$, so let~$M$ be an open subset of $P\mathbin{{\times}}_H\Sigma^\ell$ where
$\widehat\varphi$ is an isomorphism.
Note that $\widehat\varphi$ sends basic sections to ``invariant'' vector f\/ields.

Equip $TM\cong(\pi|_{M})^*({\mathcal{V}}\oplus TQ)$ with the direct sum conformal structure, so that a~conformal
connection $D$ on ${\mathcal{V}}$ induces a~conformal connection on $TM$.
As before, a~modif\/ied conformal connection ${\mathcal{D}}$ can be constructed, using the pullbacks of $\psi$ and
${\mathcal{C}}$: the torsion of this connection will be vertical-valued, and will
satisfy~\eqref{eq:tg1}--\eqref{eq:tg3}.
Since $\varphi$ is injective, the torsion will be selfdual if $(\alpha,\varphi)$ satisfy the gauge f\/ield equation on
$Q$.
The calculation of the curvature of ${\mathcal{D}}$ carries over immediately to this more general setting, hence Joyce's
Lemma can be applied to establish the following result.

\begin{thm}
\label{th:sdgf}
Let $H$ be a~transitive group of diffeomorphisms of an $\ell$-manifold $\Sigma^\ell$.
Suppose that $(\alpha,\varphi)$ is a~solution of the gauge field equation $($i.e., the Nahm, Hitchin or Bogomolny
equation$)$ on a~principal $H$-bundle $P\to Q$, where $Q$ is a~$(4-\ell)$-dimensional background geometry
$($i.e., an Einstein--Weyl, spinor-vortex or Riccati space$)$.
Then on the open subset $M$ of $\pi\colon P\mathbin{{\times}}_H\Sigma^\ell\to Q$ where $\widehat\varphi$ is an
isomorphism, $(\alpha,\varphi)$ identifies $TM$ with $\pi^*({\mathcal{V}}\oplus TQ)$, and the direct sum conformal
structure is selfdual.
\end{thm}

For later work, it will be useful to have a~more explicit description of the construction of this theorem.
Choose a~local conformal frame $e_i$ for $TQ\oplus{\mathcal{V}}$ over $Q$ compatible with the direct sum decomposition,
and a~local section of $P$.
Then, by identifying $M$ locally with $Q\mathbin{{\times}}\Sigma^\ell$ and viewing the connection $\alpha$ as
a~Vect$(\Sigma^\ell)$-valued $1$-form on $Q$, the components of $(\alpha,\varphi)$ with respect to~$e_i$ def\/ine four
vector f\/ields $X_i$ on $\Sigma^\ell$: $X_1,\ldots X_k$ are the components of the connection and $X_{k+1},\ldots X_4$ are
the components of $\varphi$.
Since $e_i$ is a~conformal frame on $TQ\oplus{\mathcal{V}}$, the conformal structure on~$M$ is clearly represented
contravariantly by the metric
\begin{gather}
\label{eq:invmetric}
(e_1-X_1)^2+\dots+(e_k-X_k)^2+{X_{k+1}}^2+\dots+{X_4}^2.
\end{gather}
This is a~metric on $T^{*} M$.
The covariant metric on $TM$ is dual to this, and will be discussed in Section~\ref{s:bgf}.
In fact it is sometimes more convenient to use contravariant metrics, since they push forward easily.

\begin{rems}\quad
\begin{enumerate}\itemsep=0pt
\item[(i)]
Note that the calculations leading to this theorem are entirely formal and so, at least for local considerations, $H$
need not be f\/inite-dimensional, but could be any subgroup of $\Diff(\Sigma^\ell)$.
Hence the conformal aspects of Ward's construction~\cite{Ward:suc} are included in the theorem when the background
geometry is trivial.
\item[(ii)]
Theorem~\ref{th:sdgf} provides a~new interpretation of the switch map~\cite{MaWo:ist, MMW:sbm}.
Let $M$ be a~selfdual conformal $4$-manifold with freely acting $\ell$-dimensional symmetry group $H$.
Then the local quotient $Q=M/H$ is an background geometry of dimension $4-\ell$.
An invariant selfdual Yang--Mills f\/ield (on a~principal bundle $P$) with $\ell$-dimensional gauge group $G$ descends to
a~solution of the gauge f\/ield equation on $Q$, from which a~new selfdual space $\tilde M$ may be constructed.

\centerline{\includegraphics{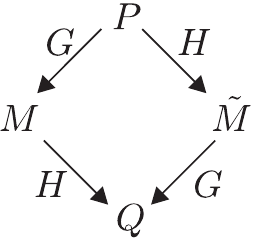}}

Note the following features of this interpretation.
\begin{itemize}\itemsep=0pt
\item
The construction avoids considering the $(4+\ell)$-dimensional manifold $P$ explicitly.
\item
The procedure decomposes into two steps of independent interest: the construction of gauge f\/ields on $Q$ from
$H$-invariant selfdual Yang--Mills f\/ields on $M$, and the construction of selfdual spaces from gauge f\/ields on $Q$.
\item
The group $G$ need not be $\ell$-dimensional, so long as it acts transitively on an $\ell$-manifold.
Hence the constructions of~\cite{DMW:2dg} f\/it into the same framework.
\end{itemize}
\item[(iii)]
In the case that $\Sigma^\ell=H$, $M$ has symmetry group $H$, but enlarging the group $H$ acting on~$\Sigma^\ell$ gives
less symmetry, rather than more symmetry, since there are fewer dif\/feomorphisms of~$\Sigma^\ell$ commuting with the
$H$-action.
In general there will be none, and~$M$ will have no symmetry.
In other words, the group~$H$ is a~\emph{structure group} rather than a~\emph{symmetry group}: only in the case of~$H$
acting freely does it happen that the commuting symmetry group is (isomorphic to)~$H$.
\end{enumerate}
\end{rems}
Selfdual conformal manifolds with symmetry groups which do not act freely are rather special: for instance they are
foliated by the surfaces of equal isotropy group, and the curvature must be invariant under the isotropy representation.
Some examples arise as very special cases of the constructions of this paper.
More precisely, if~$G$ is a~(not necessarily free or transitive) group of dif\/feomorphisms of $\Sigma^\ell$, then~$G$
will be a~symmetry group of~$M$ provided there is a~representation of~$G$ on~${\mathcal{V}}$ such that~$\varphi$ and~$\alpha$ are $G$-invariant.
In the case that~$G$ acts trivially on~${\mathcal{V}}$, this means that the structure group~$H$ reduces to the group of
dif\/feomorphisms of~$\Sigma^\ell$ commuting with~$G$.
On the other hand, there are interesting examples where~$G$ acts nontrivially on~${\mathcal{V}}$, see for
instance~\cite{Hit:hcm}.

The use of structure groups rather than symmetry groups also turns out to be natural when partial reductions are
considered in view of the following remark.
If $H$ acts freely on~$M$ and~$K$ is a~subgroup of~$H$, then $K$ also acts freely on $M$ and the structure group of
$M/K$ is $H/K_0$ acting on $H/K$, where $K_0$ is the maximal normal subgroup of $H$ lying in~$K$.
For example, if $\SU(2)$ acts freely on $M$ then $\SO(3)$ acting on $S^2$ will be the induced structure group for~$M/\Un(1)$.
In Section~\ref{s:bgf}, it will be shown that such partial reductions~$M/K$ arise directly from gauge f\/ields on~$M/H$
with gauge group~$H/K_0$.

\subsection*{Addenda: generalized constructions from null reductions}

The same principles may be applied to null reductions to obtain constructions of selfdual $4$-manifolds from twisted
f\/lat pencils and projective pairs, where the gauge group is a~transitive group of dif\/feomorphisms of a~surface
$\Sigma^2$.
The only change needed, relative to the nondegenerate case, is that the conformal structure on~$M$ is obtained from the
natural pairing between~$VM$ and~$\pi^* TQ$, rather than conformal structures on each summand.

The construction of selfdual conformal structures from f\/lat pencils of connections on a~surface is not at all new: as
discussed in the following Interlude, it underpins Plebanski's heavenly equations and interpretations of hypercomplex
and hyperk\"ahler  structures as topological chiral models.
On the other hand, the construction of selfdual conformal structures from \emph{twisted} f\/lat pencils with gauge group
$\Diff(\Sigma^2)$ has not been studied, as far as I am aware.

The analogous story for projective pairs on a~projective surface underpins the Dunajski--West construction of selfdual
conformal manifolds with a~null Killing vector~\cite{DuWe:acs}.
In~\cite{Cal:sdp}, their construction is shown to be a~reduction of gauge group from $\Diff(\Sigma^2)$ to the subgroup
commuting with a~nonvanishing vector f\/ield.

\section{Interlude: hypercomplex and hyperk\"ahler  structures}
\label{s:hchk}

From the point of view of integrable systems, hypercomplex and hyperk\"ahler  $4$-mani\-folds are considerably simpler
than the general selfdual space.

Recall that a~hypercomplex structure consists of a~triple $I$, $J$, $K$ of integrable complex structures, satisfying the
quaternionic relation $IJ=K$.
It is well known that a~hypercomplex manifold comes equipped with a~unique torsion-free connection~$D$ with
$DI=DJ=DK=0$, called the \emph{Obata connection}~\cite{Oba:cqh}.

A hypercomplex $4$-manifold possesses a~canonical conformal structure, def\/ined by requiring that $(X,IX,JX,KX)$ is
a~conformal frame for any nonzero tangent vector~$X$.
The Obata connection preserves this conformal structure, and is thus a~Weyl connection.
Since~$I$,~$J$,~$K$ are anticommuting orthogonal complex structures, their (weightless) K\"ahler  forms are either all
selfdual or all antiselfdual: I f\/ix the orientation by requiring that they are antiselfdual.
Therefore~$D$ is f\/lat on $L^2 \wedge^2_\asd T^{*} M$, and in particular, $M$ is a~selfdual conformal
manifold.
The hypercomplex structure is hyperk\"ahler  if and only if~$D$ is exact, i.e., the Obata connection preserves a~length
scale, and hence a~metric in the conformal class.

The simplicity of the hypercomplex condition manifests itself in the following local description of hypercomplex
$4$-manifolds.

\begin{thm}
\label{th:MN}
Let $V_0$, $V_1$, $V_2$, $V_3$ be linearly independent vector fields on a~$4$-manifold $M$ and let $\eta_0$, $\eta_1$, $\eta_2$, $\eta_3$
be the dual coframe of $1$-forms $($with $\eta_i(V_j)=\delta_{ij})$.
Define almost complex structures $I$, $J$, $K$ by
\begin{gather*}
I V_0=V_1,\qquad I V_2=V_3,
\qquad
J V_0=V_2,\qquad J V_3=V_1,
\qquad
K V_0=V_3,\qquad K V_1=V_2.
\end{gather*}
Then the following are equivalent.
\begin{enumerate}\itemsep=0pt
\item[$1.$]
The frame $V_i$ satisfies the equations
\begin{gather}
[V_0,V_1]+[V_2,V_3]=0,
\label{eq:MN1}
\\
[V_0,V_2]+[V_3,V_1]=0,
\label{eq:MN2}
\\
[V_0,V_3]+[V_1,V_2]=0.
\label{eq:MN3}
\end{gather}
\item[$2.$]
For each $i$, ${\rm d}\eta_i$ is selfdual with respect to the conformal structure represented by the metric
$g=\eta_0^2+\eta_1^2+\eta_2^2+\eta_3^2$ $($where $I$, $J$, $K$ are antiselfdual$)$.
\item[$3.$]
$(I,J,K)$ is hypercomplex with Obata connection $D$ and $\divg^D\eta_i=0$ for all $i$.
\end{enumerate}
Any hypercomplex $4$-manifold $M$ arises locally in this way, and is hyperk\"ahler  if and only if the vector fields
$V_i$ all preserve a~volume form $\nu$.
\end{thm}

The hyperk\"ahler  version of this theorem is due to Mason--Newman~\cite{MaNe:eym}, following a~construction of
Ashtekar--Jacobson--Smolin~\cite{AJS:hfe} which will be described below.
The hyperk\"ahler  case is simpler, because for $D$ exact, the $\eta_i$ are divergence-free with respect to $D$ if and
only if the $V_i$ preserve a~volume form~-- see~\eqref{eq:divs}.

The general construction was f\/irst written explicitly by Joyce~\cite{Joy:esd}, but the fact that all four-dimensional
hypercomplex structures arise in this way is due to Dunajski~\cite{Dun:tph}.
The description I have given dif\/fers slightly from these references, and owes a~great deal to the approaches of
Hitchin~\cite{Hit:hcm} (see below) and Grant--Strachan~\cite{GrSt:his}.
Since the role of the divergence condition has perhaps not been fully elucidated before, and will be useful later, I
give a~complete proof.

\begin{proof}
Since
\begin{gather*}
\eta_i ( [V_j,V_k] )=-{\rm d}\eta_i(V_j,V_k)
\end{gather*}
for all $i$, $j$, $k$, it is manifest that~(i) and~(ii) are equivalent formulations of the same equations.
Also,~\eqref{eq:MN2} and~\eqref{eq:MN3} clearly imply that~$I$ is integrable since they may be rewritten as
\begin{gather*}
[V_0+{\boldsymbol i} V_1,V_2+{\boldsymbol i} V_3]=[V_0,V_2]-[V_1,V_3] +{\boldsymbol i} ( [V_0,V_3]+[V_1,V_2] )=0.
\end{gather*}
Similarly~\eqref{eq:MN3} and~\eqref{eq:MN1} imply that $J$ is integrable, and~\eqref{eq:MN1} and~\eqref{eq:MN2} imply
that $K$ is integrable.
Now note that for any $1$-form $\eta$ on a~hypercomplex manifold ${\rm d}\eta$ is selfdual if and only if it is
orthogonal to the weightless K\"ahler  forms of $I$, $J$, $K$.
Since $DI=DJ=DK=0$, this is equivalent to $I\eta$, $J\eta$ and $K\eta$ being divergence-free with respect to $D$ (for
instance, $\divg^D(I\eta)=\sum\limits_i \ip{\varepsilon^i,I D_{e_i}\eta}$, and $I$ is skew).
Hence (iii) is equivalent to (i) and (ii).

On any hypercomplex manifold, the conditions $\divg^D\eta=0$ and ${\rm d}\eta^\asd=0$ form a~determined f\/irst-order linear system for a~$1$-form $\eta$, which therefore admits local (non-null) solutions: $(\eta,I\eta,J\eta,K\eta)$
is then a~divergence-free conformal coframe.

It remains to characterize the hyperk\"ahler  case in terms of volume forms~-- or length scales.
Suppose that $\mu=e^{f}\mu_g$ is a~length scale; then $\mu^{-4} V_i=\mu^{-4} \mu_g^{2}\eta_i=e^{-2f}\mu^{-2}\eta_i$ and
so
\begin{gather}
\label{eq:divs}
\divg\big(\mu^{-4}V_i\big)=\divg^D\big(e^{-2f}\mu^{-2}\eta_i\big)=\Ip{D\big(e^{-2f}\mu^{-2}\big),\eta_i} +e^{-2f}\mu^{-2}\divg^D\eta_i.
\end{gather}
Since $\divg^D\eta_i=0$, $D$ preserves the length scale $e^f\mu=e^{2f}\mu_g$ if and only if the $V_i$ all preserve the
volume form $\nu=\mu^{-4}$ (note that $\wedge^4T^{*} M\cong L^{-4}$, using the orientation of $M$).
\end{proof}

There is an equivalent way to describe the divergence-free condition on the coframe, \mbox{using}
spinors~\cite{Dun:tph,MaWo:ist}.
Recall that any conformal $4$-manifold locally admits (weightless) \emph{spin bundles}~$\$^\sdasd$, which are
$\SL(2,{\mathbb{C}})$ bundles such that $\$^\sd\otimes\$^\asd$ is isomorphic to the complexif\/ied weightless cotangent
bundle $L  T^{*} M\otimes{\mathbb{C}}$ with the metric induced by the two area forms.
The conventions are chosen so that $L^2 \wedge^2_\asd T^{*} M\otimes{\mathbb{C}}=S^2\$^\asd$.
On a~hypercomplex manifold $D$ induces a~f\/lat connection on $\$^\asd$ preserving the area form.

For Euclidean reality conditions, the real structure on the $L  T^{*} M\otimes{\mathbb{C}}$ is induced
by (parallel) quaternionic structures on $\$^\sdasd$.
If frames $(\rho_0,\rho_1)$ for $L^{-1} \$^\sd$ and $(\sigma_0,\sigma_1)$ for $\$^\asd$ are chosen so
that the quaternionic structure sends $\rho_0$ to $\rho_1$ and $\sigma_0$ to $\sigma_1$, then
\begin{gather*}
\eta_0+{\boldsymbol i}\eta_1=\rho_0\otimes\sigma_0,
\qquad
\eta_2+{\boldsymbol i}\eta_3=\rho_0\otimes\sigma_1,
\\
\eta_0-{\boldsymbol i}\eta_1=\rho_1\otimes\sigma_1,
\qquad
\eta_2-{\boldsymbol i}\eta_3=-\rho_1\otimes\sigma_0
\end{gather*}
def\/ines a~real conformal coframe $(\eta_0,\eta_1,\eta_2,\eta_3)$.

\begin{prop}
Suppose that $M$ is hypercomplex with Obata connection $D$ and that $(\sigma_0,\sigma_1)$ is a~$D$-parallel frame for
$\$^\asd$.
\begin{enumerate}\itemsep=0pt
\item[$1.$]
Let $\rho$ be a~section of $L^{-1} \$^\sd$.
Then $\rho\otimes\sigma_0$ and $\rho\otimes\sigma_1$ are divergence-free with respect to $D$ if and only if $\rho$
satisfies the Dirac--Weyl equation $\sum\limits_i \varepsilon_i\mathinner{\cdot} D_{e_i}\rho=0$.
\item[$2.$]
Let $\eta$ be a~$1$-form.
Then $\eta\mathinner{\cdot}\sigma_0$ and $\eta\mathinner{\cdot}\sigma_1$ satisfy the Dirac--Weyl equation if and only if
$\divg^D\eta=0$ and ${\rm d}\eta^\asd=0$.
\end{enumerate}
Here the dot denotes the natural $($Clifford$)$ action $T^{*} M\otimes L^{w} \$^\sdasd\to
L^{w-1} \$^\asdsd$.
\end{prop}

\begin{proof}
These are direct calculations:
\begin{enumerate}\itemsep=0pt
\item[(i)]
For $A=0,1$, $\divg^D(\rho\otimes\sigma_A)=\sum\limits_i\varepsilon^i(D_{e_i}\rho\otimes\sigma_A)
=\omega^\asd\bigl(\sum\limits_i\varepsilon^i\mathinner{\cdot} D_{e_i}\rho,\sigma_A\bigr)$.
\item[(ii)]
For $A=0,1$, $\sum\limits_i\varepsilon^i\mathinner{\cdot} D_{e_i}(\eta\mathinner{\cdot}\sigma_A)=
\sum\limits_i\varepsilon^i\mathinner{\cdot}(D_{e_i}\eta)\mathinner{\cdot}\sigma_A
=(\divg^D\eta)\sigma_A+{\rm d}\eta\mathinner{\cdot}\sigma_A$.
\end{enumerate}
In (i) $\omega^\asd$ is the area form on $\$^\asd$, and in (ii) the Clif\/ford action of ${\rm d}\eta^\sd$ is
trivial.
\end{proof}

An important class of solutions to the equations $\divg^D\eta=0$ and ${\rm d}\eta^\asd=0$ is obtained by taking
$\eta={\rm d} r$ where $\Delta^Dr:=\divg^D {\rm d} r=0$.
Evidently the equation $\Delta^Dr=0$ admits local solutions on any hypercomplex $4$-manifold.
Hence the following result is obtained, which is the original construction of Ashtekar--Jacobson--Smolin~\cite{AJS:hfe}
in the hyperk\"ahler  case, and is due to Hitchin~\cite{Hit:hcm} in general.

\begin{thm}
\label{th:s3Nahm}
Let $\Phi$ be a~solution of the Nahm equation $($on a~trivial Riccati space$)$ with gauge group
$\Diff(\Sigma^3)$ for some $3$-manifold $\Sigma^3$.
Then the selfdual space constructed from $\Phi$ is hypercomplex.
Any hypercomplex $4$-manifold arises locally in this way, and is hyperk\"ahler  if and only if there is a~reduction to
the gauge group $\SDiff(\Sigma^3)$.
\end{thm}
\begin{proof}
Let $r$ be an af\/f\/ine coordinate on a~trivial Riccati space $C$, choose a~conformal trivialization of ${\mathcal{E}}$,
and write the components of $\Phi$, which are vector f\/ields on $\Sigma^3$, as $(V_1,V_2,V_3)$.
Then the vector f\/ields $(\partial_r,V_1,V_2,V_3)$ on $M\subset C\mathbin{{\times}}\Sigma^3$
satisfy~\eqref{eq:MN1}--\eqref{eq:MN3} as a~consequence of the Nahm equation, and so the dual frame $({\rm d}
r,\eta_1,\eta_2,\eta_3)$ has ${\rm d}\eta_i$ selfdual (note that I~use the opposite orientation
to~\cite{Hit:hcm}).
Conversely any hypercomplex manifold arises locally in this way by letting~$r$ be a~solution of~$\Delta^Dr=0$ and
setting $(\eta_1,\eta_2,\eta_3)=(I{\rm d} r,J{\rm d} r,K{\rm d} r)$.

Finally note that $\partial_r$ and $V_i$ preserve a~volume form $e^{-4f} {\rm d}
r\wedge\eta_1\wedge\eta_2\wedge\eta_3$ if and only if
\begin{gather*}
{\mathcal{L}}_{\partial_r}\big(e^{-4f}\eta_1\wedge\eta_2\wedge\eta_3\big)=0
\qquad
\text{and}
\qquad
{\rm d} r\wedge {\mathcal{L}}_{V_i}\big(e^{-4f}\eta_1\wedge\eta_2\wedge\eta_3\big)=0.
\end{gather*}
Here I use the fact that ${\rm d} r(V_i)=0$, ${\rm d} r(\partial_r)=1$ and that $\iota_{\partial_r}
(e^{-4f}\eta_1\wedge\eta_2\wedge\eta_3)=0$.
The f\/irst equation says that $e^{-4f}\eta_1\wedge\eta_2\wedge\eta_3$ is an $r$-independent volume element on $\Sigma^3$
(i.e., a~parallel volume form on $C\mathbin{{\times}}\Sigma^3\to C$), and the second equation says that the $V_i$ are
volume-preserving vector f\/ields for each f\/ixed $r$ (i.e., on each f\/ibre of $C\mathbin{{\times}}\Sigma^3\to C$).
This is exactly what it means to have a~reduction to $\SDiff(\Sigma^3)$.
\end{proof}

As this construction involves taking ${\rm d}\eta_0=0$, it is natural to ask if one can f\/ind divergence-free
coframes with ${\rm d}\eta_0=0={\rm d}\eta_1$ and hence formulate the hypercomplex equations as Hitchin
equations on a~trivial spinor-vortex space.

A number of constructions of four-dimensional hyperk\"ahler  metrics from two-dimensional integrable models are known,
due to Park~\cite{Par:sdg}, Ward~\cite{Ward:suc} and (later) Husain~\cite{Hus:sdg}: see~\cite{Uen:ift} for a~review.
Unfortunately, it is not always clear in these constructions whether all hyperk\"ahler  metrics are obtained, what
choices are needed to obtain the integrable model from a~hyperk\"ahler  metric, and if they are compatible with
Euclidean reality conditions.
In particular, as far as I can tell, none of these works establish an equivalent formulation of the Euclidean
hyperk\"ahler  condition as a~Euclidean two-dimensional integrable model.

Indeed, the usual approach is to use the equation $[V_0+{\boldsymbol i} V_1,V_2+{\boldsymbol i} V_3]=0$ to introduce
coordinates $(x,y,u,v)$ such that $V_0+{\boldsymbol i} V_1=\partial_x$ and $V_2+{\boldsymbol i} V_3=\partial_y$.
There are several variations on this theme, since the meaning of $\partial_x$ and $\partial_y$ in $(x,y,u,v)$
coordinates depends on $u$ and $v$, leading to dif\/ferent forms for $V_0-{\boldsymbol i} V_1$ and $V_2-{\boldsymbol i}
V_3$.
This procedure tends to obscure the nature of the choice made to obtain the frame, making it more dif\/f\/icult to argue
that any hyperk\"ahler  metric admits such a~frame.

Fortunately there is an alternative approach, which clarif\/ies the choice of frame, is easily made compatible with any
reality conditions, and generalizes to the hypercomplex case.
The following elementary observation is very well known, at least in the hyperk\"ahler  case.

\begin{prop}
Let $z$ be a~complex function on a~hypercomplex $4$-manifold $(M,D)$ which is holomorphic with respect to one of the
complex structures.
Then $\Delta^D z=0$, so that ${\rm d} z$ is a~complex null $1$-form which is divergence-free with respect to $D$.
\end{prop}
\begin{proof}
If $z$ is $I$-holomorphic, then $I{\rm d} z:=-{\rm d} z\circ I = -{\boldsymbol i}{\rm d} z$.
Now $I$ is skew and $DI=0$, so $\divg^D(I{\rm d} z)=\sum\limits_i\ip{\varepsilon^i,I D_{e_i}{\rm d} z}=0$
since ${\rm d}^2 z=0$.
Hence $\Delta^D z=\divg^D({\rm d} z-{\boldsymbol i} I{\rm d} z)=0$.
\end{proof}
\begin{rem}
I have presented this observation using language adapted to Euclidean signature manifolds.
In Kleinian signature, some of the complex structures are imaginary, so that ${\boldsymbol i} I$ (say) is a~real
involution, inducing a~decomposition $TM=T^{+}M\oplus T^{-}M$ into its $\pm1$ eigenspaces.
These distributions are integrable, and the analogue of a~holomorphic function is a~function constant on one of the
families of integral surfaces~-- such functions can of course be real valued.
Note that the orientation $2$-forms of these integral surfaces are the (weightless) K\"ahler  forms of the null complex
structures $J\pm{\boldsymbol i} K$, which are decomposable (and up to rotation, $J$ and ${\boldsymbol i} K$ are real).

The Euclidean and Kleinian cases can be considered together by f\/irst working on a~complexif\/ied hypercomplex manifold,
then imposing reality conditions.
The above proposition applies equally in the complexif\/ied setting.
Other, more Kleinian, arguments are also available (cf.~\cite{CMN:cs4,Ple:sce}).
\end{rem}

Before discussing the non-null reduction to two dimensions, I will brief\/ly discuss the null reduction (closely related
to Plebanski's heavenly equations~\cite{Ple:sce}) which is used in the literature to relate the selfdual vacuum equation
to a~topological chiral model~\cite{Par:sdg,Ward:suc}.
I present the generalizations to the hypercomplex case, following~\cite{Dun:tph,GrSt:his}.

\subsection*{Topological models and heavenly equations}

Choose independent functions $(w,z)$, both holomorphic with respect to $I$; then ${\rm d} w,{\rm d} z$ are
null $1$-forms with $\ip{{\rm d} w,{\rm d} z}=0$, and one can take these to be $\eta_0+{\boldsymbol
i}\eta_1$ and $\eta_2+{\boldsymbol i}\eta_3$.
This can be done compatibly with Kleinian reality conditions by taking ${\boldsymbol i} I$ and $(w,z)$ real, but is
incompatible with Euclidean reality conditions.

Locally, $M$ is a~bundle of null surfaces over a~quotient surface $N$ with coordinates $w,z$.
Choosing f\/ibre coordinates amounts to choosing a~local trivialization of this bundle, and locally one can take
$M=N\mathbin{{\times}}\Sigma^2$.
Then $V_0-{\boldsymbol i} V_1=\partial_w-\alpha$, $V_2-{\boldsymbol i} V_3=\partial_z-\beta$, $V_0+{\boldsymbol i}
V_1=\phi$, $V_2+{\boldsymbol i} V_3=\psi$, where $\alpha$, $\beta$, $\phi$, $\psi$ are vector f\/ields tangent to the f\/ibres.
The equations~\eqref{eq:MN1}--\eqref{eq:MN3} now read:
\begin{gather*}
[\partial_w-\alpha,\partial_z-\beta]=0,
\qquad
[\phi,\psi]=0,
\qquad
[\partial_z-\alpha,\phi]+[\partial_w-\beta,\psi]=0,
\end{gather*}
which are the equations for a~pencil of f\/lat connections ${\rm d}+A+\lambda\Phi$ with gauge group
$\Diff(\Sigma^2)$, where $A=-\alpha {\rm d} w-\beta {\rm d} z$ and $\Phi=-\psi{\rm d}
z+\phi{\rm d} w$.
(This is also known as a~topological chiral or sigma model.)

Plebanski's f\/irst and second heavenly equations~\cite{Ple:sce}, and their generalizations to the hypercomplex case (due
to Grant--Strachan~\cite{GrSt:his} and Dunajski~\cite{Dun:tph} respectively) are obtained by f\/ixing the gauge freedom in
dif\/ferent ways.

$\bullet$ First, since ${\rm d}+A$ is f\/lat, one can set $A=0$ (i.e., $\alpha=\beta=0$),
then integrate the equation ${\rm d}\Phi=0$ (i.e., $\phi_z+\psi_w=0$) to get $V_0-{\boldsymbol i}
V_1=\partial_w$, $V_2-{\boldsymbol i} V_3=\partial_z$, $V_0+{\boldsymbol i} V_1=U_w$, $V_2+{\boldsymbol i} V_3=-U_z$ for
a~vector f\/ield $U$ tangent to the f\/ibres (with $U_z$ and $U_w$ linearly independent).
The remaining equation is $[U_w,U_z]=0$.
If $U$ is area preserving on the f\/ibres, with local hamiltonian $\Omega$, then Plebanski's f\/irst equation
$\{\Omega_w,\Omega_z\}=1$ is obtained, where $\{\cdot,\cdot\}$ denotes the Poisson bracket with respect to a~suitably
scaled area form on the f\/ibres.

$\bullet$ Second, since $[\phi,\psi]=0$, one can choose the f\/ibre coordinates $(x,y)$ so that
$\phi=\partial_x$ and $\psi=-\partial_y$, then integrate the equation $\alpha_x-\beta_y=0$ to give $V_0-{\boldsymbol i}
V_1=\partial_w-\gamma_y$, $V_2-{\boldsymbol i} V_3=\partial_z-\gamma_x$, $V_0+{\boldsymbol i} V_1=\partial_x$,
$V_2+{\boldsymbol i} V_3=-\partial_y$.
The remaining equation is $\gamma_{xw}-\gamma_{yz}=[\gamma_x,\gamma_y]$.
Again the area preserving condition reduces everything to a~single function $\Theta$, satisfying Plebanski's second
heavenly equation $\Theta_{xw}+\Theta_{yz}=\{\Theta_x,\Theta_y\}$.

All complexif\/ied hypercomplex and hyperk\"ahler  metrics are obtained from these constructions, but information about
Euclidean real slices is lost.

\subsection*{Hypercomplex structures from the Hitchin equations}

In order to obtain a~formulation compatible with Euclidean reality conditions, take $z$ to be $I$-holomorphic and
$\tilde z$ to be $(-I)$-holomorphic.
Then ${\rm d} z$ and ${\rm d}\tilde z$ are null and so ${\rm d} z+{\rm d} \tilde z$ and
$I({\rm d} z+{\rm d}\tilde z)=-{\boldsymbol i}({\rm d} z-{\rm d}\tilde z)$ are orthogonal,
closed, divergence-free $1$-forms of the same length.
Generically, this length will be nonzero on a~dense open set, and one can take $\eta_0+{\boldsymbol
i}\eta_1={\rm d} z$ and $\eta_0-{\boldsymbol i}\eta_1={\rm d}\tilde z$.
Euclidean reality conditions are easily obtained by setting $\tilde z=\bar z$ for $I$ real.

\begin{thm}
\label{th:s2Hit}
Let $(A,\Phi)$ be a~solution to the Hitchin equations $($on a~trivial spinor-vortex space$)$, with gauge group
$\Diff(\Sigma^2)$ for some $2$-manifold $\Sigma^2$.
Then the selfdual space constructed from $(A,\Phi)$ is hypercomplex.
Any hypercomplex $4$-manifold arises locally in this way, and is hyperk\"ahler  if and only if there is a~reduction to
the gauge group $\SDiff(\Sigma^2)$.
\end{thm}

\begin{proof}
Choose conformal coordinates $z=x+{\boldsymbol i} y$, $\tilde z=x-{\boldsymbol i} y$ on a~trivial spinor-vortex space
and write $A=\alpha {\rm d} z+\widetilde\alpha {\rm d}\tilde z$, $\Phi=\phi {\rm d}
z+\widetilde\phi {\rm d}\tilde z$, where $\alpha,\widetilde\alpha, \phi,\widetilde\phi$ are complex vector
f\/ields on $\Sigma^2$.
Then the Hitchin equations become
\begin{gather*}
[\partial_z-\alpha,\widetilde\phi]=0,
\qquad
[\partial_{\tilde z}-\widetilde\alpha,\phi]=0,
\qquad
[\partial_z-\alpha,\partial_{\tilde z}-\widetilde\alpha]-[\phi,\widetilde\phi]=0,
\end{gather*}
and so the vector f\/ields $V_0-{\boldsymbol i} V_1=\partial_z-\alpha$, $V_0+{\boldsymbol i} V_1=\partial_{\tilde
z}-\widetilde\alpha$, $V_2-{\boldsymbol i} V_3=\widetilde\phi$ and $V_2+{\boldsymbol i} V_3=\phi$ satisfy
equations~\eqref{eq:MN1}--\eqref{eq:MN3}.
Hence the selfdual space is hypercomplex, with a~divergence-free coframe $(\eta_0,\eta_1,\eta_2,\eta_3)$ such that
$\eta_0+{\boldsymbol i} \eta_1={\rm d} z$ and $\eta_0-{\boldsymbol i}\eta_1={\rm d}\tilde z$ so that
$\eta_0$ and $\eta_1$ are closed.
Conversely, any hypercomplex $4$-manifold admits such a~divergence-free coframe, so the distribution generated by $V_2$
and $V_3$ (i.e., annihilated by $\eta_0$ and $\eta_1$) is integrable.
Under this assumption the form of the vector f\/ields given above is entirely general, and so any hypercomplex structure
arises in this way.

The characterization of the hyperk\"ahler  case is entirely analogous to Theorem~\ref{th:s3Nahm}: $V_i$ preserve
a~volume form $e^{-4f}{\rm d} z\wedge {\rm d}\tilde z\wedge\eta_2\wedge\eta_3$ if and only if the area
form $e^{-4f}\eta_2\wedge\eta_3$ is parallel with respect to the connection $A=\alpha {\rm d}
z+\widetilde\alpha {\rm d} \tilde z$, and the vector f\/ields $\phi$, $\widetilde\phi$ are area-preserving (on each
f\/ibre), which is exactly what it means to have a~reduction to $\SDiff(\Sigma^2)$.
\end{proof}

\begin{rem}
This result can also be interpreted in Kleinian signature, when the Hitchin equations are replaced by harmonic maps into
a~Lie group or the principal chiral model.
The latter is the context for Husain's formulation~\cite{Hus:sdg}.
The non-null reduction (in the hyperk\"ahler  case) is also discussed brief\/ly by Ward~\cite{Ward:suc} and
Ueno~\cite{Uen:ift}.
\end{rem}

Although this discussion has been local, there are intriguing connections with the global geometry of elliptically
f\/ibred K3 surfaces.
Yau's solution of the Calabi problem shows that on any K3 surface there is a~unique hyperk\"ahler  metric in each
K\"ahler  class, but no explicit description is known.
Any such hyperk\"ahler  metric will correspond to a~solution of the $\SDiff(\Sigma^2)$ Hitchin equations, once
a~holomorphic function is chosen on a~suitable open subset of the K3 surface to def\/ine the dimensional reduction.
Now there are K3 surfaces which admit f\/ibrations over $\CP1$ (meromorphic functions) with elliptic curves as f\/ibres,
and, generically, 24 singular f\/ibres.
On the complement of the singular f\/ibres, there is therefore a~dimensional reduction to the Hitchin equations (on $\CP1$
minus 24 points) with gauge group $\SDiff(T^2)$.
The work of Gross and Wilson~\cite{GrWi:lcsK3} shows that this solution is well approximated by an Abelian solution
(gauge group $T^2$) def\/ining a~`semi-f\/lat' metric.

\begin{rem}
Continuing the development of this section, it is natural to ask if the hypercomplex equations are equivalent to the
$\Diff(S^1)$ Bogomolny equation on ${\mathbb{R}}^3$.
In fact, it is shown in~\cite{Cal:sde} that the selfdual space constructed from a~solution of the $\Diff(S^1)$
generalized Bogomolny equation on an Einstein--Weyl space $B$ is hypercomplex if $B$ is ``hyperCR'' (see
Section~\ref{s:s1Hit}).
However, not all hypercomplex structures arise this way, since for the hyperk\"ahler  case in particular,
$\SDiff(S^1)=\Un(1)$ and only metrics with symmetry are obtained.
\end{rem}

\begin{rem}\label{r:zdg}\sloppy
The Mason--Newman--Dunajski--Joyce construction of Theorem~\ref{th:MN} also has a~gau\-ge-theo\-retic interpretation, of
course: the $V_i$ satisfy gauge f\/ield equations on a~trivial zero-dimen\-sio\-nal geometry~(!) with gauge group
$\Diff(\Sigma^4)$ for some $4$-manifold $\Sigma^4$.
It has been observed in many places (in particular~\cite{MaNe:eym}) that these equations are the reduction of the
selfdual Yang--Mills equations on ${\mathbb{R}}^4$ by four translations, and it would therefore be natural to extend the
integrable background geometries programme to the zero-dimensional case.
The background geometry is a~four-dimensional conformal vector space ${\mathcal{M}}$ together with an element
${\mathcal{Y}}$ of $\wedge^2_\asd{\mathcal{M}}^{*}\otimes{\mathcal{M}}$ which acts as a~right hand side for
equations~\eqref{eq:MN1}--\eqref{eq:MN3}.
The background equation is an unpleasant quadratic condition on ${\mathcal{Y}}$, which I have left as an exercise for
the enthusiastic reader: including all the constructions involving this zero-dimensional geometry would have added
unnecessarily to the length of this paper.
\end{rem}

\section{Background geometries from gauge f\/ields}
\label{s:bgf}

Theorem~\ref{th:sdgf} admits the following generalization.

\begin{thm}
\label{th:bgf}
Let $H$ be a~transitive group of diffeomorphisms of an $\ell$-manifold $\Sigma^\ell$.
Suppose that $(\alpha,\varphi)$ is a~solution of the gauge field equation on a~principal $H$-bundle $P\to Q$, where $Q$
is a~$k$-dimensional background geometry and $k+\ell\leqslant 4$.
Then the open subset of $\pi\colon P\mathbin{{\times}}_H\Sigma^\ell\to Q$ where $\widehat\varphi$ is surjective carries
naturally the structure of a~$(k+\ell)$-dimensional background geometry.
\end{thm}

\begin{proof}
The idea is to apply Theorem~\ref{th:sdgf} using the group $H\mathbin{{\times}}{\mathbb{R}}^{4-k-\ell}$ acting on
$\Sigma^\ell\mathbin{{\times}}{\mathbb{R}}^{4-k-\ell}$.
Suppose that $\widehat\varphi$ is surjective at some point $x$ of the f\/ibre of $P\mathbin{{\times}}_H\Sigma^\ell$ over
$q\in Q$.
Let $K\leqslant{\mathcal{V}}$ be the kernel of $\widehat\varphi$ at $x$.
Then there exists a~solution $(\alpha_0,\varphi_0)$ of the gauge f\/ield equation, with gauge group
${\mathbb{R}}^{4-k-\ell}$, def\/ined on a~neighbourhood of $q$, such that $\kernel\varphi_0\cap K=\{0\}$ at $q$: the
linear gauge f\/ield equation can be solved with any initial condition.
The pair $(\widehat\varphi,\widehat\varphi_0)$ is therefore an isomorphism on a~neighbourhood
$M=U\mathbin{{\times}}{\mathbb{R}}^{4-k-\ell}$ of $\{x\}\mathbin{{\times}}{\mathbb{R}}^{4-k-\ell}$ in
$\bigl(P\mathbin{{\times}}_H\Sigma^\ell\bigr)\mathbin{{\times}}{\mathbb{R}}^{4-k-\ell}$.

The resulting selfdual conformal structure on $M$ clearly admits ${\mathbb{R}}^{4-k-\ell}$ as a~symmetry group.
Hence the quotient $U\subseteq P\mathbin{{\times}}_H\Sigma^\ell$ is a~$(k+\ell)$-dimensional background geometry by the
results of Section~\ref{s:sdbg}.
However, the conformal metric on this quotient is clearly independent of the choice of $(\alpha_0,\varphi_0)$, since the
pushforward of the inverse metric~\eqref{eq:invmetric} from $\Sigma^\ell\mathbin{{\times}}{\mathbb{R}}^{4-k-\ell}$ to
$\Sigma^\ell$ kills the components of the vector f\/ields $X_i$ in ${\mathbb{R}}^{4-k-\ell}$.
Hence the conformal metric is uniquely def\/ined wherever $\widehat\varphi$ is surjective.
One also sees that the other f\/ields def\/ining the background geometry are well def\/ined, but the details here depend on
the geometry and are rather complicated.
The general formulae are given in the following subsections: these will complete the proof, since they are manifestly
well def\/ined.
\end{proof}

I now give some explicit formulae, using the notation $(A,\Phi)$ for the gauge f\/ields, rather than $(\alpha,\varphi)$.
From the above description it is clear at least that the conformal metric on the $(k+\ell)$-dimensional background
geometry may be represented contravariantly by
\begin{gather*}
(e_1-X_1)^2+\cdots+(e_k-X_k)^2+{X_{k+1}}^2+\cdots+{X_4}^2,
\end{gather*}
but note that $X_{k+1},\ldots X_4$ are $4-k$ vector f\/ields on an $\ell$-manifold $\Sigma^\ell$ with $\ell\leqslant 4-k$,
so inverting this metric is only straightforward when $k+\ell=4$.
It will be convenient therefore to introduce a~volume form $\nu$ on $\Sigma^\ell$ and hence present the explicit
formulae in an `SDif\/f-gauge'.
This will also make it easy to understand the volume-preserving case.

\subsection{Riccati space constructions}

First suppose that $\Phi$ is a~generalized Nahm f\/ield with values in $\Vect(\Sigma^3)$ and that $\nu$ is a~volume form
on $\Sigma^3$.
Then the contravariant metric
\begin{gather*}
\partial_r^2+\ip{\Phi,\Phi}
\qquad
\text{is dual to}
\quad
{\rm d} r^2+\ip{\eta,\eta},
\end{gather*}
{where}
\begin{gather*}
\eta=\frac{\nu(\Phi\mathbin{{\times}}\Phi,\cdot)}{\nu(\Phi\mathbin{{\times}}\Phi\mathbin{{\times}}\Phi)}
\end{gather*}
is a~section of ${\mathcal{E}}^{*}$ with values in $\Omega^1(\Sigma^3)$ (the space of $1$-forms).
Here and in the following~$\mathbin{{\times}}$ denotes the cross product ${\mathcal{E}}^{*}\otimes{\mathcal{E}}^{*}\to
L^{-1} {\mathcal{E}}^{*}$ (given by the wedge product and Hodge star operator): the $\Vect(\Sigma^3)$
values of $\Phi$ are then contracted into the entries of the volume form~$\nu$.
${\rm d} r^2+\ip{\eta,\eta}$ is the covariant form of the selfdual conformal structure obtained from~$\Phi$.

Now suppose that $\Phi$ is a~generalized Nahm f\/ield with values in $\Vect(\Sigma^2)$ and that $\nu$ is an area form on
$\Sigma^2$.
Let $F$ be an Abelian Nahm f\/ield acting on ${\mathbb{R}}$ with coordinate $\theta$.
Then $\Phi+F\partial_\theta$ is a~generalized Nahm f\/ield with values in $\Vect(\Sigma^2\mathbin{{\times}}{\mathbb{R}})$
and $\nu\wedge {\rm d}\theta$ is a~volume form on $\Sigma^2\mathbin{{\times}}{\mathbb{R}}$.
The selfdual conformal structure is therefore represented by ${\rm d} r^2+\ip{\eta,\eta}$ where now
\begin{gather*}
\eta=\frac{\nu(\Phi\mathbin{{\times}}\Phi){\rm d}\theta+F\mathbin{{\times}}\nu(\Phi,\cdot)}
{\ip{F,\nu(\Phi\mathbin{{\times}}\Phi)}}.
\end{gather*}
Straightforward manipulations and triple cross product identities may be used to rediagonalize the conformal metric
\begin{gather*}
{\rm d} r^2+\ip{\eta,\eta}= {\rm d} r^2+\frac{|\nu(\Phi\mathbin{{\times}}\Phi)|^2{\rm d}\theta^2
+2\ip{\nu(\Phi\mathbin{{\times}}\Phi),F\mathbin{{\times}}\nu(\Phi,\cdot)}{\rm d}\theta
+|F\mathbin{{\times}}\nu(\Phi,\cdot)|^2}{\ip{F,\nu(\Phi\mathbin{{\times}}\Phi)}^2}
\\
\phantom{{\rm d} r^2+\ip{\eta,\eta}}
={\rm d} r^2+\frac{|\nu(\Phi\mathbin{{\times}}\Phi)|^2
|F\mathbin{{\times}}\nu(\Phi,\cdot)|^2-\ip{\nu(\Phi\mathbin{{\times}}\Phi),F\mathbin{{\times}}\nu(\Phi,\cdot)}^2}
{|\nu(\Phi\mathbin{{\times}}\Phi)|^2\ip{F,\nu(\Phi\mathbin{{\times}}\Phi)}^2}
\\
\phantom{{\rm d} r^2+\ip{\eta,\eta}=}
{}+\frac{|\nu(\Phi\mathbin{{\times}}\Phi)|^2}{\ip{F,\nu(\Phi\mathbin{{\times}}\Phi)}^2}
\biggl({\rm d}\theta+\frac{\ip{\nu(\Phi\mathbin{{\times}}\Phi),F\mathbin{{\times}}\nu(\Phi,\cdot)}}
{|\nu(\Phi\mathbin{{\times}}\Phi)|^2}\biggr)^2
\\
\phantom{{\rm d} r^2+\ip{\eta,\eta}}
={\rm d} r^2+\frac{\ip{\nu(\Phi,\cdot),\nu(\Phi,\cdot)}}{|\nu(\Phi\mathbin{{\times}}\Phi)|^2}
+\frac{|\nu(\Phi\mathbin{{\times}}\Phi)|^2}{\ip{F,\nu(\Phi\mathbin{{\times}}\Phi)}^2}
\biggl({\rm d}\theta+\frac{\ip{\nu(\Phi\mathbin{{\times}}\Phi),F\mathbin{{\times}}\nu(\Phi,\cdot)}}
{|\nu(\Phi\mathbin{{\times}}\Phi)|^2}\biggr)^2.
\end{gather*}
Hence the conformal structure on the quotient by $\partial_\theta$ is represented by the metric
\begin{gather*}
{\rm d} r^2+\frac{\ip{\nu(\Phi,\cdot),\nu(\Phi,\cdot)}}{|\nu(\Phi\mathbin{{\times}}\Phi)|^2},
\end{gather*}
which is, of course, inverse to $\partial_r^2+\ip{\Phi,\Phi}$: note in particular that
$\ip{\nu(\Phi\mathbin{{\times}}\Phi),\Phi}=0$, expressing the fact that neither the components of $\Phi$, nor the dual
$1$-form components of $\nu(\Phi,\cdot)$ are linearly independent (pointwise on $\Sigma^2$).

I want to give the Weyl structure in an SDif\/f-gauge, with representative metric
\begin{gather*}
|\nu(\Phi\mathbin{{\times}}\Phi)|^2{\rm d} r^2+\ip{\nu(\Phi,\cdot),\nu(\Phi,\cdot)}.
\end{gather*}
This metric can be conveniently diagonalized as
\begin{gather*}
g=|\nu(\Phi\mathbin{{\times}}\Phi){\rm d} r+\nu(\Phi,\cdot)|^2
\end{gather*}
using the fact that $\ip{\nu(\Phi\mathbin{{\times}}\Phi),\nu(\Phi,\cdot)}=0$.
It takes quite a~bit of calculation to compute the Jones--Tod Weyl structure $\omega$ in this gauge (i.e.,
$D=D^g+\omega$), but the result is
\begin{gather*}
\omega=\frac{\Ip{2{\mathcal{B}}\bigl(\nu(\Phi\mathbin{{\times}}\Phi)\bigr)
-\nu(\Phi\mathbin{{\times}}\Phi)\mathbin{{\times}}\divg_\nu\Phi, \nu(\Phi\mathbin{{\times}}\Phi){\rm d}
r+\nu(\Phi,\cdot)}}{|\nu(\Phi\mathbin{{\times}}\Phi)|^2}.
\end{gather*}
Writing $\nu={\rm d} p\wedge {\rm d} q$ and expanding the cross products in components gives a~fuller
expression
\begin{gather}
\label{eq:RicEW}
g= \eta_1^2+\eta_2^2+\eta_3^2,
\\
\nonumber
\omega= \frac{2\sum\limits_{j,k} {\mathcal{B}}_{jk} \nu_j\eta_k
-\sum\limits_{i,j,k}\varepsilon_{ijk}\nu_i(\phi^j_p+\psi^j_q)\eta_k} {\nu_1^2+\nu_2^2+\nu_3^2},
\end{gather}
where
\begin{gather*}
\eta_i= \nu_i {\rm d} r + \phi^i {\rm d} q - \psi^i {\rm d} p,
\\
\nu_1= \phi^2\psi^3-\phi^3\psi^2,
\qquad
\nu_2= \phi^3\psi^1-\phi^1\psi^3,
\qquad
\nu_3= \phi^1\psi^2-\phi^2\psi^1,
\end{gather*}
 and
\begin{gather*}
\Phi=\big(\phi^1,\phi^2,\phi^3\big)\partial_p+\big(\psi^1,\psi^2,\psi^3\big)\partial_q.
\end{gather*}
Well, nobody said it was going to be easy!
The reward is the knowledge that this Weyl structure is Einstein--Weyl if
${\mathcal{B}}$ satisf\/ies the matrix Riccati equation and $\Phi$ is a~generalized Nahm f\/ield on this Riccati space.

The construction of spinor-vortex spaces from Riccati spaces is perhaps the most awkward to make explicit, because of
the gauge freedom in the bundle ${\mathcal{W}}$ on a~spinor-vortex space.
If $\Phi$ is a~generalized Nahm f\/ield on a~Riccati space with values in $\Vect(\Sigma^1)$ for a~$1$-manifold $\Sigma^1$
with coordinate $t$, then the only natural way to proceed is to take ${\mathcal{W}}$ to be the kernel of $\Phi$ in the
pullback of ${\mathcal{E}}$ to $C\mathbin{{\times}}\Sigma^1$.
This kernel is not preserved, in general, by the connection $D$ on ${\mathcal{E}}$, but is preserved by the conformal
connection
\begin{gather*}
\nabla=D+\frac{\Phi_r\mathinner{\vartriangle}\Phi{\rm d} r+\dot\Phi\mathinner{\vartriangle}\Phi {\rm d}
t}{|\Phi|^2},
\end{gather*}
where an af\/f\/ine coordinate $r$ and a~$D$-parallel trivialization of ${\mathcal{E}}$ have been introduced.
The complex structure on ${\mathcal{W}}$ is given by cross product with $\Phi/|\Phi|$, while the holomorphic structure
is def\/ined using the connection
\begin{gather*}
\nabla+\frac{\ip{{\mathcal{B}}\Phi,\Phi}\iden}{|\Phi|^2}{\rm d} r
\end{gather*}
on ${\mathcal{W}}$.
The other two f\/ields on the spinor-vortex space are
\begin{gather*}
{\mathcal{C}}={\mathcal{B}}-\frac{{\mathcal{B}}\Phi\otimes\Phi+\Phi\otimes{\mathcal{B}}\Phi}{|\Phi|^2}
+\frac{\ip{{\mathcal{B}}\Phi,\Phi}}{2|\Phi|^2}\left(\iden+\frac{\Phi\otimes\Phi} {|\Phi|^2}\right),
\\
\psi=\frac{\Phi\mathbin{{\times}}(2{\mathcal{B}}\Phi+\Phi\mathbin{{\times}}\dot\Phi)}{|\Phi|^2}.
\end{gather*}
These f\/ields satisfy the spinor-vortex equations if $\Phi$ is a~generalized Nahm f\/ield on a~Riccati space.
(Since $\Phi/|\Phi|$ is $D$-parallel, it is reasonably straightforward to check~\eqref{eq:sv1}--\eqref{eq:sv2} directly,
although~\eqref{eq:sv3} is harder.)

\subsection{Spinor-vortex space constructions}

The approach here is the same as in the previous subsection, and the details are slightly less complicated.
For explicitness, introduce conformal coordinates $z=x+{\boldsymbol i} y$, $\tilde z=x-{\boldsymbol i} y$ on~$N$ and let
$\nu$ be an area form on~$\Sigma^2$.
First suppose that $(\Phi,\widetilde\Phi, \alpha{\rm d} z+\widetilde\alpha{\rm d}\tilde z)$ is
a~generalized Hitchin f\/ield with values in~$\Vect(\Sigma^2)$.
Then the contravariant metric
\begin{gather*}
4(\partial_z-\alpha)(\partial_{\tilde z}-\widetilde\alpha)+4\ip{\Phi,\widetilde\Phi}
\qquad
\text{is dual to}
\quad
{\rm d} z {\rm d}\tilde z+\ip{\eta,\widetilde\eta},
\end{gather*}
{where}
\begin{gather*}
\eta=\frac{\nu(\Phi,\cdot)+\nu(\Phi,\alpha){\rm d} z+\nu(\Phi,\widetilde\alpha) {\rm d}\tilde
z} {\nu(\Phi,\widetilde\Phi)}
\qquad
\text{and}
\qquad
\widetilde\eta=\frac{\nu(\widetilde\Phi,\cdot)+\nu(\widetilde\Phi,\alpha){\rm d} z
+\nu(\widetilde\Phi,\widetilde\alpha){\rm d}\tilde z} {\nu(\widetilde\Phi,\Phi)}.
\end{gather*}

Now suppose that $(\Phi\otimes\partial_t,\widetilde\Phi\otimes\partial_t, (\alpha{\rm d} z+\widetilde\alpha
{\rm d}\tilde z)\otimes \partial_t)$ is a~generalized Hitchin f\/ield with values in $\Vect(\Sigma^1)$ and that
$(F\partial_\theta,\widetilde F\partial_\theta,(\beta{\rm d} z+\widetilde\beta{\rm d}\tilde
z)\otimes\partial_\theta)$ is an Abelian Hitchin f\/ield acting on~${\mathbb{R}}$ with coordinate $\theta$.
Adding these together produces a~generalized Hitchin f\/ield with values in
$\Vect(\Sigma^1\mathbin{{\times}}{\mathbb{R}})$ and ${\rm d} t\wedge{\rm d}\theta$ is an area form on
$\Sigma^1\mathbin{{\times}}{\mathbb{R}}$.
Direct substitution gives
\begin{gather*}
\eta=\frac{F({\rm d} t+\alpha{\rm d} z+\widetilde\alpha{\rm d}\tilde z)
-\Phi({\rm d}\theta+\beta{\rm d} z+\widetilde\beta{\rm d}\tilde z)}
{F\widetilde\Phi-\Phi\widetilde F},
\\
\widetilde\eta= \frac{\widetilde F({\rm d} t+\alpha{\rm d} z+\widetilde\alpha{\rm d}\tilde
z) -\widetilde\Phi({\rm d}\theta+\beta{\rm d} z+\widetilde\beta{\rm d}\tilde z)} {\widetilde
F\Phi-\widetilde\Phi F}.
\end{gather*}
It is straightforward to rediagonalize $\ip{\eta,\widetilde\eta}$ to obtain the metric
\begin{gather*}
{\rm d} z{\rm d}\tilde z+\frac{({\rm d} t + \alpha{\rm d} z + \widetilde\alpha
{\rm d}\tilde z)^2} {4\Phi\widetilde\Phi}
\\
\qquad
{}-\frac{\Phi\widetilde\Phi}{(F\widetilde\Phi-\Phi\widetilde F)^2} \biggl({\rm d}\theta+\beta{\rm d}
z+\widetilde\beta{\rm d}\tilde z-\frac{F \widetilde\Phi+\Phi\widetilde F}{2 \Phi\widetilde\Phi}({\rm d}
t+\alpha{\rm d} z+\widetilde\alpha{\rm d}\tilde z)\biggr)^2.
\end{gather*}
As in the Riccati space construction, the conformal structure is easy to obtain, and is dual to
$4(\partial_z-\alpha)(\partial_{\tilde z}-\widetilde\alpha) +4\Phi^2\partial_t^2$, while more work is required to
compute the Jones--Tod Weyl structure~$\omega$.
The result, again in an SDif\/f-gauge (i.e., the gauge given by ${\rm d} t$) is reasonably simple, however,
\begin{gather}
g= 4\Phi\widetilde\Phi {\rm d} z{\rm d}\tilde z +({\rm d} t + \alpha{\rm d} z +
\widetilde\alpha {\rm d}\tilde z)^2,
\nonumber
\\
\omega= \biggl( \dot\alpha-\frac{2\widetilde{\mathcal{C}}\Phi}{\widetilde\Phi}\biggr){\rm d} z +\biggl(
\dot{\widetilde\alpha}-\frac{2{\mathcal{C}}\widetilde\Phi}{\Phi}\biggr){\rm d}\tilde z
-\frac12\biggl(\frac{\psi+\dot\Phi}{\Phi}+ \frac{\widetilde\psi+\dot{\widetilde\Phi}}{\widetilde\Phi}\biggr)
({\rm d} t + \alpha{\rm d} z + \widetilde\alpha {\rm d}\tilde z).
\label{eq:SVtoEW}
\end{gather}
This Weyl structure is Einstein--Weyl if $(\Phi,\widetilde\Phi,\alpha{\rm d} z + \widetilde\alpha
{\rm d}\tilde z)$ is a~generalized Hitchin f\/ield on a~spinor-vortex space.

\subsection{Einstein--Weyl constructions}

For completeness, I record here the explicit form of the generalized Jones--Tod construction of selfdual spaces from
Einstein--Weyl spaces~\cite{JoTo:mew,Cal:sde}.
The conformal structure on $M$ is obtained from the $\Diff(\Sigma^1)$ monopole $(A,\Phi)\partial_t$ on $B$ by the
formula
\begin{gather*}
\mathsf{c}=\mathsf{c}_B + \Phi^{-2}({\rm d} t+A)^2,
\end{gather*}
where $t$ is a~coordinate on $\Sigma^1$.
Compatible metrics for $\mathsf{c}$ are easily obtained by introducing a~compatible metric $g_B=\mu^{-2}\mathsf{c}_B$ on
$B$ and writing $\Phi=V\mu^{-1}$.
Then $g_B+V^{-2}({\rm d} t+A)^2$, $Vg_B+V^{-1}({\rm d} t+A)^2$ and $V^2g_B+({\rm d} t+A)^2$ are all
possibilities.
The latter may be written more invariantly as $\Phi^2\mathsf{c}_B +({\rm d} t+A)^2$: it is the SDif\/f-gauge
determined by ${\rm d} t$.

\subsection*{Addendum: null and non-null reductions}

Unlike the preceding addenda, the main observation here is a~negative one: the methods of this section do not extend
readily to relate nondegenerate background geometries to the~$\alpha$ and~$\beta$ surface reductions.
One might hope to obtain closer links by considering intermediate null reductions, in which the radical of~$VM$ is both
proper and nontrivial.
However, this is beyond the scope of this paper.

\section{Interlude: the Dif\/f(1) Hitchin equation}
\label{s:s1Hit}

\subsection*{HyperCR Einstein--Weyl spaces}

An important class of Einstein--Weyl spaces are the \emph{hyperCR} Einstein--Weyl spaces~\cite{GaTo:hms}.
In~\cite{Tod:sew}, Paul Tod presented a~way of reducing the hyperCR Einstein--Weyl equation to a~single second-order
dif\/ferential equation for a~complex function of three variables.
Since the equation is expected to be integrable, he posed the problem of identifying it.
In this section, I will follow a~similar approach to Tod and identify the hyperCR Einstein--Weyl equation with the
$\Diff(S^1)$ Hitchin equation on a~trivial spinor-vortex space.

An Einstein--Weyl structure $(\mathsf{c},D)$ on $B$ is said to be~\emph{hyperCR} if it admits an orthonormal frame
$\chi_1$, $\chi_2$, $\chi_3$ for the weightless (co)tangent bundle $L  T^{*} B\cong L^{-1} TB$ such that $D\chi_i=\kappa {*\chi_i}$ for some section $\kappa$ of $L^{-1}$ and each~$i$~-- indeed any~$\chi$ in
the unit sphere generated by $\chi_1$, $\chi_2$, $\chi_3$ satisf\/ies the same equation; these $\chi$'s are called the
\emph{hyperCR congruences} of $B$.

Since the Weyl connection is torsion-free and conformal, it is easy to see that these equations are implied by their
skew parts
\begin{gather}
\label{eq:hcrEW1}
{\rm d}^D\chi_1 = 2\kappa \chi_2\wedge\chi_3,
\\
{\rm d}^D\chi_2 = 2\kappa \chi_3\wedge\chi_1,
\\
{\rm d}^D\chi_3 = 2\kappa \chi_1\wedge\chi_2.
\label{eq:hcrEW3}
\end{gather}
Tod noticed that $\chi_i$ satisfying~\eqref{eq:hcrEW1}--\eqref{eq:hcrEW3} determine the Einstein--Weyl space: the
conformal metric is $\chi_1^2+\chi_2^2+\chi_3^2$ and the Einstein--Weyl equation follows
from~\eqref{eq:hcrEW1}--\eqref{eq:hcrEW3}~-- see also~\mbox{\cite{CaPe:sdc,GaTo:hms}}.
Now introduce a~gauge $\mu$ with $D\mu=\omega\mu$ and def\/ine $1$-forms $\alpha_i=2\mu^{-1}\chi_i$.
Then the equations~\eqref{eq:hcrEW1}--\eqref{eq:hcrEW3} may be rewritten in the form given by Tod~\cite{Tod:sew}:
\begin{gather}
\label{eq:hcrg1}
{\rm d}\alpha_1=-\omega\wedge\alpha_1+\kappa\alpha_2\wedge\alpha_3,
\\
{\rm d}\alpha_2=-\omega\wedge\alpha_2+\kappa\alpha_3\wedge\alpha_1,
\\
{\rm d}\alpha_3=-\omega\wedge\alpha_3+\kappa\alpha_1\wedge\alpha_2.
\label{eq:hcrg3}
\end{gather}
These equations are easier to interpret after complexif\/ication, so that the conformal structure is determined by its
null lines, which form a~bundle of conics in $P(T^{*} B)$.
This bundle is trivial, since $T^{*} B$ is trivialized by $\alpha_1$, $\alpha_2$, $\alpha_3$.
The pullback of the tautological $1$-form on $P(T^{*} B)$ by a~constant section is a~constant linear combination of
$\alpha_1$, $\alpha_2$ and $\alpha_3$, which is null for a~section of the bundle of conics.
Then~\eqref{eq:hcrg1}--\eqref{eq:hcrg3} are equivalent to the integrability of the distributions def\/ined by these null
$1$-forms, i.e., to the integrability of a~rank~$2$ distribution~${\mathcal{H}}$ on the bundle of conics.
The integral surfaces are the null surfaces which motivated Cartan~\cite{Car:cew} to study $3$-dimensional
Einstein--Weyl geometry, and the quotient of the bundle of conics by ${\mathcal{H}}$ is the \emph{minitwistor space}
${\mathcal{S}}$ of $B$~\cite{Hit:cme}.

On a~general Einstein--Weyl space, ${\mathcal{H}}$ is def\/ined by the Weyl connection: the hyperCR case is special in
that there is a~preferred trivialization of the bundle of conics with respect to which the distribution is horizontal.
Using this trivialization, the bundle of conics is $B\mathbin{{\times}}{\mathbb P}^1$, and after choosing a~projective
coordinate $\zeta$ on ${\mathbb P}^1$, the null $1$-forms are
\begin{gather*}
\alpha_\zeta = \alpha_1+{\boldsymbol i}\alpha_2 + 2\zeta \alpha_3-\zeta^2(\alpha_1-{\boldsymbol i}\alpha_2).
\end{gather*}
The system~\eqref{eq:hcrg1}--\eqref{eq:hcrg3} is equivalent to $\alpha_\zeta\wedge {\rm d}\alpha_\zeta=0$ for all
$\zeta$.
If $X_1$, $X_2$, $X_3$ is the dual frame, this means that the vector f\/ields $X_1+{\boldsymbol i} X_2+\zeta X_3$ and
$X_3-\zeta(X_1-{\boldsymbol i} X_2)$ span an integrable distribution (tangent to the null surfaces) for each f\/ixed~$\zeta$.
This is the hyperCR analogue of the Mason--Newman--Dunajski--Joyce description of hypercomplex structures.

To see explicitly what this means, put $\varphi=\alpha_1+{\boldsymbol i}\alpha_2$ and $\omega=-\tau\alpha_3+\gamma$,
with $\ip{\gamma,\alpha_3}=0$, so that~\eqref{eq:hcrg1}--\eqref{eq:hcrg3} become
\begin{gather}
\label{eq:hcrh}
{\rm d}\varphi=\bigl((\tau+{\boldsymbol i}\kappa)\alpha_3-\gamma\bigr)\wedge\varphi,
\\
{\rm d}\alpha_3=-\gamma\wedge\alpha_3+\tfrac{{\boldsymbol i}}2\kappa \varphi\wedge\overline\varphi.
\label{eq:hcrc}
\end{gather}
The f\/irst equation implies the integrability of the distribution def\/ined by $\varphi$.
Tod~\cite{Tod:sew} uses~\eqref{eq:hcrh} to put $\varphi=w{\rm d} z$.
I will not repeat this here.
Instead I want to analyse these equations from the point of view of integrable background geometries.
The idea is that the foliation determined by~$\chi_3$ is a~generalized dimensional reduction.
To see this, break the equations into horizontal and vertical parts by writing $\alpha_3={\rm d} t+A$ for some
f\/ibre coordinate~$t$ so that~$A$ is horizontal (i.e., in the span of~$\alpha_1$ and~$\alpha_2$).
Then for any $1$-form~$\beta$,
\begin{gather*}
{\rm d}\beta=d_N\beta+{\rm d} t\wedge\dot\beta=d_N\beta-A\wedge\dot\beta +\alpha_3\wedge\dot\beta,
\end{gather*}
where $d_N\beta$ is a~multiple of $\alpha_1\wedge\alpha_2$, $\dot\beta=\partial_t\beta$,
$\alpha_1(\partial_t)=0=\alpha_2(\partial_t)$, and ${\rm d} t(\partial_t)=1$.
Hence~\eqref{eq:hcrh}--\eqref{eq:hcrc} become
\begin{gather*}
d_N\varphi-A\wedge\dot\varphi=-\gamma\wedge\varphi,
\qquad
\dot\varphi=(\tau+{\boldsymbol i}\kappa)\varphi,
\\
d_NA -A\wedge\dot A=\tfrac{{\boldsymbol i}}2\kappa\varphi\wedge\overline\varphi,
\qquad
\dot A=\gamma.
\end{gather*}
The equations on the right simply def\/ine $\tau+{\boldsymbol i}\kappa$ and $\gamma$, so after computing that
$\dot\varphi\wedge\overline\varphi -\varphi\wedge\dot{\overline\varphi} =2{\boldsymbol
i}\kappa\varphi\wedge\overline\varphi$, the equations on the left reduce to
\begin{gather*}
d_N\varphi+\dot A\wedge\varphi-A\wedge\dot\varphi=0,
\\
d_NA+\dot A\wedge A=\tfrac14\bigl(\dot\varphi\wedge\overline\varphi -\varphi\wedge\dot{\overline\varphi}\bigr).
\end{gather*}
The conformal structure on $N$ has representative metric $\varphi\overline\varphi$ so the orientation can be chosen so
that $\varphi$ has type $(1,0)$.
It is now easy to see that these equations are Hitchin's equations with gauge group $\Diff(S^1)$: $\varphi$ and $A$ are
$1$-forms on $N$ with values in $\lie{g}=\Vect(S^1)$; $\varphi$ is the Higgs f\/ield and $A$ is the connection $1$-form,
satisfying
\begin{gather*}
F^A=[\varphi,\overline\varphi]_{\lie{g}}
\qquad
\overline\partial{}^{A}\varphi=0.
\end{gather*}
These are equivalent to the~\eqref{eq:hcrh}--\eqref{eq:hcrc} and hence to~\eqref{eq:hcrg1}--\eqref{eq:hcrg3}: since
$\dot\varphi$ has type $(1,0)$, one can locally write $\dot\varphi=(\tau+{\boldsymbol i}\kappa)\varphi$ and thus def\/ine
$\alpha_1$, $\alpha_2$, $\alpha_3$, $\omega$.
\begin{thm}
Let $(A,\varphi)$ be a~solution of the $\Diff(S^1)$ Hitchin equations $($on a~trivial spinor-vortex space
$N)$.
Then the Einstein--Weyl space defined by $(A,\varphi)$ is hyperCR, and one of its hyperCR congruences defines the
foliation over $N$.
Any hyperCR Einstein--Weyl space arises in this way.
\end{thm}
Tod's simplif\/ication of the $\Diff(S^1)$ Hitchin equations amounts to a~f\/ixing the $\Diff(S^1)$ gauge via
$e^{2{\boldsymbol i} t}=\varphi/\overline\varphi$.

{\sloppy In order to obtain new examples of hyperCR Einstein--Weyl spaces, consider Hitchin f\/ields where the gauge group is
a~f\/inite-dimensional subgroup of $\Diff(S^1)$.
The Abelian gauge group~$\Un(1)$ yields only f\/lat Einstein--Weyl spaces, but the af\/f\/ine and projective groups,
Af\/f$({\mathbb{R}})$ and~PSL$(2,{\mathbb{R}})$ are more interesting.
The former is more tractable, since Af\/f$({\mathbb{R}})$ is a~solvable group, meaning that the nonlinear Hitchin
equations can be solved by integrating a~sequence of linear equations.
Indeed, writing $\varphi=\varphi_0+\varphi_1 t$ and $A=A_0+A_1 t$ gives
\begin{gather*}
\begin{split}
& d \varphi_1=0,
\qquad
d \varphi_0+A_1\wedge \varphi_0-A_0\wedge \varphi_1=0,
\\
& d A_1=0,
\qquad
d A_0+A_1\wedge A_0=\tfrac14(\varphi_1\wedge\overline\varphi_0 -\varphi_0\wedge \overline\varphi_1).
\end{split}
\end{gather*}
Locally, the af\/f\/ine gauge freedom can be used to eliminate the linear term $A_1$ of the connection~$A$, while the
conformal gauge freedom can be used to make the linear term~$\varphi_1$ of the Higgs f\/ield~$\varphi$ equal to~$\lambda{\rm d} z$ with $\lambda$ constant (without loss of generality $\lambda=1$ unless it vanishes, which is
the Abelian case).
Let $\varphi_0= {\boldsymbol i} f {\rm d} z$ so that the equations reduce to
\begin{gather*}
{\boldsymbol i} f_{\bar z} {\rm d}\bar z\wedge {\rm d} z = \lambda A_0\wedge {\rm d} z,
\qquad
d A_0=-\tfrac12 {\boldsymbol i} f{\rm d} z\wedge {\rm d}\bar z.
\end{gather*}
For $\lambda\neq 0$, the general solution up to gauge transformation is therefore determined by a~real function
$f(z,\bar z)$ satisfying $\Delta f+2\lambda^2 f=0$: $\varphi=(t+{\boldsymbol i} f){\rm d} z$ and
$A=-{*{\rm d} f}/\lambda$.
Hence (local) eigenfunctions of the Laplacian on ${\mathbb{R}}^2$ give rise to hyperCR Einstein--Weyl spaces.

}

\subsection*{Einstein--Weyl spaces with a~geodesic generalized symmetry}

A \emph{shear-free geodesic congruence} on an Einstein--Weyl space is a~weightless unit vector f\/ield
$\chi\in\mathrm{C}^\infty(B,L^{-1}  TB)$ such that
\begin{gather*}
D\chi=\tau(\iden-\chi\otimes\chi)+\kappa {*\chi}
\end{gather*}
and for sections $\tau$, $\kappa$ of $L^{-1}$ called the \emph{divergence} and \emph{twist} of $\chi$.
On a~hyperCR Einstein--Weyl space, the hyperCR congruences are examples: they are also divergence-free and in fact this
characterizes them.
In the previous subsection it was found that the foliation def\/ined by a~hyperCR congruence is a~generalized symmetry, so
it is natural to ask, more generally, \textit{when does a~shear-free geodesic congruence define a~generalized symmetry?}
Note that the weightless unit vector f\/ield tangent to a~generalized symmetry over a~spinor-vortex space is always
shear-free, so this question can be rephrased: \textit{when is a~generalized symmetry geodesic?} Since the
Einstein--Weyl space $B$ is completely explicit~\eqref{eq:SVtoEW} in terms of the Hitchin f\/ield on the spinor-vortex
space $N$, this question is easily answered: the generalized symmetry is geodesic if and only if ${\mathcal{C}}$ (and
$\widetilde{\mathcal{C}}$, which is the complex conjugate in the Euclidean case) vanishes.
The Einstein--Weyl structure is then given by
\begin{gather}
g= 4\Phi\overline\Phi {\rm d} z{\rm d}\bar z +({\rm d} t + \alpha{\rm d} z +
\overline\alpha {\rm d}\bar z)^2,
\nonumber
\\
\omega=\dot\alpha{\rm d} z+\dot{\overline\alpha}{\rm d}\bar z
-\frac12\biggl(\frac{\psi+\dot\Phi}{\Phi}+ \frac{\overline\psi+\dot{\overline\Phi}}{\overline\Phi}\biggr)
({\rm d} t + \alpha{\rm d} z + \overline\alpha {\rm d}\bar z).
\label{eq:ewggs}
\end{gather}
If $\psi$ vanishes, then the spinor-vortex space is (locally) trivial and this is the case of the previous subsection.
Otherwise, $\psi$ is a~holomorphic trivialization of ${\mathcal{W}}$ on the open set where it is nonzero, and the
spinor-vortex space is given by a~spherical metric on $N$.

Two special cases of this construction have already been studied: the case that $\partial_t$ is a~genuine
symmetry~\cite{CaPe:sdc}, and the case that the congruence is also twist-free~\cite{CaTo:emh}.
The f\/irst class arises by supposing that the gauge group reduces to $\Un(1)$ or ${\mathbb{R}}$ and $B$ is said to be
Einstein--Weyl \emph{with a~geodesic symmetry}.
The Abelian Hitchin equations are easily solved on the spherical spinor-vortex space yielding the explicit formula
\begin{gather}
g=|h|^{2} g_{S^2}+\beta^2,
\qquad
\omega=-\frac{{\boldsymbol i}(h-\overline h)}{2|h|^2}\beta,
\qquad
{\rm d}\beta=\tfrac12(h+\overline h)\vol_{S^2},
\label{eq:ewgs}
\end{gather}
where $h$ is a~holomorphic function on an open subset of $S^2$.

The second class, the \emph{hyperCR-Toda spaces}, are obtained by from some explicit solutions of the af\/f\/ine Hitchin
equations: the twist-free condition reduces the gauge group to Af\/f$({\mathbb{R}})$; the connection is f\/lat, the linear
part of the Higgs f\/ield is constant, while the translational part of the Higgs f\/ield is given by a~holomorphic function~$h$.
The resulting Einstein--Weyl structure is{\samepage
\begin{gather}
g=(t+h)(t+\overline h)g_{S^2}+{\rm d} t^2,
\qquad
\omega=-\frac{2t+h+\overline h}{(t+h)(t+\overline h)}{\rm d} t.
\label{eq:hcrtoda}
\end{gather}
In this case the twist of the geodesic generalized symmetry $\partial_t$ vanishes.}

Note that the general theory of this paper justif\/ies the f\/inal remarks of~\cite{CaTo:emh}, by explaining the sense in
which the spherical metric is the natural quotient geometry of both structures~\eqref{eq:ewgs} and~\eqref{eq:hcrtoda}.

In both cases, the Einstein--Weyl space is hyperCR, although the congruence generated by~$\partial_t$ is no longer one
of the hyperCR congruences.
The same holds in general.
\begin{thm}
An Einstein--Weyl space with a~geodesic generalized symmetry is hyperCR.
\end{thm}
\begin{proof}
Gauduchon and Tod~\cite{GaTo:hms} show that an Einstein--Weyl space $(B,\mathsf{c},D)$ is hyperCR with twist
$\hat\kappa$ if and only if $*D\hat\kappa=\frac12F^D$ and $\hat\kappa^2=\frac16\scal^D$.

On an Einstein--Weyl space with a~geodesic generalized symmetry, direct computation of ${\rm d}^D\chi$ in the
gauge $(g,\omega)$ for the Einstein--Weyl structure~\eqref{eq:ewggs} yields
\begin{gather*}
\kappa=\frac1{4{\boldsymbol i}}\biggl(\frac{\dot\Phi-\psi}{\Phi}-
\frac{\dot{\overline\Phi}-\overline\psi}{\overline\Phi}\biggr).
\end{gather*}
Now put
\begin{gather*}
\hat\kappa=\frac1{2{\boldsymbol i}}\biggl(\frac{\dot\Phi}{\Phi}- \frac{\dot{\overline\Phi}}{\overline\Phi}\biggr)-\kappa
=\frac1{4{\boldsymbol i}}\biggl(\frac{\dot\Phi+\psi}{\Phi}-
\frac{\dot{\overline\Phi}+\overline\psi}{\overline\Phi}\biggr).
\end{gather*}
In the gauge $(g,\omega)$, the equation $*({\rm d}\hat\kappa-\omega\hat\kappa)=\frac12{\rm d}\omega$ is
a~straightforward though tedious computation.
The equation $\hat\kappa^2=\frac16\scal^D$ also follows by direct computation, although the calculation is greatly
simplif\/ied by using the general formula $-\frac16\scal^D=D_\chi\tau+\tau^2-\kappa^2$ for the scalar curvature of an
Einstein--Weyl space with a~shear-free geodesic congruence $\chi$~\cite{CaPe:sdc, PeTo:3ew}.
In this case $-2\tau=\psi/\Phi+\overline\psi/\overline\Phi$: the verif\/ication that
$\dot\tau-\omega(\partial_t)\tau+\tau^2 =\kappa^2-\hat\kappa^2$ is now easy.
\end{proof}
Conversely, twistor methods show that the foliation def\/ined by any shear-free geodesic cong\-ruen\-ce on any hyperCR
Einstein--Weyl space is a~generalized dimensional reduction over a~trivial or spherical spinor-vortex geometry, although
this is not easy to see by direct computation.

\section{Riccati spaces}
\label{s:rs}

Riccati spaces form a~foundation on which higher-dimensional geometries can be built, so although the matrix Riccati
equation is easy to solve, it is invaluable to understand the solutions carefully.
For this reason, I will begin by tackling the Riccati equation in an invariant way, without choosing a~conformal
trivialization of ${\mathcal{E}}$ or a~coordinate on $C$.
Recall that $\wedge^3{\mathcal{E}}=(TC)^3$, and it will be convenient to write $TC=L$, although the choice of
orientation implicit in this identif\/ication is not essential.

The matrix ${\mathcal{B}}$ is a~section of $L^{-1} \Sym_0{\mathcal{E}}$ and so, at each point of $C$,
it has two obvious invariants, of weight $-2$ and $-3$ respectively: $x=\frac23\trace({\mathcal{B}}^2)$ and
$y=4\det{\mathcal{B}}$, normalized so that the characteristic polynomial of ${\mathcal{B}}$ is $4\lambda^3-3x\lambda-y$.
The discriminant of this polynomial is the section $y^2-x^3$ of $L^{-6}$: more precisely, writing $4(y^2-x^3)=-27c^2$
yields $c=\psi_1 \psi_2 \psi_3$, where $\psi_1=\frac23(\lambda_2-\lambda_3)$, $\psi_2=\frac23(\lambda_3-\lambda_1)$ and
$\psi_3=\frac23(\lambda_1-\lambda_2)$, $\{\lambda_i\}$ being the eigenvalues of ${\mathcal{B}}$; the sign of $c$ depends
on the ordering of the eigenvalues.
Note that
\begin{gather*}
\psi_1+\psi_2+\psi_3=0,
\qquad
\psi_1^2+\psi_2^2+\psi_3^2=2x
\end{gather*}
{and}
\begin{gather*}
(\psi_1-\psi_2)(\psi_2-\psi_3)(\psi_3-\psi_1)=-8\lambda_1\lambda_2\lambda_3=-2y.
\end{gather*}

\subsection*{Pencils of conics}

In order to interpret this geometrically, complexify $C$ and ${\mathcal{E}}$ so that the
conformal structure $\mathsf{c}$ on ${\mathcal{E}}$ is determined by its null lines, which form a~bundle of conics
${\mathcal{S}}({\mathcal{E}})$ in the bundle $P({\mathcal{E}})$ of projective planes over $C$.
The role of ${\mathcal{B}}$ is to determine a~pencil of conics in each f\/ibre, associated to the two-dimensional family
of bilinear forms $\mathsf{c}\circ(s\iden+t{\mathcal{B}})$: $\mathsf{c}\circ{\mathcal{B}}$ is distinguished by being
traceless with respect to the f\/ixed bilinear form~$\mathsf{c}$.
Now, two distinct conics meet in four points, counted with multiplicity, so there are six possible conf\/igurations: the
generic case (I), where the four points are distinct; the four degenerations (II, III, D, N),
when two, three, two pairs, or four points come together; and the trivial case~(0), when ${\mathcal{B}}=0$ and the `pencil' is constant.
$$
\includegraphics{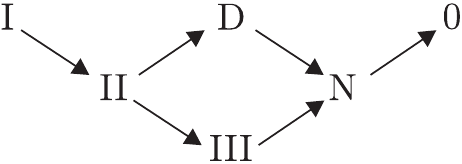}
$$

The notation here follows the well-known application of this classif\/ication to Weyl tensors in four dimensions.
Only types (I, D, 0) are compatible with Euclidean reality conditions: in this case ${\mathcal{B}}$ must be
diagonalizable, and the eigenvalues are either distinct $(\lambda_1,\lambda_2,\lambda_3)$ with
$\lambda_1+\lambda_2+\lambda_3=0$, or of the form $(\lambda,\lambda,-2\lambda)$, or all zero.

A natural way to analyse a~pencil of conics is to identify one of the conics with a~projective line~$\mathbb{P}^1$.
In the present situation, this is done by introducing a~bundle of spinors ${\mathcal{U}}$ for ${\mathcal{E}}$, i.e.,
with ${\mathcal{E}}=S^2{\mathcal{U}}$, so that ${\mathcal{S}}({\mathcal{E}})=P({\mathcal{U}})$ and ${\mathcal{U}}$
inherits a~connection from ${\mathcal{E}}$.
Trivializing ${\mathcal{U}}$ using this connection identif\/ies ${\mathcal{S}}({\mathcal{E}})$ locally with
$C\mathbin{{\times}}\mathbb{P}^1$, on which an af\/f\/ine coordinate $r$ for $C$ and a~projective coordinate $\zeta$ for
$\mathbb{P}^1$ may be introduced.
The isomorphism between sections of the inner product bundle $L^{-1} {\mathcal{E}}$ and vertical
vector f\/ields on ${\mathcal{S}}({\mathcal{E}})$ may be described concretely using the induced parallel orthonormal frame
for $L^{-1} {\mathcal{E}}$: sections induce vector f\/ields via their inner product with the
tautological null vector f\/ield $e_\zeta\otimes\partial_\zeta$ where
\begin{gather}
\label{eq:tnvf}
e_\zeta=\bigl(\tfrac12\big(\zeta^2+1\big),{\boldsymbol i}\zeta,\tfrac{{\boldsymbol i}}2\big(\zeta^2-1\big)\bigr).
\end{gather}
The base locus of the pencil of conics (at each point of $C$) consists of the four zeros of the quartic polynomial
$\ip{{\mathcal{B}}(e_\zeta),e_\zeta}\otimes\partial_\zeta^2$, which is a~section of $L^{-1}\otimes{\mathcal{O}}(4)$.
It is now straightforward and entirely classical~\cite{Todd:pag} to analyse the zeros of this quartic in the various
cases, and hence relate the six types of pencil to properties of the matrix ${\mathcal{B}}$.
In particular, the generic case is given by $c\neq0$, when the four points are distinct.

\begin{rem}
Four distinct points on $\mathbb{P}^1$ are determined up to projective transformation by their cross-ratio, and this
freedom is often used to identify the points with $0$, $1$, $\infty$, $t$, with $t$ being the cross-ratio.
However, as remarked by Yoshida in his wonderful book~\cite{Yos:hml}, ``it is not fair that only the fourth point is
allowed to move freely''.
The democratic cross-ratio used there is the point $[\psi_1,\psi_2,\psi_3]$ on the line $\psi_1+\psi_2+\psi_3=0$ in~$\mathbb{P}^2$.
\end{rem}

\subsection*{Solution of the Riccati equation}

I shall now f\/ind the solutions of the Riccati equation $D{\mathcal{B}}=2({\mathcal{B}}^2)_0$.
Here and elsewhere $D$ denotes dif\/ferentiation with respect to the af\/f\/ine structure on $C$, and takes a~section of
a~natural bundle~$F$ to a~section of $L^{-1}  F$.
In naive terms, $D$ is dif\/ferentiation with respect to an af\/f\/ine coordinate~$r$ (i.e., ${\rm d} r$ is
a~$D$-parallel section of $L^{-1}=T^{*} C$, and $L^{-1}$ is trivialized using this section).

The Riccati equation shows that $D{\mathcal{B}}$ commutes with ${\mathcal{B}}$, and so dif\/ferentiating the
Cayley--Hamilton equation $4{\mathcal{B}}^3-3x{\mathcal{B}}-y\iden=0$, gives
\begin{gather*}
\begin{split}
& 0=12 (D{\mathcal{B}}){\mathcal{B}}^2-3xD{\mathcal{B}}-3Dx{\mathcal{B}}-Dy\iden
=24{\mathcal{B}}^4-18x{\mathcal{B}}^2+3x^2\iden-3Dx{\mathcal{B}}-Dy\iden
\\
& \phantom{0}
=3(2y-Dx){\mathcal{B}}+(3x^2-Dy)\iden.
\end{split}
\end{gather*}
Hence $Dx=2y$, $Dy=3x^2$, $D^2x=6x^2$ and so $D(y^2-x^3)=0$.
In other words, the discriminant~$c^2$ is constant ($D$-parallel).
The constancy of $c$ means that if~${\mathcal{B}}$ has distinct eigenvalues at a~point, it has distinct eigenvalues
everywhere, in which case it is diagonalizable, and can be assumed diagonal, since $D{\mathcal{B}}$ commutes with~${\mathcal{B}}$.
I consider this f\/irst.
\begin{itemize}\itemsep=0pt
\item[I.] \textit{${\mathcal{B}}$ has distinct eigenvalues}.
\end{itemize}

Taking ${\mathcal{B}}$ diagonal, it is straightforward to compute the equation for the eigenvalues:
$D\lambda_1=2\lambda_1^2-\frac23(\lambda_1^2+\lambda_2^2+\lambda_3^2)$ and similarly for $\lambda_2$ and $\lambda_3$.
Hence $D\psi_1=\psi_1(\psi_2-\psi_3)=-2\lambda_1\psi_1$ and similarly for $\psi_2$ and $\psi_3$.
Now let $\chi_1=-\psi_2\psi_3$ etc., so that $\chi_1\chi_2\chi_3=-c^2$ and
$D\chi_1=-\chi_1(\psi_2-\psi_3)=2\lambda_1\chi_1$.
Squaring this gives $(D\chi_1)^2=\chi_1^2 (\psi_2^2+\psi_3^2-2\psi_2\psi_3)=\chi_1^2(\psi_1^2-4\psi_2\psi_3)
=4\chi_1^3+c^2$, since $\psi_1+\psi_2+\psi_3=0$.
Hence the $\chi_i$ all satisfy the equation
\begin{gather*}
(D\chi)^2=4\chi^3+c^2,
\\
\intertext{which is closely related to the equation for $x$, namely} (Dx)^2=4x^3-27c^2.
\end{gather*}
Thus $x$ and the $\chi$'s are equianharmonic Weierstra\ss\;elliptic functions of an af\/f\/ine coordinate, the period
lattice for $x$ being the barycentric subdivision of the lattice for the $\chi$'s.

\begin{rem}
The advantage of an invariant description is the f\/lexibility in the choice of coordinates: one does not have to use the
af\/f\/ine coordinate $r$.
Indeed, in Section~\ref{s:bm}, a~dif\/ferent gauge choice was motivated: there is a~natural projective structure
$D^2+\frac16\trace{\mathcal{B}}^2=D^2+\frac14 x$ and with respect to a~projective coordinate $t$, $D=\partial_t+a$,
where $\dot a-\frac12 a^2=-\frac12 x$.
Now $\dot x=2(ax+y)$ and $\dot y=3(ay+x^2)$, so that $\ddot a=a\dot a-\frac12 \dot x=a\dot a-ax-y$ and
\begin{gather*}
\dddot a= a \ddot a+\dot a^2-\dot a x-\tfrac52a\dot x+\bigl(\tfrac32a\dot x-\dot y\bigr)
= a \ddot a+\dot a^2-\dot a (a^2-2\dot a)-5a(a\dot a-\ddot a)+3(a^2-x)x
\\
\phantom{\dddot a}
=6a\ddot a+3\dot a^2-6a^2\dot a+6\dot a(a^2-2\dot a)
=6a\ddot a-9\dot a^2,
\end{gather*}
which is the \emph{Chazy equation}.
Of course it is well-known that the Chazy equation arises in the study of selfdual metrics in a~scalar-f\/lat gauge.

The generic solution, given by elliptic functions of the af\/f\/ine coordinate, may instead be presented in terms of modular
functions of the projective coordinate~\cite{AbCl:sis,Dub:g2d}.
\end{rem}

I turn next to the degenerate cases, when $c=0$, i.e., $y^2=x^3=64\lambda^6$, where the eigenvalues of ${\mathcal{B}}$
are $\lambda,\lambda$ and $-2\lambda$.
Assume f\/irst that $\lambda$ is not identically zero.
Then, on an open set at least, the (generalized) eigenspaces of ${\mathcal{B}}$ are constant (again using the fact that
$D{\mathcal{B}}$ commutes with ${\mathcal{B}}$) and hence ${\mathcal{B}}$ can be assumed to take the the form
\begin{gather*}
{\mathcal{B}}=\left[
\begin{matrix}
\lambda+\mu&{\boldsymbol i}\mu& 0
\\
{\boldsymbol i}\mu&\lambda-\mu&0
\\
0&0&-2\lambda
\end{matrix}
\right].
\end{gather*}
The matrix Riccati equation now yields
\begin{gather*}
D\lambda=-2\lambda^2,
\qquad
D\mu=4\lambda\mu.
\end{gather*}
Since $\lambda$ is not identically zero, it is given (up to translation) by $\lambda=1/2r$.
This is def\/ined for $r$ nonzero, and $\lambda$ has no zeros.
In this coordinate $\mu=b r^2$ for constant $b$, which is either identically zero, or nonzero for all nonzero $r$.
There are thus two cases.{\samepage
\begin{itemize}\itemsep=0pt
\item[II.] \textit{$y^2=x^3$ is nowhere zero, and ${\mathcal{B}}$ is nowhere diagonalizable}.
\item[D.] \textit{$y^2=x^3$ is nowhere zero, and ${\mathcal{B}}$ is everywhere diagonalizable}.
\end{itemize}
The notation is justif\/ied by noting that $\ip{{\mathcal{B}}(e_\zeta),e_\zeta}$ is the quartic
$\frac14(\zeta-1)^2(\mu(\zeta-1)^2+3\lambda(\zeta+1)^2)$: $\zeta=1$ is a~repeated root, as is $\zeta=-1$ when $\mu=0$.}

\begin{rem}
\label{r:Dproj}
The type (D) solution is even simpler with respect to a~projective coordinate.
Since $2{\mathcal{B}}^2_0=-2\lambda{\mathcal{B}}$, setting $D=\partial_t+2\lambda$ gives $\partial_t{\mathcal{B}}=0$, so
that $\lambda$ is constant in this gauge and $a=2\lambda$ satisf\/ies $\dot a-\frac12
a^2=-2\lambda^2=-\frac13\trace{\mathcal{B}}^2=-\frac12 x$; thus $t$ is a~projective coordinate with respect to the
natural projective structure, and ${\mathcal{B}}$ and $a$ are constant.
\end{rem}

It remains to consider the case that $y^2=x^3$ is identically zero, i.e., the eigenvalues of ${\mathcal{B}}$ are all
zero, so that ${\mathcal{B}}^3=0$.
If ${\mathcal{B}}^2=0$ then ${\mathcal{B}}$ is constant, so that in any case the kernel and image of ${\mathcal{B}}$ are
constant; ${\mathcal{B}}$ may thus be assumed to take the form
\begin{gather*}
{\mathcal{B}}=\left[
\begin{matrix}
\mu& {\boldsymbol i}\mu& \beta
\\
{\boldsymbol i}\mu&-\mu&{\boldsymbol i}\beta
\\
\beta&{\boldsymbol i} \beta&0
\end{matrix}
\right]
\end{gather*}
and the matrix Riccati equation yields
\begin{gather*}
D\beta=0,
\qquad
D\mu=2\beta^2.
\end{gather*}
Hence $\beta$ is a~constant, zero if and only if ${\mathcal{B}}$ is constant, and if $\beta$ is nonzero, then in
a~suitably translated af\/f\/ine coordinate $\mu=2\beta^2r$.
Thus there are three more cases.
\begin{itemize}\itemsep=0pt
\item[III.] \textit{${\mathcal{B}}^3=0$ and ${\mathcal{B}}^2$ is nowhere zero}.
\item[N.] \textit{${\mathcal{B}}^2=0$ and ${\mathcal{B}}$ is constant and nonzero}.
\item[0.] ${\mathcal{B}}=0$.
\end{itemize}
Again the notation is justif\/ied by computing $\ip{{\mathcal{B}}(e_\zeta),e_\zeta}=
\frac14(\zeta-1)^3(\mu(\zeta-1)+2{\boldsymbol i} \beta(\zeta+1))$.

\subsection*{Isomonodromic deformations} Recall that, for
$\Phi\in\mathrm{C}^\infty(C,{\mathcal{E}}^{*}\otimes\lie{g}_C)$, the generalized Nahm equation~\eqref{eq:nahm} is
\begin{gather*}
D\Phi-{*[\Phi,\Phi]_{\lie{g}}}={\mathcal{B}}\mathinner{\cdot}\Phi.
\end{gather*}
I claim that this equation, for ${\mathcal{B}}$ not identically zero, is equivalent to the fact that
\begin{gather*}
{\rm d}+\frac{\Phi}{{\mathcal{B}}}={\rm d}+\frac{\Phi(e_\zeta)
{\rm d}\zeta}{\ip{{\mathcal{B}}(e_\zeta),e_\zeta}}
\end{gather*}
is an isomonodromic family of connections on $\CP1$ parameterized by $C$, with four poles along the base locus of the
pencil of conics, where I use the tautological null vector f\/ield $e_\zeta$~\eqref{eq:tnvf} to clarify the meaning of the
connection $1$-form.
More precisely, this follows from the fact that the meromorphic connection
\begin{gather*}
{\rm d}+\frac{\Phi+\ip{*{\mathcal{B}},\Phi}}{{\mathcal{B}}} ={\rm d}+\frac{\Phi(e_\zeta)
{\rm d}\zeta+\ip{*e_\zeta\wedge{\mathcal{B}}(e_\zeta),\Phi}{\rm d} r}
{\ip{{\mathcal{B}}(e_\zeta),e_\zeta}},
\end{gather*}
def\/ined on the pullback of $\lie{g}_C$ to the bundle of conics ${\mathcal{S}}({\mathcal{E}})$ over the Riccati space
$C$, is f\/lat if and only if $\Phi$ satisf\/ies the generalized Nahm equation.

This is a~straightforward verif\/ication.
Equivalently, if the gauge algebra $\lie{g}_C$ is represented as a~Lie algebra of vector f\/ields, then this observation
can be reformulated as the integrability of a~rank two distribution on a~bundle over ${\mathcal{S}}({\mathcal{E}})$, and
the poles of the connection can be viewed as points of ${\mathcal{S}}({\mathcal{E}})$ over which the distribution is
tangent to the f\/ibres.
Explicitly, in the above coordinates, two vector f\/ields generating this distribution are
\begin{gather}
\ip{{\mathcal{B}}(e_\zeta),e_\zeta}\partial_\zeta+\Phi(e_\zeta),
\nonumber
\\
\partial_r+\ip{{\mathcal{B}}(e_\zeta),e_\zeta'}\partial_\zeta+\Phi(e_\zeta'),
\label{eq:NLP}
\end{gather}
where $e_\zeta'=(\zeta,{\boldsymbol i},{\boldsymbol i}\zeta)$.
One can easily check that the distribution is closed under Lie bracket.

An advantage of this vector f\/ield interpretation is that it also makes sense on the trivial Riccati space
${\mathcal{B}}=0$, when~\eqref{eq:NLP} reduces to the standard Lax pair for the Nahm equation.

\subsection*{Addendum: Riccati spaces and integrability\\ by the method of hydrodynamic reductions}

In~\cite{FHZ:cqim}, Ferapontov et al.~independently discover the integrability of a~system of equations equivalent to
the generalized Nahm equations on Riccati spaces.
Their study concerns equations of the form
\begin{gather*}
\frac{\partial^2 (A^{\alpha\beta}(u))}{\partial x^\alpha \partial x^\beta}=0,
\end{gather*}
where $u=u(x^1,x^2,x^3)$ is a~function of three variables, and $A^{\alpha\beta}$ is a~$3\times 3$ symmetric matrix of
functions of one variable.
(Here I use the summation convention on Greek indices; this may also be viewed as abstract index notation.) Expanding
one derivative yields the equivalent formulation
\begin{gather}
\label{eq:fhz1}
\frac{\partial}{\partial x^\alpha} \left( V^{\alpha\beta}(u) \frac{\partial u}{\partial x^\beta} \right) =0,
\end{gather}
where $V^{\alpha\beta}$ is the derivative of $A^{\alpha\beta}$ (i.e., $V=A'$).

The approach in~\cite{FHZ:cqim} characterizes integrability using the method of hydrodynamic reductions; this tests for
the existence of suf\/f\/iciently many multi-phase solutions $u=U(R^1,R^2,\ldots R^N)$, where the $R^j(x^1,x^2,x^3)$ are
solutions to arbitrarily many commuting $(1+1)$-dimensional systems of hydrodynamic type.
The result of this analysis is that~\eqref{eq:fhz1} is integrable in this sense if and
only if there is a~scalar function $k(u)$ of one variable such that
\begin{gather}
\label{eq:fhz2}
V''(u) = (\log\det V )'(u) V'(u) + k(u) V(u)
\end{gather}
or equivalently $(V'/\det(V))' = k V$.

The relation to the Riccati space equation is subtle, and involves introducing a~$u$-dependent triple of vectors
$\theta^i_\alpha(u)$ ($i=1,2,3$) such that
\begin{gather}
\label{eq:V}
V=(\det\theta)\theta^{-1}(\theta^{\scriptscriptstyle\mathrm T})^{-1},
\qquad
\text{i.e.,}
\qquad
V^{\alpha\beta}=\sum\limits_{i=1}^3(\det\theta)\big(\theta^{-1}\big)_i^\alpha\big(\theta^{-1}\big)_i^\beta,
\end{gather}
where $\theta=(\theta^i_\alpha)$ is assumed invertible.
Now write ${\mathcal{C}}_{ij}:=(\theta^{-1})_i^\alpha(\theta^j_{\alpha})'$ so that $\theta^i_\alpha =\sum\limits_{j=1}^3
{\mathcal{C}}_{ij}\theta^j_\alpha$.
There is a~gauge freedom in the choice of $\theta$, which may be used to suppose ${\mathcal{C}}$ is symmetric.

Substituting~\eqref{eq:V} into~\eqref{eq:fhz2}, straightforward computation, using the symmetry of ${\mathcal{C}}$
together with standard matrix identities such as $(A^{-1})'=-A^{-1}A'A^{-1}$ and $(\det A)^{-1}(\det A)'=\trace(A^{-1}
A')$, then yields
\begin{gather*}
{\mathcal{C}}' = 2 {\mathcal{C}}^2 - (\trace {\mathcal{C}}) {\mathcal{C}} + c I,
\end{gather*}
where $c$ is an unknown scalar function of $u$ (related to $k$).

Now decompose ${\mathcal{C}}={\mathcal{B}}+a I$ with ${\mathcal{B}}$ tracefree; then the trace part of the above
equation determines $c$, while the trace-free part is the Riccati equation for ${\mathcal{B}}$.
Thus, remarkably, the hydrodynamic integrability condition agrees with the twistor-theoretic Riccati equation.

After classifying solutions of~\eqref{eq:fhz2} (equivalently, solutions of the Riccati equation), Ferapontov et al.
compute the form of the equation~\eqref{eq:fhz1} for $u$ arising from each solution.
They then observe that the central quadric Ansatz
\begin{gather*}
x^\alpha M_{\alpha\beta}(u) x^\beta =1
\end{gather*}
for $u$ yields all Painlev\'e equations~\cite{FHZ:cqim}.

It is this geometry that originally led to the identif\/ication of~\eqref{eq:fhz2} with the Riccati equation.
The key observation is that~\eqref{eq:fhz1} is a~divergence form equation, meaning that $u$ is in the kernel of the
Laplace--Beltrami operator of the metric $g_{\alpha\beta}(u)$ with $g={\rm adj}(V)=(\det V)V^{-1}$ so that
$V=\sqrt{\det g} g^{-1}$.
Hence the equations in this class have the form
\begin{gather}
\label{eq:ibg}
{\rm d}{*_g {\rm d} u}=0,
\end{gather}
where $g=g(u)$, and $*_g$ is the associated Hodge star operator.
Now $\theta^i_\alpha(u)$ has an interpretation, as a~\emph{framing} of the metric $g$:
\begin{gather*}
g_{\alpha\beta}(u)=\sum\limits_{i=1}^3 \theta^i_{\alpha}(u)\theta^i_{\beta}(u).
\end{gather*}
In terms of the geometry of $g$, $\theta$ is thus an ${\mathbb{R}}^3$-valued $1$-form: on vector f\/ields $X,Y$, the
metric is $g(u)(X,Y)=\ip{\theta(u)X,\theta(u)Y}$.
Thus there is a~gauge freedom to rotate $\theta$ by an $\SO(3)$-valued function of $u$: this is the freedom used to make
${\mathcal{C}}_{ij}$ symmetric above (which in turn f\/ixes $\theta$ up to a~rigid rotation).

As will be explained in~Section~\ref{s:sdsn}
(the last part of the next and f\/inal interlude),~\eqref{eq:fhz1} is a~hodograph
transformation of the generalized Nahm equation with gauge group $\SDiff(\Sigma^2)$.
The central quadric Ansatz corresponds to a~reduction of gauge group from $\SDiff(\Sigma^2)$ to $\SU(2)$ (which acts by
area preserving dif\/feomorphisms on the $2$-sphere).
This in turn explains the appearance of the Painlev\'e equations: generalized Nahm equations with gauge group $G$
describe isomonodromic deformation problems for $G$-connections on ${\mathbb{C}} P^1$ with four poles, and, when
$G=\SU(2)$ (rank~$2$ bundles over ${\mathbb{C}} P^1$), the connection between such isomonodromy problems and the
Painlev\'e equations is well known~\cite{JMU:mpd}.

\section{Interlude: the Dif\/f(2) generalized Nahm equation}

This interlude is devoted to Weyl structures of the form~\eqref{eq:RicEW}.
It follows from the work of Section~\ref{s:bgf}, that these Weyl structures are Einstein--Weyl provided that
${\mathcal{B}}$ satisf\/ies the matrix Riccati equation and $\Phi^i=\phi^i\partial_p+\psi^i\partial_q$ ($i=1,2,3$) def\/ines
a~generalized Nahm f\/ield on this Riccati space, the gauge group being a~subgroup of $\Diff(\Sigma^2)$ where $\Sigma^2$
is a~surface with coordinates $p$, $q$.

However, as details were omitted in Section~\ref{s:bgf}, I present a~self-contained study of these Weyl structures in
the following two cases:
\begin{enumerate}\itemsep=0pt
\item[(i)]
${\mathcal{B}}=0$ (so the Riccati space is trivial);
\item[(ii)]
the $\Phi^i$ are area preserving vector f\/ields.
\end{enumerate}
These examples have an interest that stretches beyond the proscenium of Einstein--Weyl geo\-met\-ry, since they are closely
related both to well-known integrable systems, and also to hypercomplex, hyperk\"ahler  and scalar-f\/lat K\"ahler
structures in four dimensions.

\subsection{HyperCR Einstein--Weyl spaces and the Dif\/f(2) Nahm equation}

In Section~\ref{s:s1Hit}, hyperCR Einstein--Weyl spaces were related to the $\Diff(S^1)$ Hitchin equation.
Now it is known~\cite{Cal:sde, GaTo:hms} that any solution of the $\Diff(S^1)$ Einstein--Weyl Bogomolny equation on
a~hyperCR Einstein--Weyl space gives rise to a~hypercomplex $4$-manifold (which is sometimes hyperk\"ahler).

This two step construction ties in with the construction of hypercomplex and hyperk\"ahler  structures from the
$\Diff(\Sigma^2)$ Hitchin equation: see Section~\ref{s:hchk}.
On the other hand, also in Sec\-tion~\ref{s:hchk}, hypercomplex and hyperk\"ahler  structures were related to the
$\Diff(\Sigma^3)$ Nahm equation, so it is natural to expect that hyperCR Einstein--Weyl spaces may be constructed from
the $\Diff(\Sigma^2)$ Nahm equation.

In fact this is not hard to see.
Suppose that $D=D^g+\omega$ is a~Weyl structure given by
\begin{gather*}
g= \eta_1^2+\eta_2^2+\eta_3^2,
\qquad
\omega= -\frac{\sum\limits_{i,j,k}\varepsilon_{ijk}\nu_i(\phi^j_p+\psi^j_q)\eta_k} {\nu_1^2+\nu_2^2+\nu_3^2},
\end{gather*}
{where}
\begin{gather*}
\eta_i= \nu_i {\rm d} r + \phi^i {\rm d} q - \psi^i {\rm d} p,
\\
\nu_1= \phi^2\psi^3-\phi^3\psi^2,
\qquad
\nu_2= \phi^3\psi^1-\phi^1\psi^3,
\qquad
\nu_3= \phi^1\psi^2-\phi^2\psi^1,
\end{gather*}
{and}
\begin{gather*}
\Phi=\big(\phi^1,\phi^2,\phi^3\big)\partial_p+\big(\psi^1,\psi^2,\psi^3\big)\partial_q.
\end{gather*}
The Nahm equation expands to give
\begin{gather*}
\phi^1_r=\phi^2\phi^3_p-\phi^3\phi^2_p+\psi^2\phi^3_q-\psi^3\phi^2_q,
\qquad
\psi^1_r=\phi^2\psi^3_p-\phi^3\psi^2_p+\psi^2\psi^3_q-\psi^3\psi^2_q
\end{gather*}
and cyclic permutations.
From this it follows easily that ${\rm d}\eta_1+\omega\wedge\eta_1=\kappa \eta_2\wedge\eta_3$ (and similarly
for the cyclic permutations), where
\begin{gather*}
\kappa=\frac{\nu_1\big(\phi^1_p+\psi^1_q\big)+\nu_2\big(\phi^2_p+\psi^2_q\big) +\nu_3\big(\phi^3_p+\psi^3_q\big)}{\nu_1^2+\nu_2^2+\nu_3^2}.
\end{gather*}
Hence the Weyl structure is hyperCR Einstein--Weyl.
Note that if the Nahm f\/ield preserves the (\textit{a priori} arbitrary) area form ${\rm d} p\wedge{\rm d}
q$ then $\kappa=0$ and the Einstein--Weyl space is f\/lat.
Ward~\cite{Ward:sut} uses this fact to linearize the Nahm equation in this case.

In fact any hyperCR Einstein--Weyl space $B$ arises locally from this construction.
One way to see this is to choose an Abelian monopole on $B$.
This def\/ines a~hypercomplex structure with a~triholomorphic vector f\/ield $K$ tangent to the f\/ibres over $B$.
For example, $\kappa$ is an Abelian monopole on any hyperCR Einstein--Weyl space and the corresponding hypercomplex
structure is hyperk\"ahler  with a~triholomorphic homothetic vector f\/ield~\cite{GaTo:hms}.

Now let $r$ be a~solution of $\Delta r=0$, which is constant on the f\/ibres over $B$, where $\Delta$ is the Laplacian of
the Obata connection.
This def\/ines a~divergence-free coframe ${\rm d} r$, $I{\rm d} r$, $J{\rm d} r$, $K{\rm d} r$ and
a~dual frame of vector f\/ields $V_0$, $V_1$, $V_2$, $V_3$.
Since $r$ is constant on the f\/ibres over $B$, $K$ is tangent to the level surfaces of $r$, so it is in the span of
$V_1$, $V_2$, $V_3$ and commutes with them all (since it is triholomorphic).
Therefore the gauge group of the Nahm f\/ield def\/ined by~$V_1$, $V_2$, $V_3$ reduces to the group of dif\/feomorphisms of~$\Sigma^3$ commuting with the f\/low of~$K$.
Taking the local quotient by~$K$ gives a~Nahm f\/ield with gauge group~$\Diff(\Sigma^2)$ for some surface~$\Sigma^2$, and
this is clearly a~Nahm f\/ield giving rise to the Einstein--Weyl space~$B$.

\begin{thm}
Let $\Phi$ be a~solution of the Nahm equation with gauge group $\Diff(\Sigma^2)$ for some surface~$\Sigma^2$.
Then the Einstein--Weyl space defined by $\Phi$ is hyperCR, and any hyperCR Einstein--Weyl space arises locally in this way.
If there is a~reduction to~$\SDiff(\Sigma^2)$, the Einstein--Weyl space is flat.
\end{thm}

There is a~great deal more freedom in this construction than in the $\Diff(S^1)$ Hitchin construction of hyperCR
Einstein--Weyl spaces, which can be an advantage or a~drawback, since the same hyperCR Einstein--Weyl space will arise
in many ways.
There are some interesting special cases one could consider, such as the group $\SL(2,{\mathbb{C}})$ acting on~$\CP1$
(which does not preserve an area form).
The Nahm equations in this case reduce to a~(complexif\/ied) Euler top equation, solvable in terms of elliptic functions.

\subsection{The SDif\/f(2) generalized Nahm equation}
\label{s:sdsn}

In addition to hyperCR Einstein--Weyl spaces, another important class consists of the Einstein--Weyl spaces arising from
solutions of the $\SU(\infty)$ Toda f\/ield equation \mbox{$u_{xx}+u_{yy}+(e^u)_{zz}=0$} \cite{LeBr:cp2,Tod:p3,Ward:sut}.
The terminology originates by regarding the equation as a~dispersionless limit of the $\SU(N)$ Toda equation as
$N\to\infty$: the Toda equation is a~two-dimensional system with independent variables $(x,y)$; the $z$ variable is
normally discrete, but becomes continuous in the limit.
The large $N$ limit of $\SU(N)$ may be interpreted as the group $\SDiff(S^1\mathbin{{\times}} S^1)$ of area preserving
dif\/feomorphisms of a~torus, with its Lie algebra ${\mathbb{Z}}$-graded by the Fourier components in one of the circles.
There is a~potential source of confusion here, however, since the $\SU(\infty)$ Toda equation is related to area
preserving dif\/feomorphisms in another way.
Namely, each solution gives rise to a~hyperk\"ahler  metric with a~Killing vector $K$~\cite{BoFi:kve}, and hyperk\"ahler
metrics are in turn obtained from the $\SDiff(\Sigma^2)$ Hitchin equations on ${\mathbb{R}}^2$, as shown in
Theorem~\ref{th:s2Hit} (cf.~\cite{Ward:suc}).
Hence one expects the $\SU(\infty)$ Toda equation to be equivalent to a~symmetry reduction of the $\SDiff(\Sigma^2)$
Hitchin equations to one dimension.
$$
\includegraphics{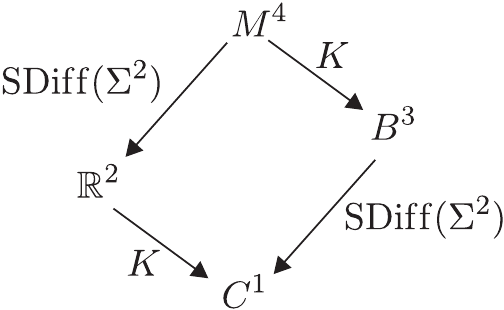}
$$

The Einstein--Weyl structure on a~Toda Einstein--Weyl space $B^3$ is given by the metric and Weyl $1$-form
\begin{gather*}
h=e^u\big({\rm d} x^2+{\rm d} y^2\big)+{\rm d} z^2,
\qquad
\omega_h=-u_z{\rm d} z.
\end{gather*}
In these geometric terms, the $\SU(\infty)$ Toda equation may be written ${\rm d}{*_h{\rm d} u}=0$, and so
$*_h{\rm d} u$ is a~closed $2$-form, which is therefore locally of the form ${\rm d} p\wedge{\rm d}
q$.
Now $(p,q,u)$ will be functionally independent unless $0={\rm d} u\wedge{*_h {\rm d}
u}={*_h|{\rm d} u|^2}$.
Hence (in Euclidean signature), $(p,q,u)$ may be used as coordinates, except in the trivial case ($u$ constant).
The solution of the Toda equation is now given implicitly by the functions $(x,y,z)$ of $(p,q,u)$.
These functions will turn out to satisfy a~generalized Nahm equation on the type (D) Riccati space with gauge group
$\SDiff(\Sigma^2)$.

In order to see this, and relate it to the work of Ferapontov et al.~\cite{FHZ:cqim}, consider the generalized Nahm
equation on any Riccati space with af\/f\/ine coordinate $r$, where the Higgs f\/ields
$\Phi^i=\phi^i\partial_p+\psi^i\partial_q$ preserve the area form ${\rm d} p\wedge{\rm d} q$ on a~surface
$\Sigma^2$ with coordinates $(p,q)$.
Motivated by the general construction of Section~\ref{s:bgf}, introduce the Weyl structure
\begin{gather*}
g= \eta_1^2+\eta_2^2+\eta_3^2,
\qquad
\omega= 2\frac{\lambda_1\nu_1\eta_1+\lambda_2\nu_2\eta_2+\lambda_3\nu_3\eta_3} {\nu_1^2+\nu_2^2+\nu_3^2},
\end{gather*}
{where}
\begin{gather*}
\eta_i= \nu_i {\rm d} r + \phi^i {\rm d} q - \psi^i {\rm d} p,
\\
\nu_1= \phi^2\psi^3-\phi^3\psi^2,
\qquad
\nu_2= \phi^3\psi^1-\phi^1\psi^3,
\qquad
\nu_3= \phi^1\psi^2-\phi^2\psi^1.
\end{gather*}
The metric $g$ can be simplif\/ied considerably using the fact that area preserving vector f\/ields are locally hamiltonian:
this means one can write $\phi^i={F^i}_q$ and $\psi^i=-{F^i}_p$ for some functions $F^1$, $F^2$, $F^3$ of $(p,q,r)$.
The ($r$-dependent) constants of integration for the hamiltonians may be chosen so that the $F^i$ satisfy the
generalized Nahm equation in the Lie algebra of functions under Poisson bracket: ${F^1}_r-\{F^2,F^3\}=\sum\limits_j
{\mathcal{B}}_{1j} F^j$ and so on.
(The system ${p^i}_r-\sum\limits_j {\mathcal{B}}_{ij} p^i=q^i$ has local solutions for any $q^i(r)$, and these can be
added to $F^i$.
However, if the $\Phi^i$ take values in a~Lie subalgebra of the area preserving vector f\/ields, then the $F^i$ may take
values in a~central extension of this Lie subalgebra.) The metric now reduces to
\begin{gather*}
\nonumber g=\sum\limits_i \left({\rm d} F^i-\sum\limits_{j}{\mathcal{B}}_{ij} F^j {\rm d} r\right)^2
=\sum\limits_i \theta^i_{\alpha}(r) \theta^i_\beta(r) {\rm d} x^\alpha {\rm d} x^\beta,
\end{gather*}
where $\theta^i_\alpha$ is a~basis of solutions to the linear system ${\theta^i}_r=\sum\limits_j {\mathcal{B}}_{ij}
\theta^j$ and, using the summation convention on Greek indices, $F^i=\theta^i_\alpha x^\alpha$, for some functions
$x^\alpha(p,q,r)$ ($\alpha=1,2,3$).
The equation ${F^i}_r-\frac12\sum\limits_{j,k}\varepsilon_{ijk}\{F^j,F^k\} =\sum\limits_j {\mathcal{B}}_{ij} F^j$ is
immediately equivalent to
\begin{gather}
\label{eq:gnH}
\theta^i_{\alpha} {x^\alpha}_r =\sum\limits_{j,k} \varepsilon_{ijk} \theta^j_\alpha \theta^k_\beta\big\{x^\alpha,x^\beta\big\}.
\end{gather}
Now for \emph{any} functions $x^\alpha(p,q,r)$ direct calculations yield the general jacobian identities
\begin{gather*}
\{x^2,x^3\}{\rm d} x^1+\{x^3,x^1\}{\rm d} x^2+\{x^1,x^2\}{\rm d} x^3=J(x){\rm d} r,
\\
{x^1}_r{\rm d} x^2\wedge {\rm d} x^3+{x^2}_r{\rm d} x^3\wedge {\rm d} x^1
+{x^3}_r{\rm d} x^1\wedge {\rm d} x^2 =J(x){\rm d} p\wedge{\rm d} q,
\end{gather*}
where $J(x)$ is the determinant of the jacobian of $(x^1,x^2,x^3)$ with respect to $(p,q,r)$.
Hence~\eqref{eq:gnH} holds if and only if $*_g {\rm d} r={\rm d} p\wedge {\rm d} q$.

This argument was carried out using an af\/f\/ine coordinate $r$.
However, any coordinate $t$ can be used, simply by writing $\partial_t=\partial_r+a$.
The following result summarizes the general construction.

\begin{prop}
Consider the $3$-dimensional metric
\begin{gather*}
g=\sum\limits_{i} \bigl(\theta^i_{\alpha}(t){\rm d} x^\alpha\bigr)^2,
\end{gather*}
where $\theta^i_\alpha$ is a~basis of solutions of the linear system ${\theta^i}_t-a \theta^i =\sum\limits_j
{\mathcal{B}}_{ij} \theta^j$ for some symmetric traceless matrix ${\mathcal{B}}_{ij}(t)$, and $x^\alpha(p,q,t)$ are
arbitrary functions.
Then $*_g {\rm d} t={\rm d} p\wedge{\rm d} q$ if and only if the functions $F^i=\theta^i_\alpha
x^\alpha$ satisfy
\begin{gather*}
F^i_t-aF^i-\frac12\sum\limits_{j,k}\varepsilon_{ijk}\big\{F^j,F^k\big\}=\sum\limits_j {\mathcal{B}}_{ij} F^j.
\end{gather*}
\end{prop}

\begin{rem}
We pause brief\/ly to reiterate (conversely) the link between this analysis and the work of Ferapontov et
al.~\cite{FHZ:cqim}.
The framing $\theta^i_\alpha$ of the metric $g$ appearing in ${\rm d} *_g{\rm d} u =0$, i.e., the
formulation~\eqref{eq:ibg} of~\eqref{eq:fhz1}, provides a~diagonalization
\begin{gather*}
g_{\alpha\beta} {\rm d} x^\alpha {\rm d} x^\beta= \sum\limits_{i=1}^3 \bigl(\theta^i_{\alpha}
{\rm d} x^\alpha\bigr)^2 =\sum\limits_{i=1}^3 \bigl({\rm d} F^i - ({\mathcal{C}} F)^i {\rm d} u)^2,
\end{gather*}
where $F^i=\theta^i_{\alpha} x^\alpha$, $({\mathcal{C}} F)^i=\sum\limits_{i=1}^3 {\mathcal{C}}_{ij}F^j$.
Writing $*_g {\rm d} u = {\rm d} p \wedge {\rm d} q$ yields
\begin{gather*}
(F^i)'-\frac12 \sum\limits_{i,j,k} \varepsilon_{ijk}\big\{F^j,F^k\big\} =\sum\limits_j {\mathcal{C}}_{ij} F^j,
\end{gather*}
which reduces to the generalized Nahm equation after splitting ${\mathcal{C}}$ into its tracefree and tracelike parts.
\end{rem}

Returning to the $\SU(\infty)$ Toda equation, $(\theta^i_\alpha)$ is now a~diagonal matrix with eigenvalues
$(e^{u/2},e^{u/2},1)$.
The case that $u$ is constant (so that ${\mathcal{B}}=0$) is the overlap with the previous subsection, and this
construction recovers Ward's linearization of the $\SU(\infty)$ Nahm equation~\cite{Ward:sun}.
On the other hand if~$u$ is not constant, $t=u$ can be used as a~coordinate, so that~$a=1/3$ and~${\mathcal{B}}_{ij}$ is
diagonal with constant eigenvalues $1/6$, $1/6$, $-1/3$.
This is the type~(D) solution of the matrix Riccati equation with $\lambda=1/6$ (see Remark~\ref{r:Dproj}).

There is another class of Einstein--Weyl spaces, governed by the dKP equation~\cite{DMT:dKP}.
The Weyl structure in this case is Lorentzian:
\begin{gather*}
h={\rm d} y^2-4{\rm d} x {\rm d} t-4u {\rm d} t^2,
\qquad
\omega_h=2u_x{\rm d} t.
\end{gather*}
The Einstein--Weyl condition may again be written ${\rm d}{*_h{\rm d} u}=0$, which now reduces to the dKP
equation $u_{yy}=(u_t-uu_x)_x$.
If the solutions $(\theta^i_\alpha)$ are taken to be $(1+u,2{\boldsymbol i} u,0)$, $(2{\boldsymbol i} u,1-u,0)$ and
$(0,0,1)$, then $h= \sum\limits_{i,\alpha,\beta} \theta^i_{\alpha}(r)\theta^i_\beta(r) {\rm d} x^\alpha
{\rm d} x^\beta$ with $-\sqrt2t=x^1+{\boldsymbol i} x^2$, $2\sqrt2 x=x^1-{\boldsymbol i} x^2$, $y=x^3$, whereas
${\mathcal{B}}$ is the type (N) solution of the Riccati equation.
The equivalence works as long as $|{\rm d} u|_h\neq0$.

\begin{thm}
The $\SU(\infty)$ Toda and dKP equations are generically locally equivalent to $\SDiff(\Sigma^2)$ generalized Nahm
equations, on the type {\rm (D)} and {\rm (N)} Riccati space respectively.
\end{thm}

The equivalence in each case is obtained by a~``hodograph transformation'', i.e., dependent and independent variables
are exchanged.
The beauty of such transformations is that they interchange coordinate-invariance and gauge-invariance, which is why
they are potentially useful for studying geometric dif\/ferential equations, where the independent variables are often not
well def\/ined.

This result was f\/irst obtained in for the $\SU(\infty)$ Toda equation and type (D) Nahm equation.
The extension to more general Nahm equations resulted from discussions with Maciej Dunajski.
Together with Paul Tod, he has given a~detailed discussion of the dKP case~\cite{DuTo:p12}.

Solutions of the $\SU(\infty)$ Toda equation may be divided into two classes according to whether they are most easily
presented explicitly or implicitly: in the former class, there are very few examples~-- to the best of my knowledge, the
only examples are the separable solutions~\cite{BoFi:kve,GeDa:sst} and the hyperCR-Toda solutions~\cite{CaTo:emh} (the
latter were discussed in Section~\ref{s:s1Hit}).

The class of implicit solutions is much richer and all of them (that I know of) are obtained by the above hodograph
transformation.
There are two general classes.

\subsection*{Solutions with a~Killing vector~\cite{CaPe:sdc, Ward:sut}}

Suppose that the generalized Nahm f\/ield is invariant under a~one-dimensional symmetry group of~$\Sigma^2$, generated by
a~hamiltonian vector f\/ield $X$, with momentum map $\eta$, and let $\psi$ be a~function on~$\Sigma^2$ with
${\rm d}\psi(X)=1$.
Then $\psi,\eta$ can be used as coordinates on $\Sigma^2$, $X=\partial_\psi$ and the area form is
${\rm d}\psi\wedge{\rm d}\eta$.
The Lie algebra of hamiltonian vector f\/ields commuting with $\partial_\psi$ consists of vector f\/ields of the form
$b\partial_\eta+W(\eta)\partial_\psi$, where $b$ is constant (on $\Sigma^2$), and the hamiltonian is $-b\psi+V(\eta)$
where $V_\eta=W$.
The solutions found by Ward~\cite{Ward:sut} are given by the metric
\begin{gather*}
h=\rho^2\bigl({\rm d}(U_\eta)^2+{\rm d}\psi^2\bigr)+{\rm d}(\rho U_\rho)^2,
\end{gather*}
where $U(\rho,\eta)$ satisf\/ies $(\rho U_\rho)_\rho+\rho U_{\eta\eta}=0$, the equation for axisymmetric harmonic
functions on~${\mathbb{R}}^3$.
In this case one of the hamiltonians is simply~$\psi$.
The general solutions were found in~\cite{CaPe:sdc} following similar ideas.

\subsection*{Solutions constant on central ellipsoids~\cite{Tod:p3}}

Motivated by the Pedersen--Poon Ansatz for scalar-f\/lat K\"ahler  metrics~\cite{PePo:kzs}, Tod~\cite{Tod:p3} observed
that solutions of the Toda equation which are constant on central ellipsoids may be obtained from an ordinary
dif\/ferential equation, namely Painlev{\'e}'s third equation.
Painlev{\'e} III is equivalent to the isomonodromic deformation problem for an $\SU(2)$ connection with two double
poles, and hence to the generalized Nahm equation with gauge group $\SU(2)$ on the type (D) Riccati space.
Tod's solutions are obtained by reducing to the gauge algebra $\lie{so}(3)$ of inf\/initesimal isometries of $S^2$.
Solutions constant on central ellipsoids arise in this way because ellipsoids are precisely the quotients of left
invariant metrics on $S^3=\SU(2)$ by a~$\Un(1)$ subgroup.

Similar reductions are known in the dKP case~-- see~\cite{DMT:dKP,DuTo:p12}~-- and the latter construction is the context
for the central quadric Ansatz and Painlev\'e equations of~\cite{FHZ:cqim}.

\section{Further speculations}\label{section12}

The integrable background geometries discussed in this paper are those arising from the selfduality equation for
conformal structures in four dimensions.
Roughly speaking, they are characte\-ri\-zed as geometries associated to a~twistor space containing rational curves with
degree $2$ normal bundle, although this requires some interpretation for the one or two-dimensional geometries, where
the twistor space is not Hausdorf\/f and has dimension zero or one respectively.

Null reductions lead to twistor spaces with lower degree normal bundles.
For instance, a~surface with a~projective structure may be viewed as a~moduli space of rational curves in a~complex
surface with normal bundle ${\mathcal{O}}(1)$ (see Hitchin~\cite{Hit:gse} and LeBrun~\cite{LeB:phd}).

These considerations lead to the idea of extending the integrable background geometry concept to higher degree normal
bundles.
This would encompass the Einstein--Weyl and selfdual hierarchies (where the normal bundle is ${\mathcal{O}}(n)$ or
${\mathcal{O}}(n)\otimes{\mathbb{C}}^2$ respectively), as well as quaternionic geometry (normal bundle
${\mathcal{O}}(1)\otimes{\mathbb{C}}^{2k}$).

In the language of integrable systems, the concept of a~background geometry appears to be related to Lax systems
involving a~derivative with respect to the spectral parameter $\zeta$.
The zeros of the coef\/f\/icient of $\partial_\zeta$ are associated with poles of the Lax system: the prototype here is the
one-dimensional case (isomonodromic deformations).
In its most general form, the background geometry idea might be regarded as a~study of these extra $\partial_\zeta$
terms: general principles for the introduction of such terms have been developed by
Burtsev--Zakharov--Mikhailov~\cite{BMZ:isv}.
At the very least, it is desirable to know when these terms can be eliminated.
For the geometries studied in this paper, this question has an answer: they can be eliminated provided that the
background geometry is a~trivial Riccati space, a~trivial or spherical spinor-vortex space, a~hyperCR Einstein--Weyl
space, or a~hypercomplex selfdual space.

Returning to the general context of this paper, recent work of Ferapontov and Kruglikov~\cite{FeKr:disew} suggest deep
links between selfduality (hence twistor theory) and integrability in low dimensions, going beyond the inspirational
analysis by Mason and Woodhouse.
For example, in the particular case of equations of the form~\eqref{eq:ibg}, the Einstein--Weyl structure may be
recovered from the linearization ${\rm d}{*_{g(u)} {\rm d} v}+{\rm d}(v {*_g}'(u){\rm d}
u)=0$.
The leading term in $v$ is
\begin{gather*}
g^{\alpha\beta}(u)\frac{\partial^2 v}{\partial x^\alpha \partial x^\beta},
\end{gather*}
so the symbol of the linearized equation is the inverse metric $g^{-1}$ (up to a~conformal factor) of the Einstein--Weyl
structure.
Ferapontov and Kruglikov also show how to recover (the $1$-form of) the Weyl structure.
They obtain a~similar result for several classes of PDEs, relating the integrability of such equations (by the method of
hydrodynamic reductions) to Einstein--Weyl or selfdual conformal structures on their moduli spaces.

\appendix
\section{Table of integrable background geometries}

The following table summarizes the background equations and gauge f\/ield equations in each dimension (given in Euclidean
signature for convenience).

In this table, I indicate the data def\/ining the background geometry, starting with the conformal manifold $Q$ and the
conformal vector bundle ${\mathcal{V}}=L,{\mathcal{W}}$, or ${\mathcal{E}}$ (which has rank $1,2$, or $3$ respectively)
followed by the other data: $D$ is a~conformal connection on ${\mathcal{V}}$, $\overline\partial{}^a$ a~holomorphic
structure, whereas ${\mathcal{C}}$, $\psi$, ${\mathcal{B}}$ are sections of ${\mathcal{W}}^{*}\otimes_{\mathbb{C}}
{\mathcal{W}}^{*}\otimes_{\mathbb{C}} TN$, ${\mathcal{W}}^{*}\otimes_{\mathbb{C}} TN$ and $T^{*}
C\otimes\Sym_0{\mathcal{E}}$ respectively.

The gauge f\/ield equations are for a~pair $(A,\Phi)$, where $A$ is a~$G$-connection and the Higgs f\/ield $\Phi$ is an
$\lie{g}_Q$-valued section of ${\mathcal{V}}^{*}$, except that in four dimensions there is no Higgs f\/ield, and in one
dimension, the connection $A$ can be assumed trivial.

\medskip

\centerline{\begin{tabular}
{|l|c|c|c|}
\hline
Dimension & Data & Background equations & Gauge f\/ield equations
\\
\hline
4 Selfdual & $M$ & $\strut^{\strut} W^\asd = 0$ & $*F^A = F^A$
\\[1mm]
\hline
3 Einstein--Weyl & $(B,L),D$ & $\strut^{\strut} r^{D}_0=0$ & $*D^A\Phi=F^A$
\\[1mm]
\hline
2 Spinor-vortex &
$(N,{\mathcal{W}})$,
&
$\overline\partial{}^{\strut a}{\mathcal{C}}=0$,
&
$\strut F^A-[\Phi,\overline\Phi]_{\lie{g}} =\psi\wedge\overline\Phi+\overline\psi\wedge\Phi$,
\\[1mm]
&
$(\overline\partial{}^a,{\mathcal{C}},\psi)$
&
$\overline\partial{}^a\psi=-3{\mathcal{C}}\overline\psi$,
&
$\overline\partial{}^{a,A}\Phi={\mathcal{C}}\overline\Phi$
\\[1mm]
&
&
$s_{{\mathcal{W}}^{-1}  TN}=\psi\overline\psi-2{\mathcal{C}}\overline{\mathcal{C}}$
&
\\[1mm]
\hline
1 Riccati & $(C,{\mathcal{E}}),(D,{\mathcal{B}})$ & $\strut^{\strut} D{\mathcal{B}}=2\big(\mathcal{B}^2\big)_0$.
& $D^A\Phi-{*[\Phi,\Phi]_{\lie{g}}} ={\mathcal{B}}\mathinner{\cdot}\Phi$
\\[1mm]
\hline
\end{tabular}}

\medskip

Finally, I summarize the background geometry constructions in an $\SDiff(\Sigma)$-gauge, i.e., using an arbitrarily
chosen volume form $\nu$ on $\Sigma$ ($\nu={\rm d} t$ on $\Sigma^1$).

\subsection*{Selfdual spaces}

\begin{itemize}\itemsep=0pt
\item
From Einstein--Weyl spaces:
\begin{gather*}
g=\Phi^{2}\mathsf{c}_B + ({\rm d} t+A)^2.
\end{gather*}
\item
From spinor-vortex spaces (with conformal coordinate $z$):
\begin{gather*}
g=\big|\nu(\Phi,\overline\Phi)\big|^2 {\rm d} z {\rm d}\bar z+\ip{\nu(\Phi,\cdot)+\nu(\Phi,A),
\nu(\overline\Phi,\cdot)+\nu(\overline\Phi, A)}.
\end{gather*}
\item
From Riccati spaces (with af\/f\/ine coordinate $r$):
\begin{gather*}
g=\big|\nu(\Phi\mathbin{{\times}}\Phi\mathbin{{\times}}\Phi)\big|^2 {\rm d} r^2
+\bigl|\nu(\Phi\mathbin{{\times}}\Phi,\cdot)\bigr|^2.
\end{gather*}
\end{itemize}

\subsection*{Einstein--Weyl spaces}
\begin{itemize}\itemsep=0pt
\item
From spinor-vortex spaces (with conformal coordinate $z$):
\begin{gather*}
g= 4\Phi\overline\Phi {\rm d} z{\rm d}\bar z +({\rm d} t + \alpha {\rm d}
z+\overline\alpha {\rm d}\bar z)^2,
\\
\omega=\dot A - \frac{2\overline{\mathcal{C}}\Phi}{\overline\Phi}{\rm d} z
-\frac{2{\mathcal{C}}\overline\Phi}{\Phi}{\rm d}\bar z -\frac12\biggl(\frac{\psi+\dot\Phi}{\Phi}+
\frac{\overline\psi+\dot{\overline\Phi}}{\overline\Phi}\biggr) ({\rm d} t + \alpha {\rm d}
z+\overline\alpha {\rm d}\bar z).
\end{gather*}
\item
From Riccati spaces (with af\/f\/ine coordinate $r$):
\begin{gather*}
g=\big|\nu(\Phi\mathbin{{\times}}\Phi){\rm d} r+\nu(\Phi,\cdot)\big|^2,
\\
\omega=\frac{\Ip{2{\mathcal{B}}\bigl(\nu(\Phi\mathbin{{\times}}\Phi)\bigr)
-\nu(\Phi\mathbin{{\times}}\Phi)\mathbin{{\times}}\divg_\nu\Phi, \nu(\Phi\mathbin{{\times}}\Phi){\rm d}
r+\nu(\Phi,\cdot)}}{|\nu(\Phi\mathbin{{\times}}\Phi)|^2}.
\end{gather*}
\end{itemize}

\subsection*{Spinor-vortex spaces from Riccati spaces}

${\mathcal{W}}$ is the kernel of $\Phi$ in the pullback of ${\mathcal{E}}$ to $C\mathbin{{\times}}\Sigma^1$ with complex
structure given by cross product with $\Phi/|\Phi|$ and holomorphic structure induced by the conformal connection
\begin{gather*}
D+\frac{(\Phi_r\mathinner{\vartriangle}\Phi+\ip{{\mathcal{B}}\Phi,\Phi}\iden){\rm d} r
+\dot\Phi\mathinner{\vartriangle}\Phi {\rm d} t}{|\Phi|^2},
\end{gather*}
the conformal structure is represented by $g=|\Phi|^2{\rm d} r^2+{\rm d} t^2$ and the other two f\/ields are
\begin{gather*}
{\mathcal{C}}={\mathcal{B}}-\frac{{\mathcal{B}}\Phi\otimes\Phi+\Phi\otimes{\mathcal{B}}\Phi}{|\Phi|^2}
+\frac{\ip{{\mathcal{B}}\Phi,\Phi}}{2|\Phi|^2}\left(\iden+\frac{\Phi\otimes\Phi} {|\Phi|^2}\right),
\\
\psi=\frac{\Phi\mathbin{{\times}}(2{\mathcal{B}}\Phi+\Phi\mathbin{{\times}}\dot\Phi)}{|\Phi|^2}.
\end{gather*}

\subsection*{Acknowledgements}

I would like to thank Harry Braden, Maciej Dunajski, Evgeny Ferapontov, Paul Gauduchon, Nigel Hitchin, Boris Kruglikov,
Claude LeBrun, Lionel Mason, Ian Strachan, Paul Tod and Nick Woodhouse for helpful discussions.
This work was supported by the EPSRC, the Leverhulme Trust, and the William Gordon Seggie Brown Trust.

\pdfbookmark[1]{References}{ref}
\LastPageEnding

\end{document}